\documentclass[12pt]{article}
\newcommand{\R}{\mathbb{R}}
\newcommand{\N}{\mathbb{N}}

\newcommand{\Om}{\Omega}
\newcommand{\bfv}{\mathbf{v}}
\newcommand{\bfJ}{\mathbf{J}}
\newcommand{\dt}{\partial_t}

\newcommand{\dx}{\,\textrm{d}x}
\newcommand{\red}[1]{\textcolor{black}{#1}} % change to red for marked-up version

\newcommand{\trho}{\widetilde{\rho}}

\newcommand{\bfT}{\mathbf{T}}

\newcommand{\bff}{\mathbf{f}}
\newcommand{\bfh}{\mathbf{h}}
\newcommand{\bfb}{\mathbf{b}}
\newcommand{\bfm}{\mathbf{m}}
\newcommand{\bfn}{\mathbf{n}}

\newcommand{\bfD}{\mathbf{D}}
\newcommand{\vprel}{\bfv_{\mathrm{p,rel}}}

\newcommand\restr[2]{{% we make the whole thing an ordinary symbol
  \left.\kern-\nulldelimiterspace % automatically resize the bar with \right
  #1 % the function
  \littletaller % pretend it's a little taller at normal size
  \right|_{#2} % this is the delimiter
  }}

\newcommand{\littletaller}{\mathchoice{\vphantom{\big|}}{}{}{}}

\newcommand{\yes}{\checkmark}
\newcommand{\no}{--}

\usepackage[utf8]{inputenc}
\usepackage{graphicx}
\usepackage{mathtools}
\usepackage{amssymb}
\usepackage{amsthm} %for proof environment
\usepackage{mathabx} %for widebar
\usepackage[sc]{mathpazo}
\usepackage[table]{xcolor} %for \textcolor
\usepackage{tabularx}
\usepackage[colorlinks=true,linkcolor=blue,citecolor=blue]{hyperref}
\usepackage{appendix}
\usepackage{lipsum}
\usepackage{rotating}
\usepackage{adjustbox}
\usepackage{nicematrix, booktabs, siunitx}
\usepackage[labelfont=bf]{caption}
\usepackage{subcaption}
\usepackage{afterpage}
\usepackage{float}
\usepackage{verbatim}
\usepackage{fancyhdr} 
\usepackage[nottoc]{tocbibind}
\usepackage{pdflscape}
\usepackage{rotating}
\usepackage{placeins}
\usepackage{enumitem}
\usepackage{algorithm}
\usepackage{algpseudocode}
\usepackage{pgfplots}
\usepackage{booktabs}
\usepackage{pgfplotstable}
\usepackage{makecell}
\usepackage{rotating}

\definecolor{Dark2_1}{HTML}{1b9e77} % dark teal
\definecolor{Dark2_2}{HTML}{d95f02} % strong orange
\definecolor{Dark2_3}{HTML}{7570b3} % deep purple

\definecolor{mix1}{HTML}{e41a1c} 
\definecolor{mix2}{HTML}{4daf4a} 
\definecolor{mix3}{HTML}{377eb8} 
\definecolor{mix4}{HTML}{984ea3} 

\providecommand{\keywords}[1]{\textbf{Keywords:} #1}

\graphicspath{{./figures/}}

\allowdisplaybreaks
\newtheorem{theorem}{Theorem}[section]

\newtheorem{remark}[theorem]{Remark}

\usepackage[skins,breakable]{tcolorbox}

\newtcolorbox{myframe}[2][]{%
  enhanced,colback=white,colframe=black,coltitle=black,
  sharp corners,boxrule=1pt,
  fonttitle=\bfseries,
  attach boxed title to top left={yshift=-0.3\baselineskip-0.4pt,xshift=2mm},
  boxed title style={tile,size=minimal,left=0.5mm,right=0.5mm,
    colback=white,before upper=\strut},
  title=#2,#1
}

\usepackage{PRIMEarxiv}
\usepackage{booktabs}       
\usepackage{nicefrac}      
\usepackage{microtype}   

\author{
  Eberhard Bänsch,\; Jonas Knoch\thanks{Corresponding Author},\; Nicolas Neuß,\; Maria Neuss-Radu \\
  Department Mathematik \\
  Friedrich-Alexander-Universität Erlangen-Nürnberg \\
  Cauerstraße 11, 91058 Erlangen\\
  \texttt{baensch@math.fau.de, jonas.knoch@fau.de}, \\ \texttt{neuss@math.fau.de, maria.neuss-radu@math.fau.de}}

\title{Derivation and Numerical Simulation of a Thermodynamically Consistent Magneto Two-phase Flow Model}

\begin{document}
\maketitle

%% Abstract
\begin{abstract}
\red{
Motivated by Magnetic Drug Targeting (MDT), we derive a thermodynamically consistent two-phase model for the coupled dynamics of superparamagnetic iron oxide nanoparticles (SPIONs), a carrier fluid, and the magnetic field. The model includes a convection-diffusion equation for the SPION dispersion volume density and a modified \textit{Navier-Stokes} system for the mixture velocity with SPION-dependent density and viscosity, an additional momentum flux due to relative SPION motion, and a magnetic body force. Additionally, the SPIONs interact with a quasi-static \textit{Maxwell} system through their magnetization which is of \textit{Langevin}-type. The coupled system satisfies a global energy-dissipation balance. Moreover, we develop a semi-implicit finite element scheme with explicit flux-corrected transport stabilization and simulate the attraction of an initially localized SPION dispersion towards an external magnet in a system closed to mass exchange, using parameter values representative of MDT settings. The resulting SPION redistribution and induced flow are analyzed through physical observables, including energy and dissipation contributions and iron mass within a target region. Temporal convergence studies assess the ability of the numerical scheme to reproduce the energy-dissipation balance.
Finally, we study the influence of individual coupling mechanisms on particle accumulation, maximum iron density, and peak flow velocity.}
\end{abstract}

\keywords{Magnetohydrodynamics \and Superparamagnetism  \and  Nanofluid  \and \textit{Onsager}'s Variational Principle  \and Finite Element Simulations \and Magnetic Drug Delivery}

%%%%%%%%%%%%%%%%%%%%%%%%%%%%%%%%%%%%%%%%%%%%%%%%%%%%%%%%%%%%%%%%%%%%%%%%%%%%%%%%
\section{Introduction}%%%%%%%%%%%%%%%%%%%%%%%%%%%%%%%%%%%%%%%%%%%%%%%%%%%%%%%%%%
%%%%%%%%%%%%%%%%%%%%%%%%%%%%%%%%%%%%%%%%%%%%%%%%%%%%%%%%%%%%%%%%%%%%%%%%%%%%%%%%

Magnetic Drug Targeting (MDT) is a technique that uses external magnetic field gradients to deliver drug-loaded magnetic (nano-)particles specifically to their destination in the body, increasing the local drug concentration while reducing systemic distribution \cite{sun2008,liu2019}. Such treatment is particularly attractive for chemotherapy of solid tumors, where often far less than 1\% of the amount of drug used actually ends up in the tumor, while the rest leads to collateral damage in healthy tissue \cite{reszka1997,wilhelm2016}. In addition to various animal studies, see e.g. \cite{huang2012,tietze2013,marie2015}, proof of practical feasibility was also provided in clinical trials on humans \cite{lubbe1996,lubbe2001}, thus underlining the potential of MDT in cancer therapy. However, due to the complex interplay between fluid dynamics, magnetic forces and particle transport in complex vascular environments, a quantitative understanding and direct experimental optimization of MDT remains both time-consuming and limited in scope. Mathematical modeling and simulation can play a crucial role in understanding, predicting, and optimizing these processes \cite{shapiro2015}. 
\par
Accordingly, various mathematical models have been developed in the recent years in order to describe the dynamics of magnetic (nano-)particles in  drug targeting settings. In this context, two main paradigms can be found in the literature: particle-based, \textit{Lagrangian} approaches, in which trajectories of a limited number of discrete particles in the carrier medium are calculated based on \textit{Newton}'s second law, and concentration-based, \textit{Eulerian} approaches, in which the dynamics of the particles are described by means of a partial differential equation of convection-diffusion type. \red{Neither approach is inherently superior; the choice primarily depends on the quantities of interest and the physical mechanisms to be resolved. \textit{Lagrangian} models are particularly well suited for studying individual particle trajectories, particle-specific properties, and capture or deposition events. \textit{Eulerian} models, in turn, provide a natural framework for describing collective transport and spatial concentration distributions, as well as for incorporating concentration-dependent multiphysics couplings between the particle dispersion, the carrier fluid, and the magnetic field.}
\par
First results following the \textit{Lagrangian} approach can be found in \cite{furlani2006} and \cite{furlani2007}, where analytical expressions for the trajectories of magnetic particles in a fluid flow under the influence of an infinite cylindrical magnet were derived and subjected to parametric analysis. Here, both the flow field of the carrier medium as well as the magnetic field were fully prescribed and did not interact with the particles. Besides analytical expressions for the trajectories in a simplified setting where the magnetic force profile was described by a step function, the authors of \cite{barnsley2017} also simulate trajectories for more complicated setups involving a linear \textit{Halbach} array and compare their results with experimental measurements. Further work in this direction deals with the simulation of particle trajectories in the case of endoscopic placement of magnetic stents, see e.g. \cite{aviles2007,aviles2008} for an investigation and experimental validation including a flow field computed from the stationary \textit{Navier}-\textit{Stokes} equations and \cite{lindemann2021} for instationary \textit{Navier}-\textit{Stokes} equations in a multi-branched vessel. In \cite{cregg2009} and \cite{cregg2010}, additional forces on the particles due to magnetic dipole-dipole interactions as well as hydrodynamic interactions were modeled and their effect on the capture efficiency was studied. The authors of \cite{calandrini2018} incorporate magnetic particle tracking for a prescribed magnetic field into a fluid-structure model to account for the interaction between blood flow and vessel in the MDT setting.\par
An early work in the field of concentration-based description of particles in MDT applications is \cite{grief2005}, where the evolution of the particle concentration is given by a convection-diffusion equation for  prescribed blood velocity profile and magnetic force and agglomeration effects on the vessel wall are modeled by a surface concentration. 
In \cite{nacev2011}, a similar setting (prescribed flow field and constant magnetic force) for a convection-diffusion model is considered, and after non-dimensionalization the authors identify three relevant parameter regimes (magnetic force dominated, flow dominated and boundary layer regime) on the basis of which they categorize experimental setups available in the literature. More advanced models include the response of the fluid to the nanoparticle dynamics by means of a concentration-dependent body force term in the \textit{Navier}-\textit{Stokes} system describing the flow, see e.g. \cite{li2008} for simulations of different flow parameters and magnetic properties of particles and \cite{gitter2012} for a similar model featuring saturation effects for the magnetization.
On the other hand, the effect of the magnetic nanoparticles on the magnetic field is considered in \cite{himmelsbach2018}, where first analytical results for radially symmetrical solutions to a system of an elliptic equation for magnetic potential and a convection-diffusion equation with prescribed flow, coupled via a linear magnetization relation, are obtained. The bidirectional interaction between the magnetic particles and the magnetic field is further considered in \cite{reinelt2023}, where the saturation of the particle magnetization is considered and simulations of a reduced model are validated with measurements from particle capture experiments. Interestingly, in \cite{thalmayer2024}, the aforementioned model is used to simulate a new experimental setup involving a \textit{Halbach} array and found to have superior capabilities to describe the obtained measurements in comparison to an analogous particle-based \textit{Lagrangian} model. \red{A model featuring the coupling of a carrier fluid, magnetic nanoparticles and the magnetic field was derived based on thermodynamic principles in \cite{grun2019}. There, the nanoparticle's influence on the flow is given in terms of the magnetic body force on the fluid, while the magnetization is modeled using an evolution equation.} \par 
These models do not fully take into account the interaction of the carrier fluid and the particles in terms of a mixture model and hence are only valid as long as the particle concentration remains small. However, this is a critical aspect, given that the particles tend to accumulate within the target region. Descriptions of the fluid-particle mixture within the framework of two-phase models can be found in the papers \cite{badfar2020} and \cite{boutopoulos2020}, where the spatio-temporal evolution of the nanoparticle volume fraction is computed for increasing (prescribed) magnetic forces and, in the latter case, for a two-dimensional bifurcated vessel setup. \red{More closely related to our model is \cite{weiss2021}, where a mixture model for superparamagnetic nanoparticles and a carrier fluid, combined with a quasi-static magnetic \textit{Maxwell} system using a linearized magnetization law, is derived based on \textit{Onsager}’s variational principle. However, to the best of our knowledge, there is still no mixture model that combines the interactions between the carrier fluid and the SPION dispersion with the magnetic field and is suitable for MDT applications.
}
\par
\red{
The goal of this paper is to derive a magneto two-phase flow model that addresses the specific aspects of MDT in greater detail. To this end, the SPION-containing phase is modeled as a dispersion, where the SPIONs are clusters of iron oxide cores surrounded by a surface coating. The physical interpretation of the SPION dispersion volume density within the mixture is based on two packing factors, accounting for the packing density of iron cores within one SPION, as well as the maximum packing density of multiple SPIONs in the dispersion. The mathematical model for the evolution of the mixture of SPION dispersion and carrier fluid is based on general conservation laws and thermodynamic principles (like \textit{Onsager}’s variational principle), using suitable energy and dissipation functionals. Besides the classical energy of the magnetic field, the energy functionals used in the derivation include two additional magnetic contributions: one from the external magnetic field and one from the field-induced magnetization of the SPION dispersion. They allow us to obtain a global energy-dissipation law for the system that is in accordance with the second law of thermodynamics. The resulting system consists of a convection-diffusion equation for the SPION dispersion volume density and a modified \textit{Navier-Stokes} system for the mixture velocity and pressure with SPION-dependent density and viscosity, an additional momentum flux due to relative SPION motion and a magnetic body force. Moreover, the model takes into account the interaction of the SPIONs with the magnetic field by means of a nonlinear magnetization of \textit{Langevin} type within a quasi-static magnetic \textit{Maxwell} system.}
\par
\red{
For the computation of physical observables relevant for MDT applications, we develop a numerical scheme based on the finite element method and semi-implicit time-stepping involving a one-step $\theta$-method for the temporal discretization. The convection-dominated dynamics of the nanoparticle phase is stabilized using an explicit flux-corrected transport method, leading to non-negative SPION dispersion volume densities. As a numerical experiment, we simulate the attraction of an initially localized SPION dispersion towards an external magnet in a system closed to mass exchange, using parameter values representative of MDT settings. We interpret the resulting redistribution of the SPION dispersion and the induced flow using a range of physical observables, including the individual energy and dissipation functionals and the fraction of the total SPION mass accumulated within a prescribed target region. Subsequently, we investigate the proposed scheme through temporal convergence studies and assess how accurately it reproduces the global energy-dissipation balance. Last but not least, we conduct an ablation study to investigate the relative importance of five individual coupling mechanisms in the model with respect to simulation outcomes, comparing $2^5 = 32$ model configurations.}
\par
\red{The paper is structured as follows: Section \ref{sec:modeling} is dedicated to the mathematical derivation and non-dimensionalization of the magneto two-phase system. Section \ref{sec:numerical_simulation} presents the numerical scheme and its application to the magnetic attraction of an initially localized SPION dispersion in a closed system, followed by an investigation of physical observables, temporal convergence and the discrete energy-dissipation balance, and a study of the importance of the model’s coupling mechanisms. The paper wraps up with a conclusion and an outlook in Section \ref{sec:conclusion_outlook}.
}

%%%%%%%%%%%%%%%%%%%%%%%%%%%%%%%%%%%%%%%%%%%%%%%%%%%%%%%%%%%%%%%%%%%%%%%%%%%%%%%%
\section{Derivation of the Thermodynamically Consistent Magneto Two-Phase Model}%%%%%%%%%%%%%%%%%%%%%%%%%%%%%%%%%%%%%%%%%%%%%%%%%
\label{sec:modeling}%%%%%%%%%%%%%%%%%%%%%%%%%%%%%%%%%%%%%%%%%%%%%%%%%%%%%%%%%%%%
%%%%%%%%%%%%%%%%%%%%%%%%%%%%%%%%%%%%%%%%%%%%%%%%%%%%%%%%%%%%%%%%%%%%%%%%%%%%%%%%

In this section, we derive a thermodynamically consistent magneto two-phase flow model for dispersed SPIONs in a carrier fluid based on general conservation laws and \textit{Onsager}'s variational principle \cite{onsager1931a,onsager1931b,doi2011}. \red{The so-obtained model fulfills a global energy-dissipation law in accordance with the second law of thermodynamics.}

\subsection{\red{Modeling of the SPION Dispersion}}

\begin{figure}
    \centering
    \begin{subfigure}[t]{0.33\textwidth}
        \centering
        \includegraphics[width=0.9\textwidth]{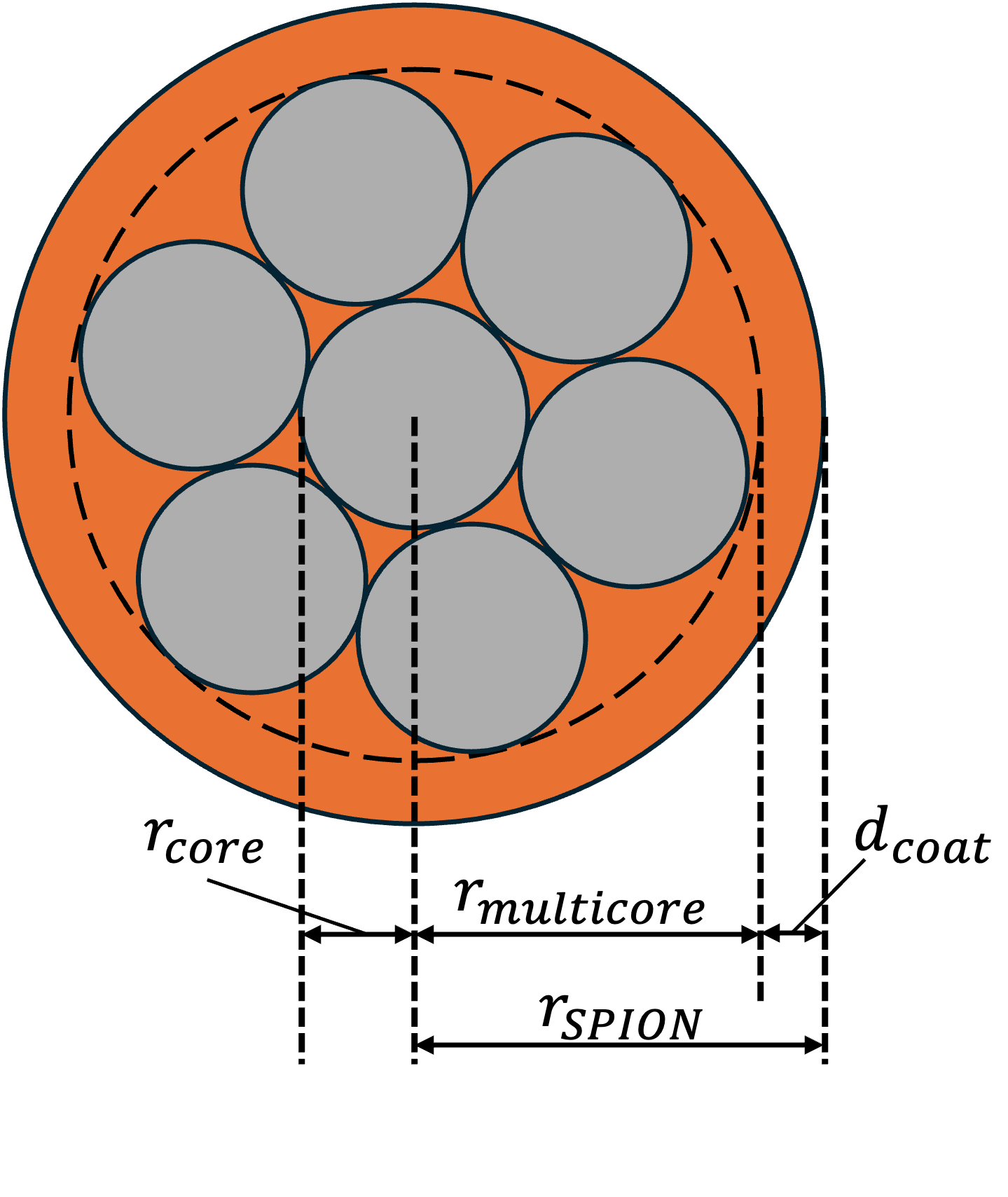}
        \caption{\red{A single SPION consisting of multiple iron oxide cores (gray), coated with a biocompatible layer (orance).}}
        \label{fig:spion_sketches:a}
    \end{subfigure}
    \begin{subfigure}[!t]{0.33\textwidth}
        \centering
        \includegraphics[width=0.8\textwidth]{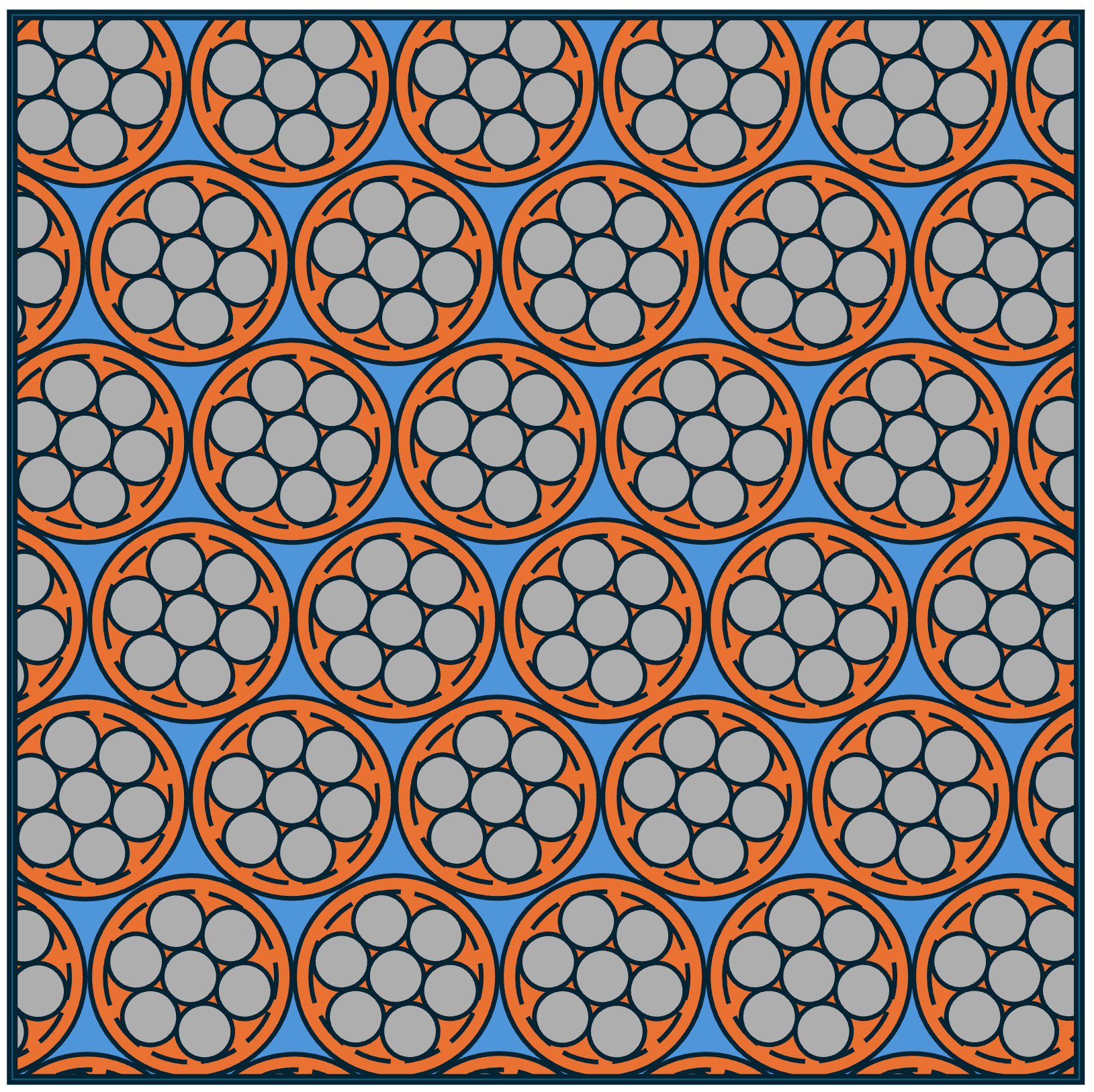}
        \caption{\red{A sketch of a densely packed SPION dispersion.}}
        \label{fig:spion_sketches:b}
    \end{subfigure}
    \begin{subfigure}{0.33\textwidth}
\begin{tikzpicture}
  \begin{axis}[
    axis lines=middle,
    xlabel={$|\bfh|$},
    ylabel={$\mathcal{L}(\xi |\bfh|)$},
    width=\textwidth,
    enlargelimits=true,
    xtick=\empty,
    ytick={0, 0.5, 1},
    legend style={at={(0.5,-0.15)}, anchor=north}
  ]
    \addplot [
  color=Dark2_1,
  domain=0:20,
  samples=400,
  ultra thick
] {abs(x) < 0.2 ? (x/3 - x^3/45) : (cosh(x)/sinh(x) - 1/x)};
    %\addlegendentry{$\mathcal{L}(x) = \coth(x) - \frac{1}{x}$}
    \addplot[black, thick, dotted, domain=0:2.5] {x/3};
    \draw[->, thick] 
  (axis cs:5,0.7) -- (axis cs:2.6,0.77);
\node[anchor=west] at (axis cs:5,0.7) {$\frac{\xi}{3} |\bfh|$};
    %\draw[->, thick] (axis cs:3,0.5) -- (axis cs:2.3,0.77) node[midway, above, sloped] {$\frac{\xi}{3} |\bfh|$};
  \end{axis}
\end{tikzpicture}
\caption{\red{Plot of the \textit{Langevin} function $\mathcal{L}(\xi |\bfh|) = \coth(\xi |\bfh|) - (\xi |\bfh|)^{-1}$.}}
    \label{fig:spion_sketches:c}
    \end{subfigure}
    \caption{\red{Schematic representation of a single SPION, a densely packed SPION dispersion, and a graphical representation of the \textit{Langevin} function and its linearization at the origin.}}
    \label{fig:spion_sketches}
\end{figure}

\red{A single SPION consists of several small iron oxide cores with radii of only a few nanometers. Together, these iron oxide cores form an approximately spherical multi-core cluster which is surrounded by a surface coating to enhance colloidal stability and biocompatibility, see Figure \ref{fig:spion_sketches:a}. Each iron core usually consists of magnetite (Fe$_3$O$_4$) and maghemite ($\gamma$-Fe$_2$O$_3$), which are ferrimagnetic materials with a high saturation magnetization. The volume fraction of the iron oxide cores inside a SPION is given by 
\begin{equation}\label{eq:varphi_cores}
\varphi_\text{cores} = P_\text{cores} \left(\frac{r_\text{multicore}}{r_\text{multicore} + d_\text{coat}}  \right)^3,
\end{equation}
where $P_\text{cores} \in [0,1]$ is the packing density of the iron oxide cores inside the SPION multi-core cluster, $r_\text{multicore}$ is the radius of the multi-core cluster and $d_\text{coat}$ is the thickness of the surface coating. Hence, we obtain the density of a SPION as 
\begin{equation} \label{eq:rho_spion}
\tilde{\rho}_\text{SPION} = \tilde{\rho}_\text{cores} \varphi_\text{cores} + \tilde{\rho}_\text{coat}(1-\varphi_\text{cores}),
\end{equation}
where $\tilde{\rho}_\text{cores}$, $\tilde{\rho}_\text{coat}>0$ are the densities of the iron oxide core compounds and the surface coating, respectively. In the MDT experiments, the SPIONs are present in a dispersion, i.e. they are surrounded by a fluid. There, the SPION's maximal packing density sets an upper limit on how much they can aggregate, see Figure \ref{fig:spion_sketches:b}. We therefore compute the density of a densely packed SPION dispersion as 
\begin{equation} \label{eq:rho_p}
\tilde{\rho}_\text{p} = P_\text{SPION} \tilde{\rho}_\text{SPION} + (1-P_\text{SPION}) \tilde{\rho}_\text{disp},
\end{equation}
with $P_\text{SPION} \in [0,1]$, the maximal packing density of SPIONs in the dispersion and $\tilde{\rho}_\text{disp}>0$, the density of the dispersion agent.
}

%%%%%%%%%%%%%%%%%%%%%%%%%%%%%%%%%%%%%%%%%%%%%%%%%%%%%%%%%%%%%%%%%%%%%%%%%%%%%%%%
\subsection{\red{Mass and Momentum Balances for a Mixture of SPION Dispersion and Carrier Fluid}}
\label{sec:mass_and_momentum_balance}
%%%%%%%%%%%%%%%%%%%%%%%%%%%%%%%%%%%%%%%%%%%%%%%%%%%%%%%%%%%%%%%%%%%%%%%%%%%%%%%%
\red{We want to model situations where a SPION dispersin is mixed with a carrier fluid (typically water or blood). To this end,} let us consider a mixture of two phases that occupy a simply connected open domain $\Om \subset \R^n$ with $n \in \{2,3\}$ for the physically relevant cases \red{and a time interval $I := (0,T_\text{end})$, $T_\text{end}>0$}. \red{The phase representing the SPION dispersion is marked by the index $\mathrm{p}$ and the phase representing the carrier fluid is marked with index $\mathrm{f}$, respectively.} We denote the specific constant density of the pure carrier fluid phase by $\trho_f>0$ and of the fully saturated SPION dispersion by $\trho_\mathrm{p}>0$ as introduced above. \red{These densities act as reference values for the pure phases.} The total mass density of the mixture is denoted by
\begin{equation}\label{density}
\rho = \rho(t,x) := \rho_\mathrm{f}(t,x) + \rho_\mathrm{p}(t,x),
\end{equation}
where $\rho_\mathrm{f}, \rho_\mathrm{p}$ are the time and space dependent densities of the phases within the mixture. Moreover, let us define the volume densities \red{$u_j \colon I \times \Omega \rightarrow [0,1]$} of the respective phases in the mixture by 
\begin{equation}
\label{def:volume_density}
u_j(t,x) := \frac{1}{\trho_j} \rho_j(t,x), \quad j \in \{\mathrm{f},\mathrm{p}\}.
\end{equation}
Hence, for example, 
\begin{equation}
    \red{V_\mathrm{p}(t_*) := \int_V u_\mathrm{p}(t_*,x) \rm{d}x}
\end{equation}
is the partial volume occupied \red{by the SPION dispersion} in a test volume $V \subset \Omega$ at time $t_*$. \red{In particular, given the volume $V_\text{SPION} = 4/3 \pi r_\text{SPION}^3$ of one SPION,
\begin{equation}
    N_\text{SPION}(t_*) = P_\text{SPION}\frac{V_\mathrm{p}(t_*)}{V_\text{SPION}}
\end{equation}
is the number of SPIONs present in volume $V$ at time $t_*$.}
Let us now introduce the \textit{volume-averaged velocity} of the mixture 
\begin{equation}\label{eq:volume_averaged_velocity}
\bfv := u_\mathrm{f} \bfv_\mathrm{f} + u_\mathrm{p} \bfv_\mathrm{p},
\end{equation}
where $\bfv_j$ denotes the individual velocity of phase $j\in \{\mathrm{f},\mathrm{p}\}$. From the mass balance equations for both phases:
\begin{equation}
\label{eq:individual_mass_balance_rho}
\dt \rho_j + \nabla \cdot (\rho_j \bfv_j) = 0, 
\end{equation}
we obtain by using definition \eqref{def:volume_density}
\begin{equation}
\label{eq:individual_mass_balance_u}
\dt u_j + \nabla \cdot (u_j \bfv_j) = 0.
\end{equation}
Upon assuming that the excess volume due to mixing of the two phases is zero, i.e.
\begin{equation}\label{eq:no_excess_volume}
u_\mathrm{f} + u_\mathrm{p} = 1,
\end{equation}
and adding equations \eqref{eq:individual_mass_balance_u}, we obtain that the volume averaged velocity $\bfv$ is divergence-free,
\begin{equation}\label{eq:v_divergence_free}
\nabla \cdot \bfv = 0.
\end{equation}
Introducing the relative velocity of phase $j$ with respect to the volume-averaged velocity
\begin{equation}
    \bfv_{j,\mathrm{rel}} := \bfv_j - \bfv
\end{equation}
yields for $j\in \{\mathrm{f},\mathrm{p}\}$
\begin{equation}
\label{eq:individual_mass_balance_with_averaged_velocity}
\dt u_j + \nabla \cdot (u_j \bfv +u_j \bfv_{j,\mathrm{rel}}) = 0.
\end{equation}
Note that \eqref{eq:no_excess_volume} renders one of the equations in \eqref{eq:individual_mass_balance_with_averaged_velocity} redundant and we rely in the following on the equation describing the SPION phase.
\par
According to \cite{abels2012}, the \textit{volume-averaged velocity} $\bfv$ fulfills the momentum equation
\begin{equation}
\label{eq:momentum_balance}
\rho \dt \bfv + (( \rho \bfv + \bfJ_\mathrm{p}) \cdot \nabla) \bfv - \nabla \cdot \bfT + \nabla p  = \bff,
\end{equation}
with stress tensor $\bfT$, pressure $p$, force density $\bff$ and $\rho = \rho(u_\textrm{p}) := (\tilde{\rho}_\mathrm{p} - \tilde{\rho}_\mathrm{f})u_\mathrm{p} + \tilde{\rho}_\mathrm{f}$. Here, the additional flux $\bfJ_\mathrm{p}$ is given as
\begin{equation}\label{eq:particle_flux}
\bfJ_\mathrm{p} = (\tilde{\rho}_\mathrm{p} - \tilde{\rho}_\mathrm{f}) u_\mathrm{p} \vprel.
\end{equation}
\red{The additional inertia term $(\bfJ_\mathrm{p} \cdot \nabla) \bfv$ accounts for momentum transport associated with mass fluxes relative to the volume-averaged mixture velocity. In our case, these relative mass fluxes will arise from diffusion and magnetophoretic drift of the SPIONs.}
The constitutive relations for $\bfT$ and $\vprel$, as well as the form of the volume force $\bff$, are derived in Section \ref{sec:derivation_constitutive_relations}.

%%%%%%%%%%%%%%%%%%%%%%%%%%%%%%%%%%%%%%%%%%%%%%%%%%%%%%%%%%%%%%%%%%%%%%%%%%%%%%%
\subsection{The Magnetic Field under the Influence of the SPION Dispersion}%%%%%%%%%%%%%%%%%%%%%%%%%%%%%%%%%%%%%%%%%%%%%%%%
\label{sec:magnetic_field}%%%%%%%%%%%%%%%%%%%%%%%%%%%%%%%%%%%%%%%%%%%%%%%%%%%%%
%%%%%%%%%%%%%%%%%%%%%%%%%%%%%%%%%%%%%%%%%%%%%%%%%%%%%%%%%%%%%%%%%%%%%%%%%%%%%%%
%
For the interaction of the SPIONs with the magnetic field, we consider \textit{Maxwell}'s equations for matter in the absence of free charges, free electric currents and time-dependent effects, i.e. 
\begin{subequations}\label{eq:maxwell}
    \begin{align}
    \nabla \cdot \bfb = 0, \quad  \nabla \times \bfh = 0, \quad \bfb &= \mu_0(\bfm + \bfh) & \text{ in } & I \times \Om,\label{eq:maxwell:a} \\
    \bfb \cdot \bfn &= \bfb_\text{ext} \cdot \bfn & \text{ on } & I \times \partial \Om, \label{eq:maxwell:b}
\end{align}
\end{subequations} 
where $\bfb, \bfh$ and $\bfm$ denote the magnetic flux density, the magnetic field strength and the magnetization, respectively. $\mu_0>0$ is the vacuum magnetic permeability and by $\bfn$ we denote the outer unit normal to $\Om$. $\bfb_\text{ext} \colon I \times \overline{\Om} \rightarrow \R^n$ is the external magnetic field that is used to target the SPIONs inside the domain. %
\par
The SPIONs have special magnetic properties, which are ideal for drug targeting applications.
\red{Due to their small size, the individual magnetic cores are single-domain particles, and thermal fluctuations can continuously change the direction of their magnetic moments, leading to a vanishing effective magnetization. In the presence of an external magnetic field however, the magnetic moments preferentially align with the field, resulting in a nonzero net magnetization. This property is called superparamagnetism and forms the foundation for the medical application of SPIONs as magnetically manipulable drug carriers \cite{vangijzegem2023}.} Accordingly, the magnetization \red{of the mixture} is given by \cite[Ch. 11.6]{cullity2011}
\begin{equation}\label{eq:magnetization}
    \bfm(u_\textrm{p},\bfh) = M_\textrm{s} u_\textrm{p} \mathcal{L}(\xi |\bfh|) \frac{\bfh}{|\bfh|},
\end{equation}
where 
\begin{equation}\label{eq:M_s}
\red{M_\textrm{s} = P_\text{SPION} \varphi_\text{cores} M_\text{Mag}}
\end{equation}
\red{is the saturation magnetization of the densely packed dispersion with $M_\text{Mag}>0$, the saturation magnetization of the iron oxide core compound magnetite/maghemite}. Moreover, $u_\textrm{p}$ is the volume density of the SPION phase, cf. Eq. \eqref{def:volume_density}, $\mathcal{L}(\alpha) = \coth(\alpha) - \frac{1}{\alpha}$ is the \textit{Langevin} function and 
\begin{equation}\label{eq:xi}
\xi = \frac{\mu_0 p}{k_\textrm{B} T}
\end{equation}
is the \textit{Langevin} parameter, with the magnetic dipole moment magnitude of a single SPION 
\begin{equation} \label{eq:p}
p = \varphi_\text{cores} M_\text{Mag} V_\text{SPION},
\end{equation}
\textit{Boltzmann}'s constant $k_\textrm{B}$ and temperature $T>0$. \red{Note that we use the continuous extension $\mathcal{L}(0)=0$ of the \textit{Langevin} function at the origin and define $\bfm(u_\textrm{p},0) = 0$}. For a graphical illustration of the \textit{Langevin} function, we refer to Figure \ref{fig:spion_sketches:c}.
\par
Since the magnetic field strength $\bfh$ is rotation free and $\Om$ is assumed to be simply connected, we deduce the existence of a magnetic potential $\phi \colon I \times \Om \rightarrow \R$ such that $\bfh = -\nabla \phi$. Using this relation and plugging the expression for the magnetization \eqref{eq:magnetization} into \eqref{eq:maxwell:a}-\eqref{eq:maxwell:b}, we obtain the nonlinear elliptic \textit{Neumann} problem
\begin{subequations}
    \begin{align}
        - \nabla \cdot \left( \mu_0 \left( 1 + M_\textrm{s} u_\textrm{p} \frac{\mathcal{L}(\xi |\nabla \phi|)}{|\nabla \phi|}\right) \nabla \phi \right) &= 0 & \text{ in } &I \times \Omega,\\
        - \mu_0 \left( 1 + M_\textrm{s} u_\textrm{p} \frac{\mathcal{L}(\xi |\nabla \phi|)}{|\nabla \phi|}\right) \nabla \phi \cdot \bfn &= \bfb_\text{ext} \cdot \bfn & \text{ in } &I \times \partial \Omega. \label{eq:magnetic_field_bc}
    \end{align}
\end{subequations}
\red{Here, the term $\mathcal{L}(\xi |\nabla \phi|)/|\nabla \phi|$ is understood by continuous extension at $|\nabla\phi|=0$, with limiting value $\xi/3$.} Since the problem depends on $\phi$ only through its gradient, the magnetic potential is only determined up to an additive constant. Uniqueness can be enforced, e.g., by demanding the normalization
\begin{equation}
    \int_\Omega \phi \dx = 0 \qquad \forall t \in I.
\end{equation}
Moreover, the prescribed magnetic field $\mathbf{b}_\text{ext}$ has to satisfy the compatibility condition
\begin{equation}\label{eq:b_ext_compatibility}
    \int_{\partial\Om} \bfb_{\mathrm{ext}}\cdot\bfn  \, \mathrm{d}\sigma = 0.
\end{equation}
By the divergence theorem, this condition is fulfilled if $\nabla \cdot \bfb_\text{ext} = 0$ in $\Omega$.
\begin{remark}
The literature is not consistent with regard to the modeling of SPION magnetization $\bfm$. For the qualitative behavior, the \textit{Langevin} relationship,  cf. \eqref{eq:magnetization}, prevails, see e.g. \cite{bean1959,cullity2011,grief2005,driscoll1984,reinelt2023}, but one also finds expressions for $\bfm$ that depend linearly \cite{rukshin2017,han2008,furlani2007} or linearly until saturation \cite{momeni2016,haverkort2009,calandrini2018} on the magnetic field. In \cite{grun2019}, an individual relaxation equation for the magnetization was considered, based on phenomenological arguments, see also \cite{shliomis2002}. Even more fundamental is the question of whether the magnetization should be a function of the magnetic flux density $\bfb$ or the magnetic field strength $\bfh$. Both approaches can be found in the literature and if the magnetic field is assumed not to interact with the SPIONs, they are equivalent. However, if one considers the influences of the SPIONs on the magnetic field, different equations for the magnetic variables result from the \textit{Maxwell} system \eqref{eq:maxwell:a}-\eqref{eq:maxwell:b}, depending on whether $\bfm$ is given in terms of  $\bfb$ or $\bfh$, see also \cite{reinelt2023}.
\end{remark}
\begin{remark}
    \red{
    The prescribed magnetic boundary data $\bfb_\mathrm{ext} \cdot \bfn $ represents the normal component of the externally applied magnetic flux density in the absence of SPIONs. By imposing this type of boundary condition, the contribution of the SPION-induced perturbation of the magnetic field on the boundary and the exterior of the domain is neglected. This problem formulation can be regarded as an approximation of the magnetostatic problem in the full space and is expected to be appropriate when the perturbation of the magnetic field due to the magnetic nanopartiles is weak. Our simulations indicate that said perturbation remains rather small ($\leq1.2\%$) for the considered SPION dispersion iron density of $1 \, \mathrm{mg /ml}$, see the first row in Figure 5. However, our investigations in Section \ref{sec:relative_importance_study} show evidence that the simulation outcome becomes more sensitive to the interaction of the SPIONs with the magnetic field for an iron density as high as $10 \, \mathrm{mg /ml}$, highlighting the need for a more accurate representation of the magnetostatic full space problem once the SPION dispersion becomes more concentrated.
    }
\end{remark}
%
%%%%%%%%%%%%%%%%%%%%%%%%%%%%%%%%%%%%%%%%%%%%%%%%%%%%%%%%%%%%%%%%%%%%%%%%%%%%%%%
\subsection{Derivation of Constitutive Relations Based on \textit{Onsager}'s Variational Principle}%%%%
\label{sec:derivation_constitutive_relations}%%%%%%%%%%%%%%%%%%%%%%%%%%%%%%%%%%
%%%%%%%%%%%%%%%%%%%%%%%%%%%%%%%%%%%%%%%%%%%%%%%%%%%%%%%%%%%%%%%%%%%%%%%%%%%%%%%
Next, we derive expressions for $\bfv_{\textrm{p,rel}}$ and $\bfT$ using \textit{Onsager}'s variational principle which can be regarded as an extension of \textit{Rayleigh}'s principle of the least energy dissipation to general irreversible processes. \red{We focus on the case where the system under consideration is closed, i.e. there is no mass exchange across $\partial \Om$:}
\begin{equation} \label{eq:velocity_SPIONs_bc}
\red{
    \mathbf{v} = 0, \qquad (u_\textrm{p}\mathbf{v} + u_\textrm{p} \vprel) \cdot \mathbf{n} = 0 \qquad \mbox{ on } I \times \partial \Omega.
}    
\end{equation}
\red{We note, however, that the system is influenced by the surrounding magnetic field through the boundary condition \eqref{eq:maxwell:b}.}
\par
Let us first state the energy and dissipation functionals that govern the physical system. We consider the kinetic energy of the two-phase system given by
\begin{equation}\label{def:kinetic_energy}
    \mathcal{E}_{\textrm{kin}}(t) := \frac{1}{2} \int_\Om \rho(u_\mathrm{p}) |\bfv|^2 \textrm{d}x.
\end{equation}
Moreover, let us consider the differential of the \textit{Gibbs} free energy density for an isothermal system \cite[Chapter 9.2.3]{ansermet2019}
\begin{equation}
\begin{aligned}
    \mathrm{d}g &= \frac{1}{V_\text{SPION}}\sum_{i \in \{\mathrm{p},\mathrm{f}\}} \mu(u_i')\mathrm{d}u_i' - \red{b} \mathrm{d}h \\
    &= \frac{1}{V_\text{SPION}}\sum_{i \in \{\mathrm{p},\mathrm{f}\}} \mu(u_i')\mathrm{d}u_i' \red{- \mu_0(h + m(h))}\textrm{d}h.
\end{aligned}
\end{equation}
Here, we used that $\bfh$ and $\bfm$ point in the same direction and hence the last equation in \eqref{eq:maxwell:a} also holds for the norms of the vectors. $\mu(u):= k_BT\ln(u)$ denotes the chemical potential and $m(h) := M_\mathrm{s} u_\mathrm{p} \mathcal{L}(\xi h)$ is the magnetization magnitude. Integration and insertion of $u_\mathrm{f} = 1-u_\mathrm{p}$ leads to
\begin{equation}
\begin{aligned}
g(t,x) &= \frac{k_\mathrm{B}T}{V_\text{SPOIN}}\sum_{i \in \{\mathrm{p},\mathrm{f}\}} \int_0^{u_i} \ln(u_i')\mathrm{d}u_i' - \mu_0 \left(  \int_0^{|\bfh|} \red{h} + M_\mathrm{s} u_\mathrm{p}\mathcal{L}(\xi h) \textrm{d}h \right)\\
&= \frac{k_\mathrm{B}T}{V_\text{SPION}} \Bigl[  u_\mathrm{p} ( \ln(u_\mathrm{p}) - 1) +  (1-u_\mathrm{p}) ( \ln(1-u_\mathrm{p}) - 1)\Bigr] \\
& \qquad \red{- \frac{\mu_0}{2}|\bfh|^2} -  \frac{\mu_0 M_\mathrm{s} u_\mathrm{p}}{\xi} \ln \left(\frac{\sinh(\xi |\bfh|)}{\xi |\bfh|} \right).
\end{aligned}
\end{equation}
From this we identify the mixing energy of the system as 
\begin{equation}\label{def:mixing_energy}
\mathcal{E}_\text{mix}(t) :=\frac{k_\mathrm{B}T}{V_\text{SPION}} \int_\Om \left(  u_\mathrm{p} ( \ln(u_\mathrm{p}) - 1) +  (1-u_\mathrm{p}) ( \ln(1-u_\mathrm{p}) - 1)\right) \textrm{d}x,
\end{equation}
as well as the magnetic energy contributions given by
\begin{equation}
    \red{\mathcal{E}_\text{mag,1}(t) := - \frac{\mu_0}{2} \int_\Om |\bfh|^2 \textrm{d}x}
\end{equation}
and
\begin{equation}\label{def:magnetic_energy}
\mathcal{E}_\text{mag,2}(t) := - \frac{\mu_0 M_\mathrm{s} }{\xi} \int_\Om u_\mathrm{p} \ln \left(\frac{\sinh(\xi |\bfh|)}{\xi |\bfh|} \right) \mathrm{d}x.
\end{equation}
\red{The apparent singularities in the energy contributions are removable. We use the continuous extensions $0\ln(0):=0$ and $\sinh(0)/0:=1$.}
\red{Additionally, we consider the contribution from the external magnetic field $\bfb_e$ in form of the energy
\begin{equation}
    \mathcal{E}_\text{mag,3}(t) := - \int_{\partial \Omega} \mathbf{b}_\text{ext} \cdot \mathbf{n} \phi\, \textrm{d}\sigma.
\end{equation}}
$\mathcal{E}_\text{mix}$ is the mixing energy of an ideal solution and \red{$\mathcal{E}_\text{mag,$i$}, i=1,2,3$, are the contributions from the magnetic field itself, the field-induced magnetization of the SPION dispersion, and the externally applied magnetic loading due to the external magnet, respectively}. Moreover, we consider the following dissipation functionals:
\begin{equation}
\label{eq:kinetic_dissipation}
\mathcal{D}_\text{kin}(t):= \int_\Om \frac{|\bfT|^2}{4\eta(u_\mathrm{p})} \dx
\end{equation}
and
\begin{equation}
\label{eq:stokes_dissipation}
\mathcal{D}_\text{drag}(t) := \frac{1}{2} \int_\Om \zeta(u_\mathrm{p}) |\vprel|^2 \dx,
\end{equation}
where we assume \red{that the fraction of the SPION dispersion within the mixture is small such} that $\eta$ is given in terms of the particle volume density by \textit{Einstein}'s equation for effective viscosity \cite{einstein1905}
\begin{equation}
    \eta(u_\mathrm{p}) = \eta_\mathrm{f}(1+2.5u_\mathrm{p}),
\end{equation}
with $\eta_\mathrm{f}>0$ denoting the viscosity of the pure fluid carrier phase. Here, $\mathcal{D}_\text{kin}$ and $\mathcal{D}_\text{drag}$ account for the energy losses due to viscous effects within the averaged fluid flow and due to drag associated with the relative movement of the particle phase with respect to the averaged velocity, respectively. In \eqref{eq:stokes_dissipation}, we assume that the friction coefficient is given by
\begin{equation}
\label{eq:friction_coefficient}
\zeta(u_\mathrm{p}) :=  \frac{6\pi \eta_\mathrm{f} r_\text{SPION}}{V_\text{SPION}}  \frac{u_\mathrm{p}}{1-u_\mathrm{p}}.
\end{equation}
Since by this relation the friction coefficient blows up for $u_\mathrm{p} \to 1$, it is expected that the volume densities in the model remain within physical constraints, i.e. $u_i \leq 1$.
\par
To apply \textit{Onsager}'s variational principle, we need the time derivatives of the energy functionals \eqref{def:kinetic_energy}-\eqref{def:magnetic_energy}. We compute, using \eqref{eq:individual_mass_balance_with_averaged_velocity}, \eqref{eq:v_divergence_free},  \eqref{eq:momentum_balance} and integration by parts that
\begin{equation}
  \mathcal{E}'_\textrm{kin}(t) = \int_\Om \bff \cdot \bfv  \dx - \int_\Om \bfT : \nabla \bfv \dx.  
\end{equation}
For the mixing energy we obtain
\begin{equation}
\begin{aligned}
\mathcal{E}'_\textrm{mix}(t) &= \frac{k_\mathrm{B} T}{V_\text{SPION}} \int_\Om \partial_t u_\mathrm{p} \ln\left( \frac{u_\mathrm{p}}{1-u_\mathrm{p}} \right) \dx\\ 
&\stackrel{\mathclap{\eqref{eq:individual_mass_balance_with_averaged_velocity}}}{=} - \frac{k_\mathrm{B} T}{V_\text{SPION}} \int_\Om  \nabla \cdot (u_\mathrm{p} \bfv + u_\mathrm{p} \vprel) \ln\left( \frac{u_\mathrm{p}}{1-u_\mathrm{p}} \right) \dx \\
&= \frac{k_\mathrm{B} T}{V_\text{SPION}} \int_\Om (\bfv + \vprel) \cdot \frac{1}{1-u_\mathrm{p}}\nabla u_\mathrm{p} \dx = \frac{k_\mathrm{B} T}{V_\text{SPION}} \int_\Om \vprel \cdot \frac{1}{1-u_\mathrm{p}}\nabla u_\mathrm{p} \dx.
\end{aligned}
\end{equation}
For the last equality, we have used the fact that $\frac{1}{1-u_\mathrm{p}} \nabla u_\mathrm{p}$ has a potential and consequently, after another integration by parts, the term involving $\nabla \cdot \bfv$ vanishes. For the time derivatives of the magnetic energies, we compute
\begin{equation}
\red{\mathcal{E}'_\textrm{mag,1}(t)} \red{= - \mu_0 \int_\Omega \mathbf{h} \cdot \partial_t \bfh \, \textrm{d}x,}
\end{equation}
\begin{equation}
    \begin{aligned}
        \mathcal{E}'_\textrm{mag,2}(t) &=  - \partial_t \left( \frac{\mu_0 M_\mathrm{s}}{\xi} \int_\Om u_\mathrm{p} \ln \left(\frac{\sinh(\xi |\bfh|)}{\xi |\bfh|} \right) \dx \right)\\
        &= - \mu_0 M_\mathrm{s}  \left( \frac{1}{\xi} \int_\Om \partial_t u_\mathrm{p} \ln \left(\frac{\sinh(\xi |\bfh|)}{\xi |\bfh|} \right) \dx +  \int_\Om u_\mathrm{p}  \frac{\mathcal{L}(\xi |\bfh|)}{|\bfh|} \bfh \cdot \partial_t \bfh \dx \right)\\ &\stackrel{\mathclap{\eqref{eq:individual_mass_balance_with_averaged_velocity}}}{=}  \mu_0 M_\mathrm{s}  \left( \frac{1}{\xi} \int_\Om \nabla \cdot (u_\mathrm{p} \bfv + u_\mathrm{p} \vprel) \ln \left(\frac{\sinh(\xi |\bfh|)}{\xi |\bfh|} \right) \dx - \int_\Om u_\mathrm{p}  \frac{\mathcal{L}(\xi |\bfh|)}{|\bfh|} \bfh \cdot \partial_t \bfh \dx \right) \\
        &= -\mu_0  \int_\Om  (\bfm(u_\mathrm{p},\bfh) \cdot \nabla) \bfh \cdot (\bfv + \vprel) \dx  - \mu_0\int_\Om \bfm(u_\mathrm{p},\bfh) \cdot \partial_t \bfh \dx, 
    \end{aligned}
\end{equation}
and
\begin{equation}
    \red{\mathcal{E}'_\textrm{mag,3}(t)} \red{= - \int_{\partial \Om} \partial_t \mathbf{b}_\text{ext} \cdot \bfn \phi \,\textrm{d}\sigma - \int_{\partial \Om} \mathbf{b}_\text{ext} \cdot \bfn \partial_t  \phi \, \textrm{d}\sigma}.
\end{equation}

\red{Adding the magnetic energies, we obtain, using $\bfb = \mu_0 (\bfh + \bfm)$, $\mathbf{h} = -\nabla \phi$, integration by parts and $\nabla \cdot \mathbf{b} = 0$, as well as the boundary condition \eqref{eq:maxwell:b}, that}
\begin{equation} \label{Energy_Dissipation_Law}
    \begin{aligned}
    \red{\mathcal{E}'_\textrm{mag,1}(t)} &  \red{+ \mathcal{E}'_\textrm{mag,2}(t) + \mathcal{E}'_\textrm{mag,3}(t)} \\ 
    & \red{=-\mu_0  \int_\Om  (\bfm(u_\mathrm{p},\bfh) \cdot \nabla) \bfh \cdot (\bfv + \vprel) \dx + \int_\Om \bfb \cdot  \nabla \partial_t \phi \textrm{d}x - \int_{\partial \Om} \partial_t \mathbf{b}_\text{ext} \cdot \bfn \phi \,\textrm{d}\sigma - \int_{\partial \Om} \mathbf{b}_\text{ext} \cdot \bfn \partial_t  \phi \, \textrm{d}\sigma}\\
    &\red{= -\mu_0  \int_\Om  (\bfm(u_\mathrm{p},\bfh) \cdot \nabla) \bfh \cdot (\bfv + \vprel) \dx - \int_{\partial \Om} \partial_t \mathbf{b}_\text{ext} \cdot \bfn \phi \,\textrm{d}\sigma}.
\end{aligned}
\end{equation}
Now, we are equipped to derive expressions for $\bfT$ and $\vprel$ using \textit{Onsager}'s variational principle, i.e. we determine $\bfT$ and $\vprel$ such that
\begin{equation}\label{eq:onsagers_vp}
\delta_{(\bfT,\vprel)} \left( \mathcal{E}'_\text{kin}  + \mathcal{E}'_\text{mix} + \red{\mathcal{E}'_\text{mag,1} + \mathcal{E}'_\text{mag,2} + \mathcal{E}'_\text{mag,3}}+ \mathcal{D}_\text{kin} +  \mathcal{D}_\text{drag}\right) \stackrel{!}{=} 0,
\end{equation}
where $\delta_{(\bfT,\vprel)}$ denotes the functional derivative with respect to  $\bfT$ and $\vprel$. From condition \eqref{eq:onsagers_vp}, we obtain
\begin{equation}\label{eq:T}
\bfT = 2\eta(u_\mathrm{p}) \bfD(\bfv),
\end{equation}
where $\bfD(\bfv):= \frac{1}{2}\left(\nabla \bfv + (\nabla \bfv)^T \right)$ denotes the symmetrized gradient and
\begin{equation} \label{def:vprel}
\vprel = \frac{1}{u_\mathrm{p}}\left( - D \nabla u_\mathrm{p} +  \kappa (1-u_\mathrm{p}) (\bfm(u_\mathrm{p},\bfh) \cdot \nabla) \bfh \right),
\end{equation}
with 
\begin{equation}\label{eq:D_and_kappa}
    D := \frac{k_\mathrm{B} T}{6 \pi \eta_\mathrm{f} r_\text{SPION}}\quad  \text{and} \quad \kappa := \frac{\mu_0 V_\text{SPION}}{6 \pi \eta_\mathrm{f} r_\text{SPION}}.
\end{equation}
\red{We emphasize that the expression for $\vprel$ is derived only formally and that, throughout the governing equations, it appears exclusively multiplied by $u_\mathrm{p}$, thereby canceling the factor $1/u_\mathrm{p}$ and ensuring that the model remains well-defined also for $u_\mathrm{p} = 0$.}
Finally, using \eqref{eq:T} and \eqref{def:vprel}, we obtain
\begin{equation}
    \begin{aligned}
    \mathcal{E}'_\text{kin}   + \mathcal{E}'_\text{mix} &+ \red{\mathcal{E}'_\text{mag,1} + \mathcal{E}'_\text{mag,2} + \mathcal{E}'_\text{mag,3}}+ 2\mathcal{D}_\text{kin} +  2\mathcal{D}_\text{drag} \\
    &=  -\mu_0  \int_\Om  (\bfm(u_\mathrm{p},\bfh) \cdot \nabla) \bfh \cdot \bfv \dx  \red{ - \int_{\partial \Om} \partial_t \mathbf{b}_\text{ext} \cdot \bfn \phi \,\textrm{d}\sigma} + \int_\Om \bff \cdot \bfv \dx.
\end{aligned}
\end{equation}

In line with the first law of thermodynamics, the right-hand side should represent the rate of work being done on the system, cf. \cite{grun2009}, which is in our case solely due to the \red{temporal change of} the external magnetic field, i.e. due to the term \red{$ - \int_{\partial \Om} \partial_t \mathbf{b}_\text{ext} \cdot \bfn \phi \,\textrm{d}\sigma$}. Accordingly, we require
\begin{equation}
    \int_\Om (\bff - \mu_0 (\bfm(u_\mathrm{p},\bfh) \cdot \nabla) \bfh) \cdot \bfv \dx = 0,
\end{equation}
and therefore set 
\begin{equation}\label{def:f}
    \bff = \mu_0 (\bfm(u_\mathrm{p},\bfh) \cdot \nabla) \bfh,
\end{equation}
\red{to obtain the overall energy-dissipation balance
\begin{equation} \label{eq:energy_balance}
    \mathcal{E}'_\text{kin}   + \mathcal{E}'_\text{mix} + \red{\mathcal{E}'_\text{mag,1} + \mathcal{E}'_\text{mag,2} + \mathcal{E}'_\text{mag,3}}+ 2\mathcal{D}_\text{kin} +  2\mathcal{D}_\text{drag} = - \int_{\partial \Om} \partial_t \mathbf{b}_\text{ext} \cdot \bfn \phi \,\textrm{d}\sigma.
\end{equation}}
%
% %
To conclude, let us summarize the derived magneto two-phase flow system. For $I= (0,T_\text{end})$, $T_\textrm{end}>0$, we look for the magnetic potential $\phi\colon I\times \Omega \rightarrow \R$, the magnetic field strength $\bfh \colon I \times \Om \rightarrow \R^n$, the volume-averaged velocity $\bfv \colon I \times \Omega \rightarrow \R^n$, the pressure $p\colon I \times \Om \rightarrow \R$, and the volume density of the nanoparticle dispersion $u_\textrm{p}\colon I \times \Omega \rightarrow [0,1]$ such that
\begin{subequations}
\label{eq:two-phase_flow_system_with_relations}
\begin{align}
- \nabla \cdot \left( \mu_0 \left( 1 + M_\mathrm{s} u_\mathrm{p} \frac{\mathcal{L}(\xi |\nabla \phi|)}{|\nabla \phi|}\right) \nabla \phi \right) &= 0 & \text{ in } &I \times \Omega,\label{eq:two-phase_flow_system_with_relations:a} \\
\bfh &= - \nabla \phi & \mbox{ in } & I \times \Om,\label{eq:two-phase_flow_system_with_relations:b} \\
\rho(u_\mathrm{p}) \dt \bfv + \biggl(\Bigl( \rho(u_\mathrm{p}) \bfv + (\trho_\mathrm{p} - \trho_\mathrm{f}) \Bigl( - D\nabla u_\mathrm{p}  + \kappa (1-u_\mathrm{p}) (\bfm(u_\mathrm{p},\bfh) \cdot &\nabla) \bfh \Bigr) \Bigr) \cdot \nabla \biggr) \bfv  & & \nonumber\\...- \nabla \cdot (2\eta(u_\mathrm{p}) \bfD(\bfv)) + \nabla p &=  \mu_0 (\bfm(u_\mathrm{p},\bfh) \cdot \nabla) \bfh  & \mbox{ in } & I \times \Om,\label{eq:two-phase_flow_system_with_relations:c}\\
\nabla \cdot \bfv &= 0 & \mbox{ in } & I \times \Om,\label{eq:two-phase_flow_system_with_relations:d}\\
\dt u_\mathrm{p} + \nabla \cdot  \Bigl(- D\nabla u_\mathrm{p} +  u_\mathrm{p} \bfv  + \kappa (1-u_\mathrm{p}) (\bfm(u_\mathrm{p},\bfh) \cdot \nabla) \bfh \Bigr) &= 0 & \mbox{ in } & I \times \Om,\label{eq:two-phase_flow_system_with_relations:e}
\end{align}
\end{subequations}
with
\begin{equation}
        \rho(u_\mathrm{p}) = (\trho_\mathrm{p}-\trho_\mathrm{f})u_\mathrm{p} + \trho_\mathrm{f}, 
\end{equation}
\begin{equation}
    \bfm(u_\textrm{p},\bfh) = M_\textrm{s} u_\textrm{p} \mathcal{L}(\xi |\bfh|) \frac{\bfh}{|\bfh|},
\end{equation}
and
\begin{equation}
    \eta(u_\mathrm{p}) = \eta_\mathrm{f}(1+2.5u_\mathrm{p}),
\end{equation}
subject to the boundary conditions given in Equations \eqref{eq:magnetic_field_bc} and \eqref{eq:velocity_SPIONs_bc}.

\begin{remark} \label{remark:couplings}
\red{
    The present model \eqref{eq:two-phase_flow_system_with_relations:a}-\eqref{eq:two-phase_flow_system_with_relations:e} features various coupling terms describing the interaction of the SPION dispersion, the flow and the magnetic field. In addition to the usual coupling terms resulting from the transport of nanoparticles by convection in the flow and magnetophoretic drift, which appear in all previous models, additional mechanisms that are modeled here, are:
    \begin{itemize}
        \item The SPION-dependent part of the coefficient in \eqref{eq:two-phase_flow_system_with_relations:a} accounts for the feedback of the local SPION distribution on the magnetic field. More precisely, the SPION dispersion volume fraction determines the local magnetization and thereby modifies the constitutive relation between magnetic field strength and magnetic flux density.
        \item The SPION-dependent mixture density in \eqref{eq:two-phase_flow_system_with_relations:c} accounts for changes in the local inertia of the mixture caused by variations of the SPION dispersion volume density.
        \item The SPION-dependent viscosity in \eqref{eq:two-phase_flow_system_with_relations:c} accounts for changes in the effctive viscous response of the mixture caused by variations of the SPION dispersion volume density.
        \item The additional inertia term in \eqref{eq:two-phase_flow_system_with_relations:c} accounts for the transport of momentum associated with the relative mass flux between the SPION and carrier-fluid phases. This relative mass flux is due to both diffusion and magnetophoretic drift of the SPION dispersion.
        \item The SPION-dependent magnetic body force on the right-hand side of \eqref{eq:two-phase_flow_system_with_relations:c} accounts for the force exerted by the magnetic field on the magnetized SPION phase and its resulting transfer of momentum to the mixture.
        \item The factor $(1-u_\mathrm{p})$ in the relative SPION flux in \eqref{eq:two-phase_flow_system_with_relations:e} accounts for finite-volume effects: it represents the locally available carrier-fluid fraction and suppresses relative particle transport as the SPION dispersion reaches its physical limit of a dense packing.
    \end{itemize}
    A comprehensive term-by-term comparison of the presented model with related models, see \cite{grief2005,gitter2012,nacev2011,reinelt2023, badfar2020, boutopoulos2020, weiss2021, grun2019}, is given in Table \ref{tab:literature_comparison}. In Section \ref{sec:relative_importance_study}, we investigate the relative importance of these model components for the outcome of our simulation.}
\end{remark}

\begin{sidewaystable}
    \centering
    \small
    \setlength{\tabcolsep}{4pt}
    %\color{red}
    \renewcommand{\arraystretch}{3}
    \caption{\red{Comparison of mathematical models for magnetic nanoparticle transport.
    A checkmark indicates that the respective coupling or effect is included,
    whereas ``--'' indicates that it is neglected.}}
    \label{tab:literature_comparison}
    \begin{tabular}{
    l
    c
    c
    cccc
    c
    l
    l
}
    \toprule

    \makecell[l]{Model /\\reference}
    &
    \makecell{Flow $\to$\\SPIONs}
    &
    \makecell{Magnetic field\\$\to$ SPIONs}
    &
    \multicolumn{4}{c}{SPIONs $\to$ Flow}
    &
    \makecell{SPIONs $\to$\\Magnetic field}
    &
    \makecell[l]{Magnetization\\model}
    &
    \makecell[l]{Crowding /\\Volume filling}
    \\

    \cmidrule(lr){4-7}

    &
    &
    &
    \makecell{rel. momentum\\flux}
    &
    \makecell{variable\\density}
    &
    \makecell{variable\\viscosity}
    &
    \makecell{magnetic\\body force}
    &
    &
    &
    \\

    \midrule
    
        Badfar et al. \cite{badfar2020}
        & \yes
        & \yes
        & \no
        & \yes
        & \yes
        & \yes
        & \no
        & \textit{Langevin}
        & \no
        \\

        Boutopoulos et al. \cite{boutopoulos2020}
        & \yes
        & \yes
        & \no
        & \yes
        & \no
        & \yes
        & \no
        & linear
        & \no
        \\

        Gitter--Odenbach \cite{gitter2012}
        & \yes
        & \yes
        & \no
        & \no
        & \no
        & \yes
        & \no
        & \textit{Langevin}
        & \no
        \\

        Grief--Richardson \cite{grief2005}
        & \yes
        & \yes
        & \no
        & \no
        & \no
        & \no
        & \no
        & \textit{Langevin}
        & \no
        \\

        Gr\"un--Wei{\ss} \cite{grun2019}
        & \yes
        & \yes
        & \no
        & \no
        & \no
        & \yes
        & \yes
        & \makecell[l]{dynamic relaxation\\to Langevin equilibrium}
        & \yes
        \\

        Nacev et al. \cite{nacev2011}
        & \yes
        & \yes
        & \no
        & \no
        & \no
        & \no
        & \no
        & \makecell[l]{linear /\\with saturation}
        & \no
        \\

        Reinelt et al. \cite{reinelt2023}
        & \yes
        & \yes
        & \no
        & \no
        & \no
        & \no
        & \yes
        & \textit{Langevin}
        & \no
        \\

        Wei{\ss}, Model B \cite{weiss2021}
        & \yes
        & \yes
        & \yes
        & \yes
        & \no
        & \yes
        & \yes
        & linear
        & \yes
        \\

        Present model
        & \yes
        & \yes
        & \yes
        & \yes
        & \yes
        & \yes
        & \yes
        & \textit{Langevin}
        & \yes
        \\

        \bottomrule
    \end{tabular}
\end{sidewaystable}

%%%%%%%%%%%%%%%%%%%%%%%%%%%%%%%%%%%%%%%%%%%%%%%%%%%%%%%%%%%%%%%%%%%%%%%%%%%%%%%
\subsection{Non-Dimensionalization}%%%%%%%%%%%%%%%%%%%%%%%%%%%%%%%%%%%%%%%%%%%%%%%%%%
\label{sec:non_dimensionalization}%%%%%%%%%%%%%%%%%%%%%%%%%%%%%%%%%%%%%%%%%%%%%%%%%%%
%%%%%%%%%%%%%%%%%%%%%%%%%%%%%%%%%%%%%%%%%%%%%%%%%%%%%%%%%%%%%%%%%%%%%%%%%%%%%%%

Let us introduce the following non-dimensional quantities:
\begin{gather} \label{eq:nondimensional_quantities}
\widehat{x} = \frac{x}{L_\textrm{ref}}, \quad \widehat{\bfv} = \frac{\bfv}{V_\textrm{ref}}, \quad \widehat{t} = \frac{t}{T_\textrm{ref}}, \quad \widehat{p} = \frac{p}{\trho_f V_\textrm{ref}^2}, \quad \widehat{\bfh} = \frac{\bfh}{H_\textrm{ref}}, \quad \red{\widehat{\phi} = \frac{\phi}{H_\textrm{ref}L_\textrm{ref}}}, 
\end{gather}
where $L_\textrm{ref}, V_\textrm{ref}$ and $H_\textrm{ref}$ are reference values for the length, velocity and magnetic field strength. The reference time is given by $T_\textrm{ref} = L_\textrm{ref}/V_\textrm{ref}$. This gives rise to the non-dimensional version of system \eqref{eq:two-phase_flow_system_with_relations:a}-\eqref{eq:two-phase_flow_system_with_relations:e}:
\begin{subequations}
\label{eq:two-phase_flow_system_with_relations_non_dimensional}
\begin{align}
- \widehat{\nabla} \cdot \left( \left( 1 + \textrm{M} u_\mathrm{p} \widehat{\mathcal{L}}(\widehat{\xi} |\widehat{\nabla}\widehat{\phi}|)\right) \widehat{\nabla} \widehat{\phi} \right) & = 0 & \mbox{ in } &  \widehat{I} \times \widehat{\Om},\label{eq:two-phase_flow_system_with_relations_non_dimensional:a}\\
\widehat{\bfh} &= -\widehat{\nabla} \widehat{\phi} & \mbox{ in } &  \widehat{I} \times \widehat{\Om},\label{eq:two-phase_flow_system_with_relations_non_dimensional:b}\\
\widehat{\rho}(u_\mathrm{p}) \partial_{\hat{t}} \widehat{\bfv} + \biggl(\biggl( \widehat{\rho}(u_\mathrm{p}) \widehat{\bfv} + \widehat{\rho}_\text{rel} \biggl(- \frac{1}{\textrm{Pe}} \nabla u_\mathrm{p} + \frac{1}{\textrm{Ke}_*}u_\mathrm{p}(1-u_\mathrm{p}) \widehat{\mathbf{K}}(\widehat{\bfh}) \biggr) \biggr) \cdot \widehat{\nabla} \biggr) \widehat{\bfv} & &  \nonumber\\...- \widehat{\nabla} \cdot \left(\frac{2\widehat{\eta}(u_\mathrm{p})}{\textrm{Re}} \widehat{\bfD}(\widehat{\bfv})\right) + \widehat{\nabla} \widehat{p} &=  \frac{1}{\textrm{Ke}} u_\mathrm{p} \widehat{\mathbf{K}}(\widehat{\bfh})  & \mbox{ in } & \widehat{I} \times \widehat{\Om},\label{eq:two-phase_flow_system_with_relations_non_dimensional:c}\\
 \widehat{\nabla} \cdot\widehat{\bfv} &= 0 & \mbox{ in } &  \widehat{I} \times \widehat{\Om},\label{eq:two-phase_flow_system_with_relations_non_dimensional:d}\\
\partial_{\hat{t}}u_\mathrm{p} + \widehat{\nabla} \cdot  \biggl(- \frac{1}{\textrm{Pe}} \widehat{\nabla}  u_\mathrm{p}  + u_\mathrm{p} \widehat{\bfv}  + \frac{1}{\textrm{Ke}_*}u_\mathrm{p}(1-u_\mathrm{p}) \widehat{\mathbf{K}}(\widehat{\bfh}) \biggr) &= 0 & \mbox{ in } &  \hat{I} \times \widehat{\Om},\label{eq:two-phase_flow_system_with_relations_non_dimensional:e}
\end{align}
\end{subequations}
with nondimensional characteristic values 
\begin{equation}\label{eq:characteristic_values}
\begin{aligned}
    \textrm{M} &= M_\mathrm{s} \xi, \quad \widehat{\xi} = \xi H_\textrm{ref}, \quad 
    \widehat{\rho}_\text{rel} = \frac{\tilde{\rho}_\mathrm{p}}{\tilde{\rho}_\mathrm{f}} - 1, \quad
    \textrm{Pe} = \frac{V_\textrm{ref}L_\textrm{ref}}{D},
    \\
    \textrm{Ke}_* &= \frac{V_\textrm{ref}L_\textrm{ref}}
    {\kappa M_\mathrm{s} \xi H_\textrm{ref}^2}, \quad
    \textrm{Re} = \frac{\tilde{\rho}_\mathrm{f} V_\textrm{ref}L_\textrm{ref}}{\eta_\mathrm{f}}, \quad
    \textrm{Ke} = \frac{\tilde{\rho}_\mathrm{f} V_\textrm{ref}^2}
    {\mu_0 M_\mathrm{s} \xi H_\textrm{ref}^2},
\end{aligned}
\end{equation}
and functions
\begin{gather}
    \widehat{\rho}(u) = \widehat{\rho}_\textrm{rel}\, u + 1, \quad \widehat{\mathbf{K}}(\widehat{\bfh}) = \widehat{\mathcal{L}}(\widehat{\xi}|\widehat{\bfh}|) (\widehat{\bfh} \cdot \widehat{\nabla}) \widehat{\bfh}, \quad \widehat{\eta}(u) = 1 + 2.5 u,
\end{gather}
where $\widehat{\mathcal{L}}(\alpha) = \mathcal{L}(\alpha) / \alpha$ \red{with continuous extension $\widehat{\mathcal{L}}(0):= 1/3$}.
In the following, we suppress the bars for the unknowns and the differential operators for a clearer presentation. 

%%%%%%%%%%%%%%%%%%%%%%%%%%%%%%%%%%%%%%%%%%%%%%%%%%%%%%%%%%%%%%%%%%%%%%%%%%%%%%%
\section{Numerical Simulation of Magentically Driven SPION Transport in a Closed System}%%%%%%%%%%%%%%%%%%%%%%%%%%%%%%%%%%%%%%%%
\label{sec:numerical_simulation}%%%%%%%%%%%%%%%%%%%%%%%%%%%%%%%%%%%%%%
%%%%%%%%%%%%%%%%%%%%%%%%%%%%%%%%%%%%%%%%%%%%%%%%%%%%%%%%%%%%%%%%%%%%%%%%%%%%%%%

After introducing a staggered finite element scheme based on a one-step $\theta$-time discretization and flux-corrected transport (FCT) for the SPION transport equation, we present numerical simulations using physical parameters representative of magnetic drug targeting applications. In addition to reporting several physically relevant observables, we investigate to what extent the energy-dissipation law derived in \eqref{eq:energy_balance} is reproduced by the numerical scheme. To this end, we perform experimental temporal convergence studies for the energy balance as well as for the individual energy and dissipation contributions involved. Finally, we conduct an ablation study to assess the influence of the various coupling mechanisms on the simulation results.

\subsection{Discretization of the Magneto Two-Phase Flow Model}
\label{sec:discretization}

Let $\mathcal{T}_h$ denote a conforming quadrilateral finite element mesh of $\Om$, with mesh size $h = \max_{K\in\mathcal{T}_h}\operatorname{diam}(K)$. For $m \in \N$, we define the general \textit{Lagrange} finite element space
\begin{equation}
   \mathcal{S}_{h,m} = \left\{ \psi_h \in \mathcal{C}^0(\overline{\Om}) \mid  \restr{\psi_h}{K}  \in \mathcal{Q}_m(K) \; \forall K \in \mathcal{T}_h \right\}, 
\end{equation}
where $\mathcal{Q}_m$ denotes the tensor-product polynomial space of degree at most $m$ in each spatial dimension, and consider specifically
\begin{gather}\label{eq:FE_spaces}
\mathcal{V}_{h} = \mathcal{S}_{h,2}, \quad \boldsymbol{\mathcal{V}}_{h} = \mathcal{V}^3_{h}, \quad \boldsymbol{\mathcal{V}}_{h,0} = \left\{ \bfv_h \in \boldsymbol{\mathcal{V}}_{h} \mid \bfv_h = 0 \text{  on  } \partial \Omega \right\}, \\
\mathcal{W}_{h} = \mathcal{S}_{h,1}, \quad \boldsymbol{\mathcal{W}}_{h} := \mathcal{W}^3_{h}.
\end{gather}
For end time $T_\text{end}>0$ determine the nondimensional end time $\widehat{T}_\text{end} = T_\text{end}/T_\textrm{ref}$ and define the time step size $\tau := \widehat{T}_\text{end}/N$ with $N \in \N$. For all unknowns $z \in \{\phi, \bfh, \bfv, p, u_\mathrm{p}\}$ we use the notation
\begin{equation}
    z^n_h \approx z(t^n)
\end{equation}
to denote the corresponding finite element approximation at time $t^n := n\tau$, $n = 0,1,...,N$. We note that we drop the index $\mathrm{p}$ of the SPION dispersion volume density $u_\mathrm{p}$ in the discrete setting for the sake of a clearer presentation.

Starting with initial data $(\mathbf{v}^0_h, u^0_h) \in \boldsymbol{\mathcal{V}}_{h,0} \times \mathcal{W}_{h}$, compute $(\phi^n_h, \bfh^n_h, \mathbf{v}^n_h, p^n_h, u^n_h) \in \mathcal{V}_{h} \times \boldsymbol{\mathcal{W}}_{h} \times \boldsymbol{\mathcal{V}}_{h,0} \times\mathcal{W}_{h} \times \mathcal{W}_{h}$ for $n = 1,2,...,N$ according to the following staggered scheme with $\theta \in [0,1]$:
\begin{enumerate}
    \item Find $\phi^n_h \in \mathcal{V}_{h}$ such that
    \begin{equation}\label{eq:FE_phi}
        \left(\big[ 1+\textrm{M}u^{n-1}_h\widehat{\mathcal{L}}(\widehat{\xi}|\nabla \phi^n_h|) \big] \nabla \phi^n_h, \nabla \psi_h\right)_{L^2} = -(\mu_0H_\textrm{ref})^{-1}\int_{\partial \Omega}\bfb_\text{ext} \cdot \bfn\; \psi_h \; \textrm{d}\sigma, \quad \forall \psi_h \in \mathcal{V}_{h}.
    \end{equation}
    \item Find $\bfh_h^n \in \boldsymbol{\mathcal{W}}_{h}$ such that
    \begin{equation}\label{eq:FE_h}
        (\bfh_h^n, \mathbf{g}_h)_{L^2} = - (\nabla \phi^n_h, \mathbf{g}_h)_{L^2}, \quad \forall \mathbf{g}_h \in \boldsymbol{\mathcal{W}}_{h}.
    \end{equation}
    \item Find $(\bfv_h^n, p^n_h) \in \boldsymbol{\mathcal{V}}_{h,0} \times \mathcal{W}_{h}$ such that
    \begin{equation} \label{eq:FE_vp}
    \begin{aligned}
        \frac{1}{\tau} (\rho_h^{n-1} \mathbf{v}^n_h, \mathbf{w}_h)_{L^2} + c(\mathbf{m}^n_{h,\theta}; \;&\bfv^n_{h,\theta}, \mathbf{w}_h) + \frac{2}{\textrm{Re}}( \eta_h^{n-1} \mathbf{D}(\bfv^n_{h,\theta}), \mathbf{D}(\mathbf{w}_h))_{L^2}  - (p^n_h, \nabla \cdot \mathbf{w}_h)_{L^2} - (\nabla  \cdot\bfv^n_h, q_h)_{L^2}\\
        &=\frac{1}{\textrm{Ke}}(u^{n-1}_h \mathbf{K}^n_h, \mathbf{w}_h)_{L^2} + \frac{1}{\tau}(\rho^{n-1} \bfv^{n-1}_h, \mathbf{w}_h)_{L^2}, \quad \forall (\mathbf{w}_h,q_h) \in \boldsymbol{\mathcal{V}}_{h,0} \times \mathcal{W}_{h},
    \end{aligned}
    \end{equation}
    where we denote 
    \begin{align}
                \rho^{n-1}_h &= \rho(u^{n-1}_h) = \widehat{\rho}_\text{rel} \; u^{n-1}_h + 1,\\
                c(\mathbf{u}; \mathbf{v}, \mathbf{w}) &= \frac12 [ ( (\mathbf{u} \cdot \nabla ) \bfv  , \mathbf{w})_{L^2} - ( (\mathbf{u} \cdot \nabla ) \mathbf{w} , \bfv )_{L^2}],\\
                \mathbf{m}^n_{h,\theta} &= \rho^{n-1}\bfv^n_{h, \theta} + \mathbf{J}^n_h, \\
                \bfv^n_{h,\theta} &= \theta \bfv^n_h + (1-\theta)\bfv^{n-1}_h, \\
                \mathbf{J}^n_h &= \widehat{\rho}_\text{rel} \left(-\frac{1}{\textrm{Pe}} \nabla u^{n-1}_h + \frac{1}{\textrm{Ke}_*}u^{n-1}_h(1-u^{n-1}_h) \mathbf{K}^n_h\right),\\
                \eta^{n-1}_h&= \eta(u^{n-1}) = 1 + 2.5 u^{n-1}_h,\\
                \mathbf{K}^n_h &= \widehat{\mathcal{L}}(\widehat{\xi} |\bfh^n_h|)(\bfh^n_h \cdot \nabla) \bfh^n_h.
    \end{align}
    \item Find $u^n_h \in \mathcal{W}_{h}$ such that
    \begin{equation}\label{eq:FE_u}
    \begin{aligned}
        \frac{1}{\tau} ( u^n_h, s_h)_{L^2} &+ \theta \left( \frac{1}{\textrm{Pe}} \nabla u^n_h - u^n_h \mathbf{v}^n_h - \frac{1}{\textrm{Ke}_*}u^n_h (1-u^n_h) \mathbf{K}^n_h, \nabla s_h \right)_{L^2} 
        = \frac{1}{\tau}( u^{n-1}_h, s_h)_{L^2}\\ &- (1-\theta) \left( \frac{1}{\textrm{Pe}} \nabla u^{n-1}_h - u^{n-1}_h \mathbf{v}^n_h - \frac{1}{\textrm{Ke}_*}u^{n-1}_h (1-u^{n-1}_h) \mathbf{K}^n_h, \nabla s_h \right)_{L^2}, \quad \forall s_h \in \mathcal{W}_{h}.
    \end{aligned}
    \end{equation}
\end{enumerate}
Here, $(.,.)_{L^2}$ denotes the inner product in $L^2(\Om)$ and is used for scalar, vector and matrix valued functions, depending on the context. 
\begin{remark}
    The proposed staggered scheme offers several practical advantages. Its modular structure allows the individual subproblems for the magnetic field, fluid flow, and SPION transport to be treated separately using established numerical techniques. This considerably simplifies the implementation and avoids the need for a fully coupled multiphysics solver. At the same time, the staggered treatment reduces the degree of nonlinearity that has to be handled within each individual solve, which is particularly advantageous for prospective computationally demanding 3D simulations of magnetic drug targeting experiments. Moreover, the one-step $\theta$ method provides some flexibility in the temporal discretization. In particular, choosing $\theta=0.5$ may yield improved temporal accuracy compared with the fully implicit choice $\theta=1$, as also suggested by the temporal convergence studies presented in Section \ref{sec:convergence_study}. These advantages come at the cost of additional approximation errors introduced by the sequential treatment of the coupled subproblems. In particular, the staggered coupling gives rise to splitting errors and introduces a temporal inconsistency between quantities that are evaluated at different stages of a time step. Consequently, the fully discrete staggered scheme does not directly inherit the energy-dissipation law \eqref{eq:energy_balance} of the continuous model. In particular, no exact discrete counterpart of the continuous energy balance is enforced by construction. The extent to which the numerical solution nevertheless reproduces the energy-dissipation behavior is assessed experimentally, see Figure \ref{fig:energy_balance} and Section \ref{sec:convergence_study}.
\end{remark}
\subsubsection{Stabilization of the SPION Transport Problem}
\label{sec:discretization:transportStabilization}

Realistic physical parameters for Magnetic Drug Targeting applications lead to Problem \eqref{eq:two-phase_flow_system_with_relations_non_dimensional:e} being dominated by convective transport with velocity $\bfv$, see also Table \ref{tab:reference_characteristic_values}. For transport-dominated problems, standard continuous \textit{Galerkin} finite element methods suffer from spurious numerical oscillations near sharp gradients and can lead to nonphysical over- and undershoots, including negative values of the solution even when the exact solution remains nonnegative. Various stabilization techniques have been developed to address these difficulties, see, e.g., \cite{john2008}. In the present work, we use the flux-corrected transport (FCT) method for the SPION transport equation arising in our model \cite{boris1973, kuzmin2002, kuzmin2003}.
\par 
For a coefficient vector $W \in \R^{N_u}$, we denote by $u_h(W) = \sum_{j=1}^{N_u} W_j \varphi_j \in \mathcal{W}_h$ the associated finite element function, where $\{\varphi_j\}_{j=1}^{N_u}$ is a basis of $\mathcal{W}_h$. We introduce the mass matrix and the stiffness matrix
$\mathbb{M}, \mathbb{A}(W)\in\mathbb R^{N_u\times N_u}$ with entries
\begin{equation}
\begin{aligned}
    \mathbb{M}_{ij} &= (\varphi_j, \varphi_i)_{L^2}, \\
    \mathbb{A}_{ij}(W) &= \left(\frac{1}{\mathrm{Pe}}\nabla\varphi_j - \varphi_j \mathbf v_h^n - \frac{1}{\mathrm{Ke}_*} \varphi_j\bigl(1- u_h(W)\bigr)\mathbf K_h^n, \nabla\varphi_i \right)_{L^2},
\end{aligned}
\end{equation}
for $i,j = 1,..., N_u$.
Then, the algebraic version of the nonlinear problem in Equation \eqref{eq:FE_u} becomes: Given $U^{n-1} \in \R^{N_u}$, find $U^{n} \in \R^{N_u}$ such that
\begin{equation}\label{eq:algebraic_transport_equation}
    (\mathbb{M} + \tau \theta \mathbb{A}(U^n))U^n = (\mathbb{M} - \tau (1-\theta) \mathbb{A}(U^{n-1}))U^{n-1}
\end{equation}
The basic idea of flux-corrected transport (FCT) is now to replace the high-order discretization \eqref{eq:algebraic_transport_equation} by a monotonicity-preserving low-order scheme and to compensate for the additional numerical diffusion by means of suitably limited antidiffusive fluxes. To this end, let $\mathbb{M}_L$ denote the lumped mass matrix and introduce the nonlinear low-order operator
\begin{equation}
\mathbb{L}(W) = \mathbb{A}(W)+\mathbb{D}(W),
\end{equation}
where $\mathbb{D}(W)$ is an artificial diffusion matrix chosen such that the off-diagonal entries of $\mathbb{L}(W)$ are non-positive. Here, the entries of $\mathbb{D}$ are given by
\begin{equation}
    \mathbb{D}_{ij}(W)
    =
    -\max\left\{
        0,\,
        \mathbb{A}_{ij}(W),\,
        \mathbb{A}_{ji}(W)
    \right\},
    \qquad i\neq j,
\end{equation}
and
\begin{equation}
    \mathbb{D}_{ii}(W)
    =
    -\sum_{j\neq i} \mathbb{D}_{ij}(W).
\end{equation}
The FCT problem originating from \eqref{eq:algebraic_transport_equation} then reads
\begin{equation}\label{eq:FCT_problem}
\left(\mathbb{M}_L +
\tau\theta\mathbb{L}(U^n) \right)U^n = \left(\mathbb{M}_L - \tau(1-\theta)\mathbb{L}(U^{n-1}) \right)U^{n-1} + F^*(U^n, U^{n-1})
\end{equation}
where $F^*(U^n, U^{n-1})\in\mathbb{R}^{N_u}$ denotes the limited antidiffusive correction. Here, 
\begin{equation}
    F^*_i = \sum_{j=1}^{N_u} \alpha_{ij} r_{ij},
\end{equation}
where
\begin{equation}
    \begin{aligned}
    r_{ij}(U^n, U^{n-1}) = \mathbb{M}_{ij}&(U^n_i - U^n_j) - \mathbb{M}_{ij}(U^{n-1}_i - U^{n-1}_j)\\
    &+ \tau \theta \mathbb{D}(U^n)_{ij}(U^n_j - U^n_i) + \tau (1-\theta)\mathbb{D}(U^{n-1})_{ij}(U^{n-1}_j - U^{n-1}_i)
    \end{aligned}
\end{equation}
denotes the raw antidiffusive flux between the degrees of freedom $i$ and $j$, and $0\leq \alpha_{ij}=\alpha_{ji}\leq 1$ are flux limiters. For all $\alpha_{ij} = 1$, the original problem \eqref{eq:algebraic_transport_equation} is recovered, whereas all $\alpha_{ij} = 0$ corresponds to the uncorrected low-order scheme. The flux limiters \(\alpha_{ij}\) are determined using \textit{Zalesak}'s algorithm \cite{zalesak1979} such that the antidiffusive fluxes are largely retained in noncritical regions, i.e. there the scheme approaches the original high-order discretization, while the antidiffusive fluxes are reduced in regions where the high-order solution would produce nonphysical over- or undershoots.
To avoid an additional nonlinearity introduced by $F^*(U^n, U^{n-1})$ in Equation \eqref{eq:algebraic_transport_equation}, we use a predictor $\widetilde{U}^n$ defined by
\begin{equation} \label{eq:predictor}
    \left(\mathbb{M}_L + \tau\theta\mathbb{L}(\widetilde{U}^n) \right)\widetilde{U}^n = \left(\mathbb{M}_L - \tau(1-\theta)\mathbb{L}(U^{n-1}) \right)U^{n-1}.
\end{equation}
The predictor is used for the computation of $F^*$ in \eqref{eq:FCT_problem}, leading to the problem
\begin{equation}\label{eq:FCT_with_predicor}
\left(\mathbb{M}_L +
\tau\theta\mathbb{L}(U^n) \right)U^n = \left(\mathbb{M}_L - \tau(1-\theta)\mathbb{L}(U^{n-1}) \right)U^{n-1} + F^*(\widetilde{U}^n, U^{n-1}).
\end{equation}
\subsubsection{Implementation Details}
\label{sec:discretization:implementation}

The numerical scheme is implemented using the \textit{Distributed and Unified Numerics Environment DUNE} in C++ \cite{blatt2025,engwer2018}, making particular use of the discretization module dune-fufem \cite{dune-fufem}. 

The nonlinear problems are treated as follows:
\begin{itemize}
    \item Problem \eqref{eq:FE_phi} is linearized using fixed-point iterations. Since the magnetic potential is determined only up to an additive constant, one degree of freedom is fixed during the solution of each linearized system. At each time step $n$, the iteration is terminated once the algebraic residual vector $\mathbf{r}^{n,k}_\phi$ associated with the nonlinear problem  \eqref{eq:FE_phi} satisfies $\|\mathbf{r}^{n,k}_\phi \| < 10^{-10}$.
    \item Problem \eqref{eq:FE_vp} is linearized using \textit{Newton}'s method. Since the pressure is determined only up to an additive constant, one pressure degree of freedom is fixed during the solution of each linearized system. At each time step $n$, the iteration is terminated once the algebraic residual vector $\mathbf{r}^{n,k}_{\mathrm{NS}}$ associated with the nonlinear problem \eqref{eq:FE_vp} satisfies $\|\mathbf{r}^{n,k}_{\mathrm{NS}}\| < 10^{-10}$.
    \item For Problem \eqref{eq:FE_u}, at each time step $n$, Equation \eqref{eq:predictor} is solved using fixed-point iterations until $\|\widetilde{U}^{n,k} - \widetilde{U}^{n,k-1} \|< 10^{-10}$ to obtain the predictor. Then, Equation \eqref{eq:FCT_with_predicor} is solved using fixed-point iterations until $\|U^{n,k} - U^{n,k-1} \|< 10^{-10}$.
\end{itemize}
Here, the superscript $k$ denotes the iteration index of the respective nonlinear solver and $\|\cdot\|$ is the \textit{Euclidean} norm. The chosen tolerances were verified to be sufficiently small such that a further reduction did not noticeably affect the reported numerical results. The arising linear problems are solved using the direct solver UMFPACK \cite{davis2004} through the interface provided by dune-istl. The visualizations of the finite element solutions are created using ParaView \cite{ayachit2015}.

To ensure a numerically stable evaluation of functions near removable singularities, we use suitable \textit{Taylor} expansions or continuous extensions in the implementation whenever the respective arguments are sufficiently close to zero. In particular, for the modified \textit{Langevin} function $\widehat{\mathcal L}(\alpha)=\mathcal L(\alpha)/\alpha$, we use
\begin{equation}
    \widehat{\mathcal L}(\alpha) = \frac{1}{3}-\frac{\alpha^2}{45} \qquad\text{for }|\alpha|<10^{-5}.
\end{equation}
As part of the post-processing, we also evaluate the energy and dissipation functionals. For the numerical evaluation of $\mathcal{E}_{\mathrm{mag},2}$, we use
\begin{equation}
    \ln\left(\frac{\sinh(\alpha)}{\alpha}\right) = \frac{\alpha^2}{6}-\frac{\alpha^4}{180} \qquad\text{for }|\alpha|<10^{-5}.
\end{equation}
For the mixing energy $\mathcal{E}_\textrm{mix}$, we consider $f(u)=u\ln u+(1-u)\ln(1-u)$. For the pure phases $u=1$ and $u=0$, we use the continuous extension $f(0)=0$ and $f(1)=0$. In our numerical computations, values of $u$ that are negative but of the order of machine precision do occur. For such values we also set $f(u) = 0$ in the evaluation of the mixing energy.
For the numerical evaluation of $\mathcal{D}_{\mathrm{drag}}$, we avoid the explicit computation of the relative particle velocity and instead use the relative particle flux $\mathbf J_{\mathrm{p, rel}} = u_{\mathrm p}\mathbf v_{\mathrm{p,rel}}$. Using the definition of the friction coefficient \eqref{eq:friction_coefficient}, the dissipation can be written as 
\begin{equation}
    \mathcal{D}_{\mathrm{drag}}(t) = \frac{3\pi\eta_{\mathrm f}r_{\mathrm{SPION}}}{V_{\mathrm{SPION}}} \int_\Omega \frac{|\mathbf J_{\mathrm{rel}}|^2} {u_{\mathrm p}(1-u_{\mathrm p})}\,\dx.
\end{equation}
To avoid division by zero near the pure states, the factor $1/[u_{\mathrm p}(1-u_{\mathrm p})]$ is regularized by replacing $u_{\mathrm p}$ with 
\begin{equation}
    u_{\mathrm p,\varepsilon} = \min\!\bigl(1-\varepsilon,\max(\varepsilon,u_{\mathrm p})\bigr), \qquad \varepsilon=10^{-10}.
\end{equation}

\begin{figure}[tbp]
    \centering
    \begin{minipage}[t]{0.32\textwidth}
        \centering
        \makebox[\linewidth][c]{%
             \resizebox{0.95\linewidth}{!}{%
                \begin{tikzpicture}%[scale=1]

    \definecolor{targetbrown}{RGB}{225,205,170}
    \definecolor{innergrid}{RGB}{170,170,170}

    \fill[targetbrown] (0,7.5) rectangle (10,10);

    \foreach \x in {1.25,2.5,3.75,5,6.25,7.5,8.75}{
        \draw[innergrid, line width=0.4pt] (\x,0) -- (\x,10);
    }

    \foreach \y in {1.25,2.5,3.75,5,6.25,7.5,8.75}{
        \draw[innergrid, line width=0.4pt] (0,\y) -- (10,\y);
    }

    % Outer boundary
    \draw[black, line width=1pt] (0,0) rectangle (10,10);
    % Lower boundary of target region
    \draw[brown!70!black, line width=0.8pt] (0,7.5) -- (10,7.5);

    \node at (5,8.75) {\LARGE $\Omega_{\mathrm{target}}$};

    \foreach \x in {0,2.5,5,7.5,10}{
        \draw (\x,0) -- (\x,-0.15);
        \node[below] at (\x,-0.15) {\small $\x$};
    }

    \foreach \y in {0,2.5,5,7.5,10}{
        \draw (0,\y) -- (-0.15,\y);
        \node[left] at (-0.15,\y) {\small $\y$};
    }

    \node[below] at (5,-0.8) {$x_1$};
    \node[rotate=90] at (-0.8,5) {$x_2$};

\end{tikzpicture}%
                }%
        }
        \captionof{figure}{Nondimensional computational domain $[0,10]^2$ with coarse mesh and target region $\Omega_{\mathrm{target}} = [0,10]\times[7.5,10]$.}
        \label{fig:domain-target-region}
    \end{minipage}
    \hfill
    \begin{minipage}[t]{0.64\textwidth}
        \centering
        \begin{subfigure}[t]{0.48\linewidth}
            \centering
            \includegraphics[width=\linewidth]{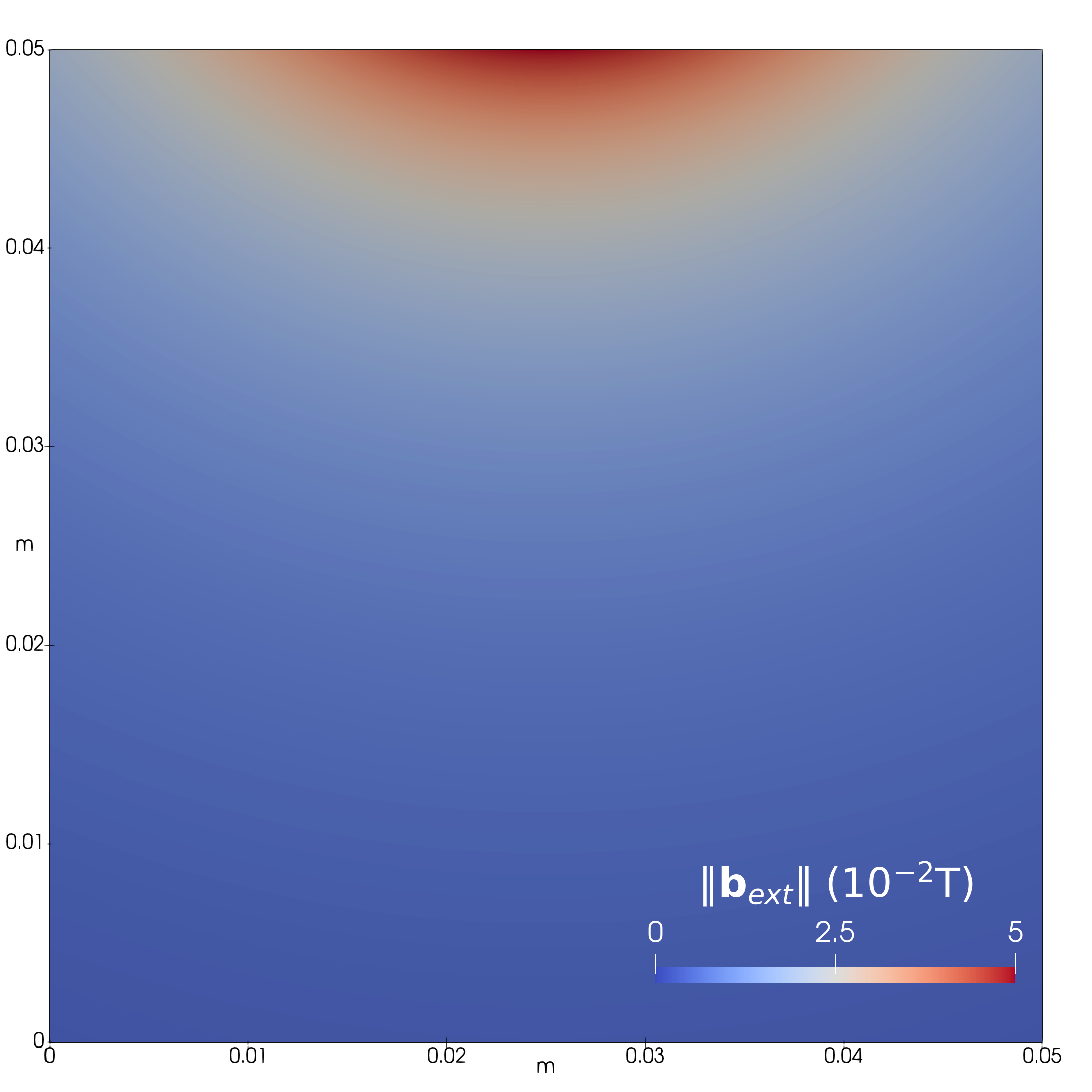}
        \end{subfigure}
        \hfill
        \begin{subfigure}[t]{0.48\linewidth}
            \centering
            \includegraphics[width=\linewidth]{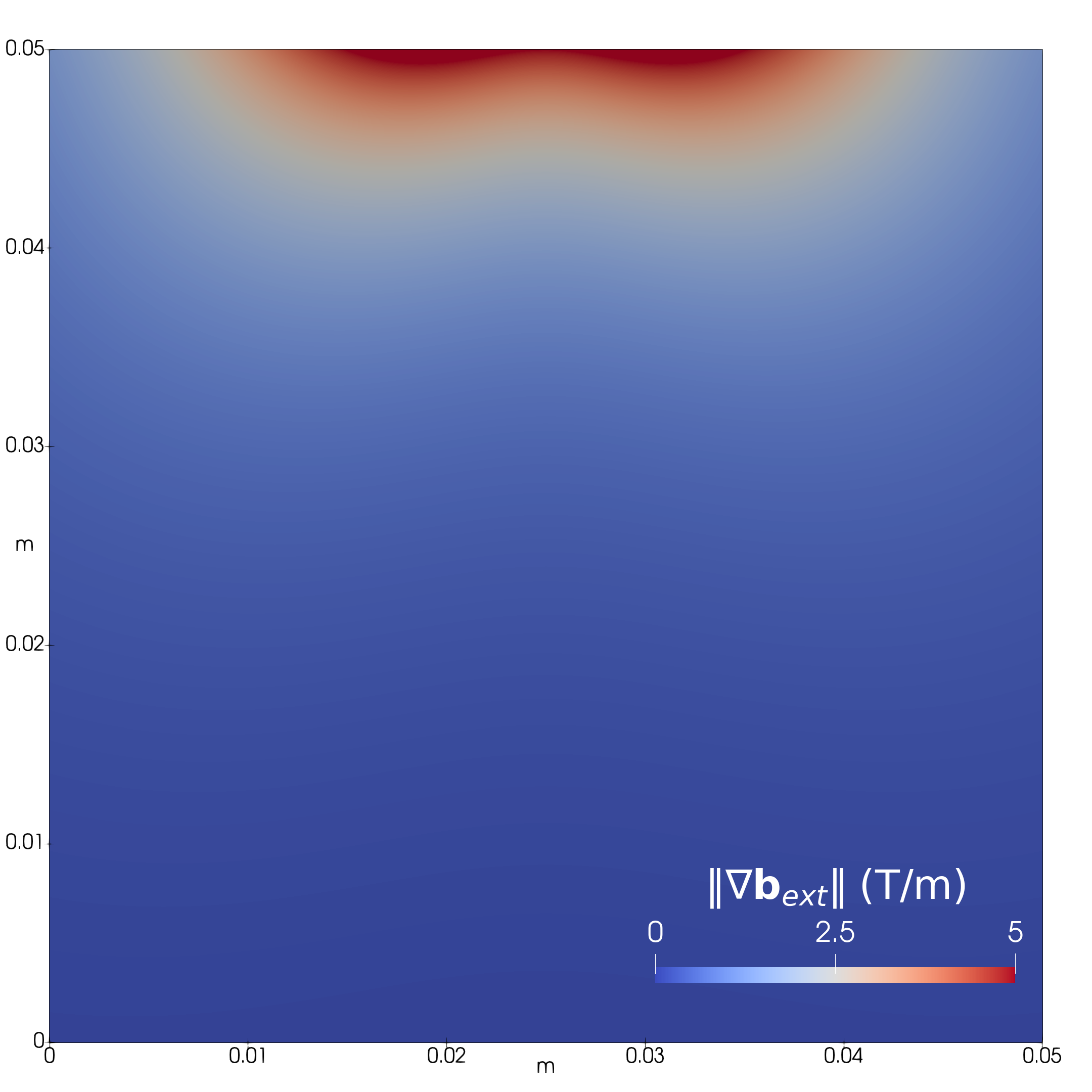}
        \end{subfigure}
        \captionof{figure}{Magnitude of the external magnetic field $\mathbf{b}_\text{ext}$ and its gradient $\nabla\mathbf{b}_\text{ext}$.}
        \label{fig:bExt}
    \end{minipage}
\end{figure}

\subsection{Physical Setting of the Simulation: Magnetic Attraction of a  SPION Dispersion in an Initially Resting Mixture}
\label{sec:setting}

As a numerical experiment, we consider a localized SPION dispersion in a closed, water-filled square domain and simulate its attraction towards a magnet located outside the domain. The physical domain is given by
\begin{equation}
    \Om = [0\textrm{m}, 0.05\textrm{m}] \times [0\textrm{m}, 0.05\textrm{m}].
\end{equation}
We assume that the averaged velocity of the mixture is initially at rest, i.e. 
\begin{equation}
    \mathbf{v}(0, \cdot) = 0 \qquad \mbox{ in } \{0 \} \times \Om.
\end{equation}
To formulate a physically interpretable initial condition for the SPION volume density $u_\mathrm{p}$, let us note that in the lab, SPION dispersions are usually characterized by their iron mass per volume $\rho_\text{Fe}>0$. In our model, this iron mass density is related to the SPION dispersion volume density via
\begin{equation}
\rho_\text{Fe} = \tilde{\rho}_\mathrm{Fe,max} u_\mathrm{p},
\end{equation}
where $\tilde{\rho}_\mathrm{Fe,max} = 0.71 \cdot P_\text{SPION} \varphi_\text{SPION} \tilde{\rho}_\text{cores}$. Here, the factor \(0.71\) represents the mass fraction of iron in the core compounds magnetite/maghemite. In order to simulate an initially localized SPION dispersion at the center of the domain, we set
\begin{equation}
    u_\mathrm{p}(0,x) = \frac{\rho_\text{Fe,initial}}{\tilde{\rho}_\mathrm{Fe,max}} \exp\left(-(0.5\|x - x_0^{u_p} \|)^8  \right)
\end{equation}
where $\rho_\text{Fe,initial} = 1 \mathrm{mg\; ml^{-1}}$, $x_0^{u_p} = (0.025, 0.025)\textrm{m}$ and $\| \cdot \|$ denotes the \textit{Euclidean} norm. As a model for the external magnet, we consider the singular magnetic field 
\begin{equation}
    \bfb_\text{ext}(x) = \frac{\mu_0}{2\pi} \left( \frac{2 (\mathbf{p}_\text{ext} \cdot (x - x_{0}^\bfb)) (x - x_0^\bfb)}{\| x-x_0^\bfb\|^4} - \frac{\mathbf{p}_\text{ext}}{\| x-x_0^\bfb\|^2}\right),
\end{equation}
where $\mathbf{p}_\text{ext} = (0, 100)^T\textrm{Am}$ and $x_0^\bfb = (0.025, 0.07)\textrm{m}$.  This field serves as a two-dimensional analogue of the magnetic field generated by a point dipole, and, in particular, it holds $\nabla \cdot \bfb_\text{ext} = 0$ in $\Omega$, i.e.  compatibility condition \eqref{eq:b_ext_compatibility} is fulfilled. For a visualization of $\bfb_\text{ext}$ and its gradient $\nabla \bfb_\text{ext}$, we refer to Figure \ref{fig:bExt}. For the velocity and the SPION volume density, we impose no-slip and no-flux boundary condition, respectively, on the whole boundary, i.e.
\begin{equation}
    \mathbf{v} = 0, \qquad \Bigl(- D\nabla u_\mathrm{p} +   u_\mathrm{p} \bfv  + \kappa (1- u_\mathrm{p}) (\bfm( u_\mathrm{p},\bfh) \cdot \nabla) \bfh \Bigr)\cdot \bfn = 0 \quad \mbox{ on } I \times \partial \Om.
\end{equation}
Consequently, the system is closed with respect to mass exchange and neither fluid nor SPIONs can enter or leave the domain through its boundary. At the same time, the dynamics inside $\Omega$ are influenced by the externally generated magnetic field.

\begin{table}[tbp]
    \caption{Physical parameters used in the simulations.}
        \centering
    
    {
    \renewcommand{\arraystretch}{1.3}
    \begin{tabularx}{\textwidth}{
        >{\centering\arraybackslash}p{1.7cm}
        >{\centering\arraybackslash}p{2.0cm}
        >{\centering\arraybackslash}p{1.8cm}
        X
        >{\raggedright\arraybackslash}p{2.75cm}
    }
        \toprule
        Parameter & Value & Unit & Description & Reference \\
        \midrule
        
        $T$
        & 293
        & $\mathrm{K}$
        & Absolute temperature of the fluid and nanoparticle dispersion
        & -
        \\

        $\mu_0$
        & $4\pi \times 10^{-7}$
        & $\mathrm{N\,A^{-2}}$
        & Vacuum permeability
        & -
        \\

        $k_B$
        & $1.38 \times 10^{-23}$
        & $\mathrm{J\,K^{-1}}$
        & Boltzmann's constant
        & -
        \\

        $r_\text{multicore}$
        & $5.9 \times 10^{-8}$
        & $\mathrm{m}$
        & Radius of the multicore-cluster in a SPION
        & \cite{reinelt2023}
        \\

        $d_\text{coat}$
        & $2 \times 10^{-9}$
        & $\mathrm{m}$
        & Thickness of SPION coating
        & \cite{reinelt2023}
        \\

        $r_\text{SPION}$
        & $6.1 \times 10^{-8}$
        & $\mathrm{m}$
        & Radius of a SPION
        & \cite{reinelt2023}
        \\

        $V_\text{SPION}$
        & $9.51 \times 10^{-22}$
        & $\mathrm{m^3}$
        & Volume of a SPION
        & computed
        \\

        $P_\text{cores}$
        & $0.7$
        & $-$
        & Packing density of single iron oxide cores inside SPION multicore-cluster
        & estimated
        \\

        $\varphi_\text{cores}$
        & $0.63338$
        & $-$
        & Volume fraction of iron oxide cores inside SPION multicore-cluster
        & computed, Eq. \eqref{eq:varphi_cores}
        \\

        $\tilde{\rho}_\text{coat}$
        & $1.0 \times 10^{3}$
        & $\mathrm{kg\,m^{-3}}$
        & Density of SPION coating agent
        & \cite{reinelt2023}
        \\

        $\tilde{\rho}_\text{cores}$
        & $5.2 \times 10^{3}$
        & $\mathrm{kg\,m^{-3}}$
        & Density of magnetite/maghemite single cores
        & \cite{reinelt2023}
        \\

        $P_\text{SPION}$
        & $0.6$
        & $-$
        & Maximum packing density of SPIONs in the dispersion
        & estimated
        \\
        
        $\tilde{\rho}_\text{SPION}$
        & $3.66 \times 10^{3}$
        & $\mathrm{kg\,m^{-3}}$
        & Density of a SPION
        & computed, Eq. \eqref{eq:rho_spion}
        \\

        $\tilde{\rho}_\text{p}$
        & $2.6 \times 10^{3}$
        & $\mathrm{kg\,m^{-3}}$
        & Density of the densely packed SPION dispersion
        & computed, Eq. \eqref{eq:rho_p}
        \\

        $\tilde{\rho}_\text{f}$
        & $1 \times 10^{3}$
        & $\mathrm{kg\,m^{-3}}$
        & Density of the pure carrier fluid phase
        & -
        \\

        $M_s$
        & $1.74 \times 10^{5}$
        & $\mathrm{A\,m^{-1}}$
        & Saturation magnetization of the densely packed SPION dispersion
        & computed, Eq. \eqref{eq:M_s}
        \\

        $\xi$
        & $8.57 \times 10^{-2}$
        & $\mathrm{m\,A^{-1}}$
        & \textit{Langevin} parameter
        & computed, Eq. \eqref{eq:xi}
        \\

        $p$
        & $2.76 \times 10^{-16}$
        & $\mathrm{A\,m^{2}}$
        & Magnetic moment of a single SPION
        & computed, Eq. \eqref{eq:p}
        \\

        $\eta_f$
        & $1 \times 10^{-3}$
        & $\mathrm{kg\,m^{-1}\;s^{-1}}$
        & Viscosity of the pure carrier fluid phase
        & -
        \\

        $D$
        & $3.51 \times 10^{-12}$
        & $\mathrm{m^2\,s^{-1}}$
        & SPION diffusion coefficient
        & computed, Eq. \eqref{eq:D_and_kappa}
        \\

        $\kappa$
        & $1.04 \times 10^{-18}$
        & $\mathrm{m^4\,A^{-2} \;s^{-1}}$
        & SPION magnetophoresis coefficient
        & computed, Eq. \eqref{eq:D_and_kappa}
        \\
        
        \bottomrule
    \end{tabularx}
    } % end red
    \label{tab:physical_parameters}
\end{table}

\begin{table}[htbp]
    \caption{Reference and characteristic values for the numerical simulations.}
    \centering
    
    \begin{tabular}{lccc}
        \toprule
        Variable & Value & Unit & Equation \\
        \midrule
        $L_\textrm{ref}$
        & $5 \times 10^{-3}$
        & $\mathrm{m}$
        & \eqref{eq:nondimensional_quantities}
        \\

        $V_\textrm{ref}$
        & $1 \times 10^{-2}$
        & $\mathrm{m\;s^{-1}}$
        & \eqref{eq:nondimensional_quantities}
        \\

        $H_\textrm{ref}$
        & $7.96 \times 10^{5}$
        & $\mathrm{A\,m^{-1}}$
        & \eqref{eq:nondimensional_quantities}
        \\

        $\textrm{M}$
        & $2.18\times10^{-1}$
        & $\mathrm{-}$
        & \eqref{eq:characteristic_values}
        \\

        $\widehat{\xi}$
        & $6.82\times10^{4}$
        & $\mathrm{-}$
        & \eqref{eq:characteristic_values}
        \\

        $\widehat{\rho}_\text{rel}$
        & $1.6$
        & $\mathrm{-}$
        & \eqref{eq:characteristic_values}
        \\

        $\textrm{Pe}$
        & $1.42 \times 10^{7} $
        & $\mathrm{-}$
        & \eqref{eq:characteristic_values}
        \\

        $\textrm{Ke}_*$
        & $5.1 \times 10^{-3} $
        & $\mathrm{-}$
        & \eqref{eq:characteristic_values}
        \\

        $\textrm{Re}$
        & $5 \times 10^1 $
        & $\mathrm{-}$
        & \eqref{eq:characteristic_values}
        \\

        $\textrm{Ke}$
        & $8.43 \times 10^{-12} $
        & $\mathrm{-}$
        & \eqref{eq:characteristic_values}
        \\
        
        \bottomrule
    \end{tabular}
    \label{tab:reference_characteristic_values}
\end{table}

\subsection{Simulation Results for the Magnetic Attraction of the SPION Dispersion in a Closed System}

\begin{figure}
    \include{figures/code/numericalExperiment}
    \caption{Snapshots of the magnetic attraction of a SPION dispersion in an initially quiescent mixture. The iron mass density $\rho_\mathrm{Fe}$ is color-coded throughout the computational domain, while the mixture velocity $\bfv$ is represented by white arrows whose lengths are proportional to the velocity magnitude. Black lines represent the magnetic field lines associated with the magnetic field strength $\bfh$.}
    \label{fig:simulation-snapshots}
\end{figure}

\begin{figure}
    \centering
    
    \includegraphics[width=0.24\textwidth]{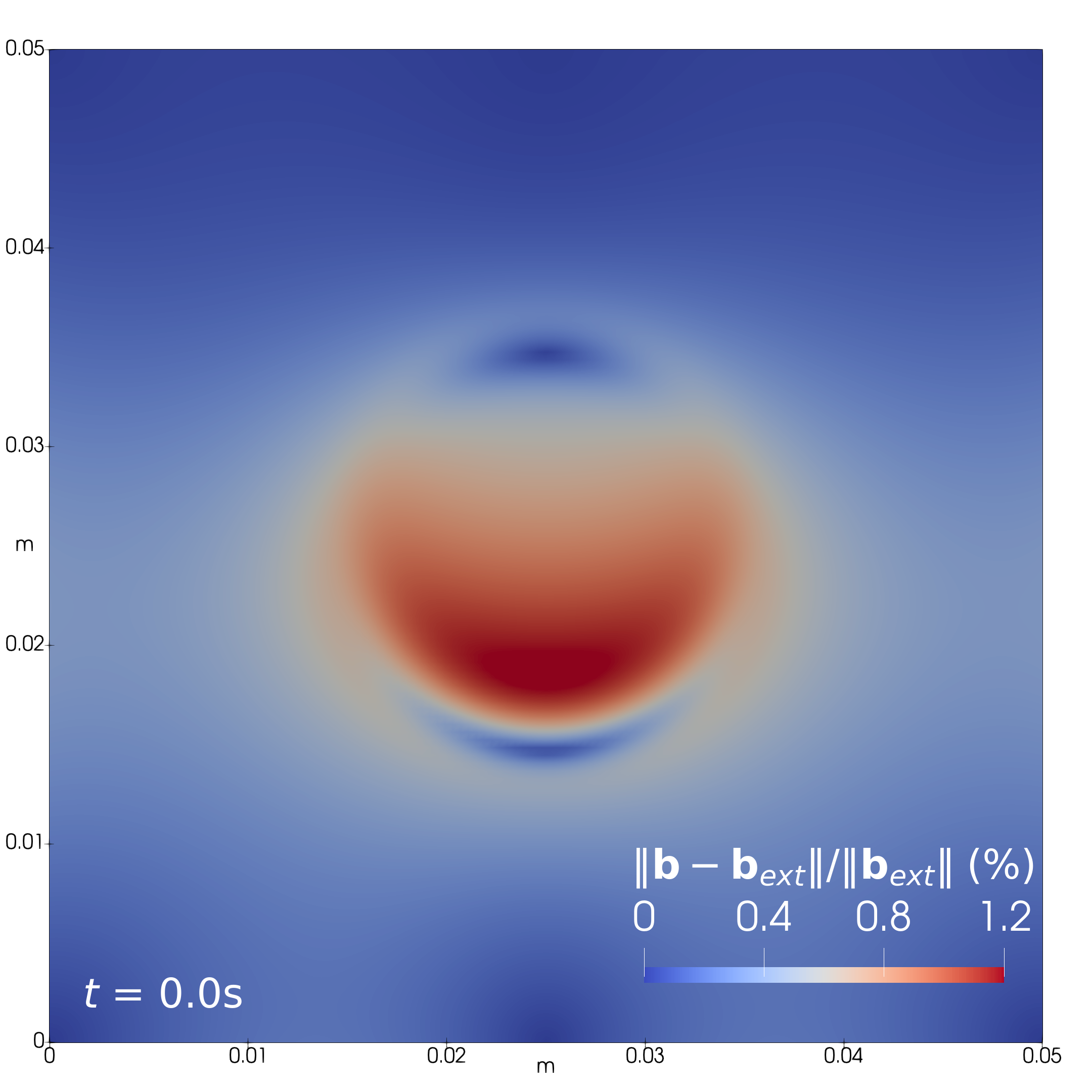}
    \includegraphics[width=0.24\textwidth]{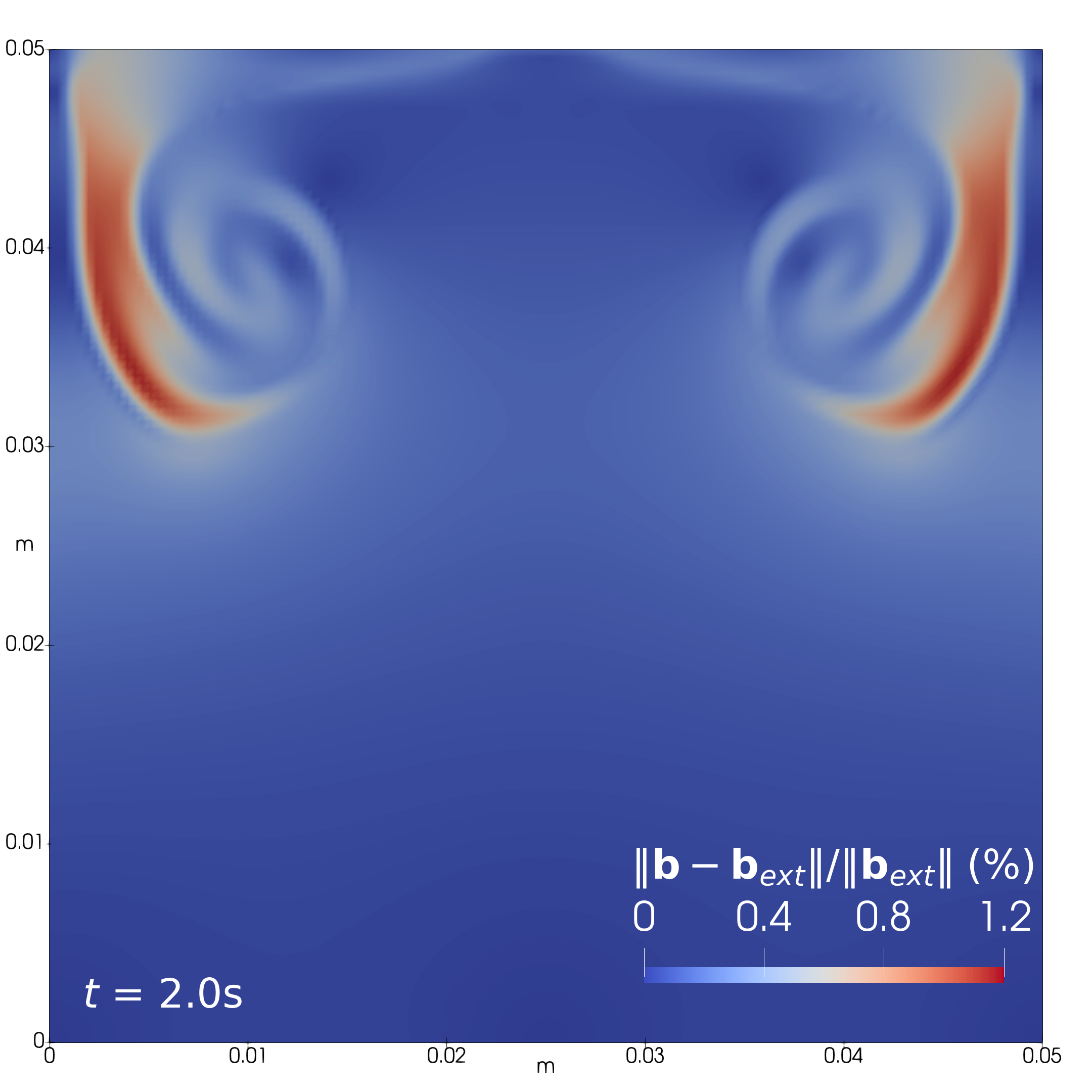}
    \includegraphics[width=0.24\textwidth]{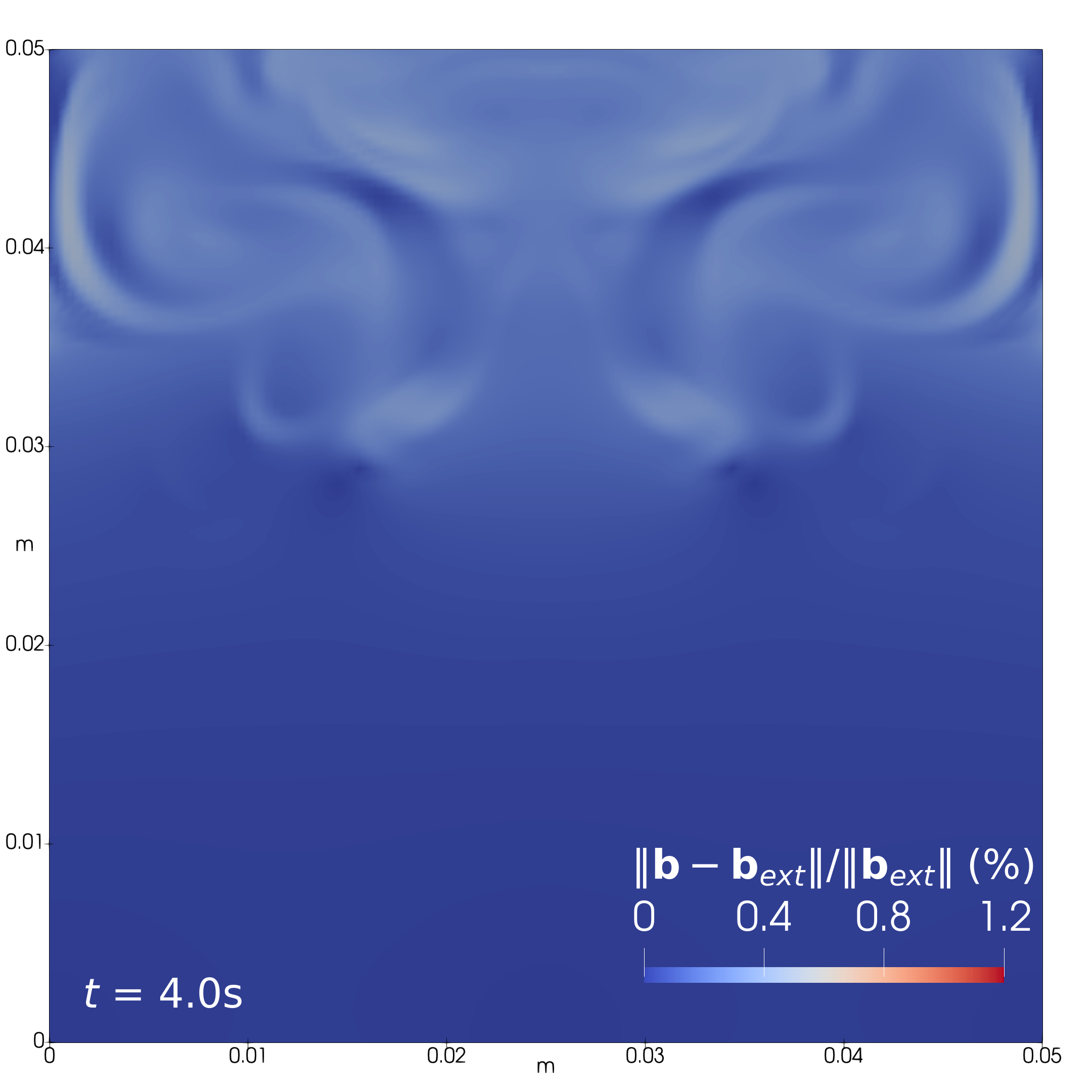}
    \includegraphics[width=0.24\textwidth]{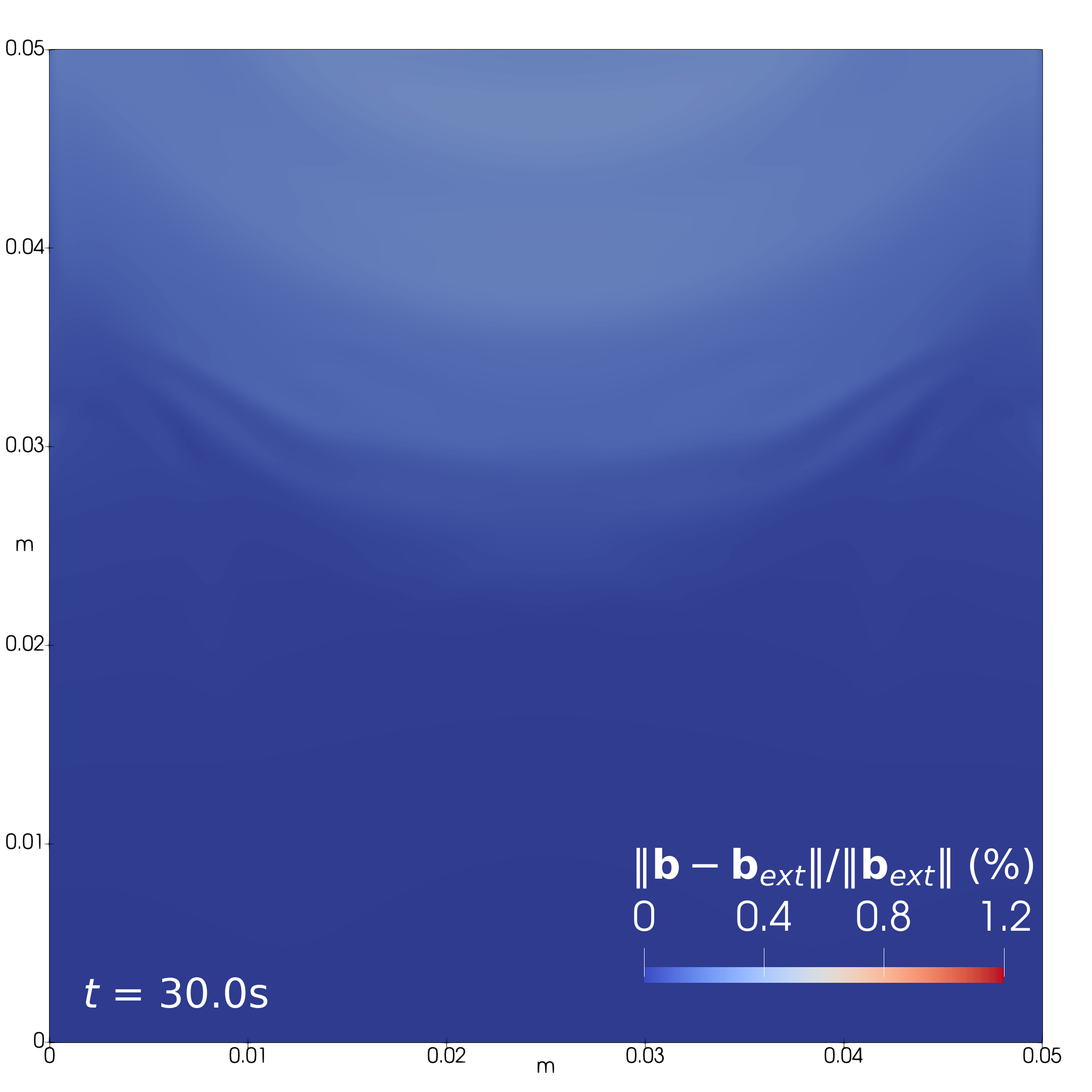}

    \includegraphics[width=0.24\textwidth]{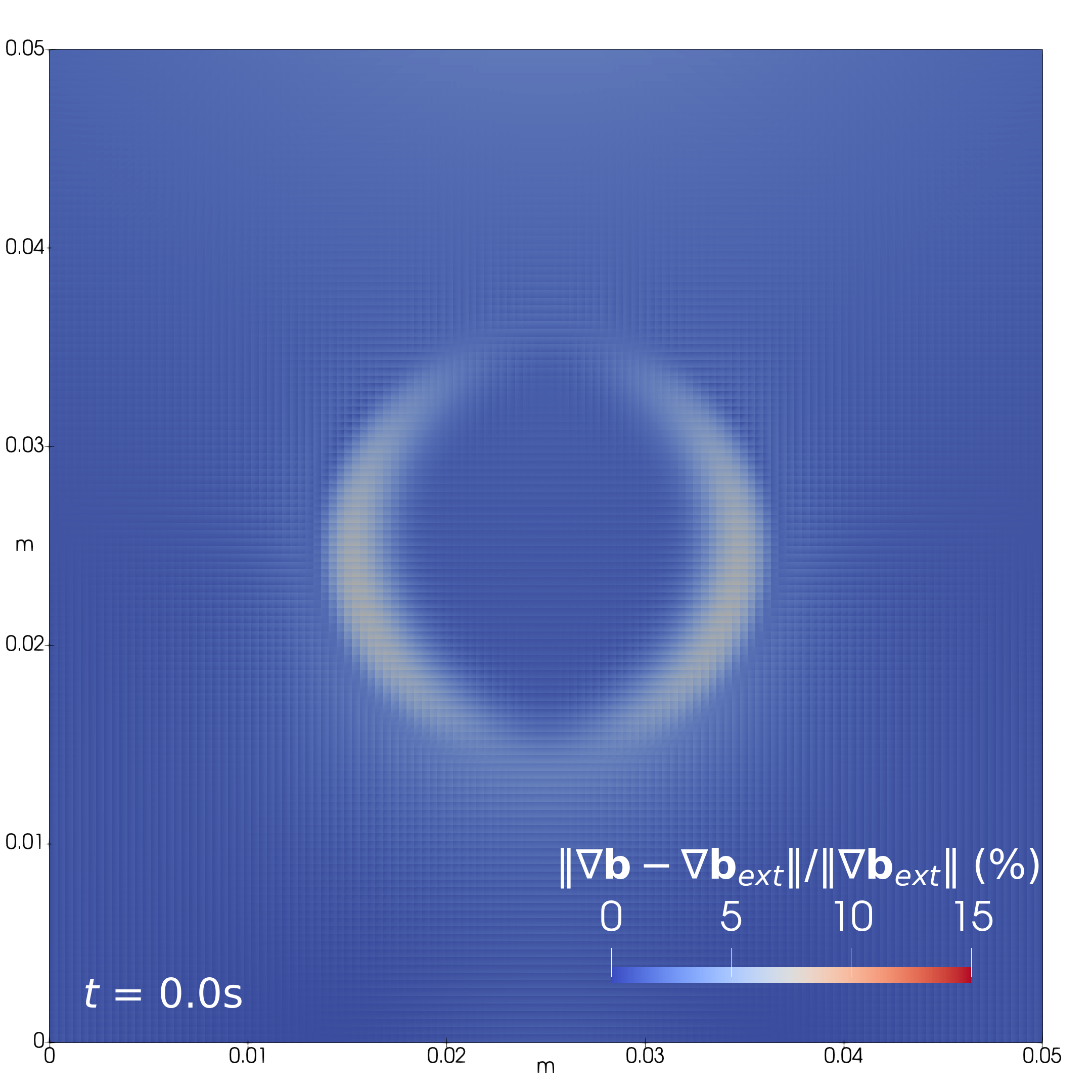}
    \includegraphics[width=0.24\textwidth]{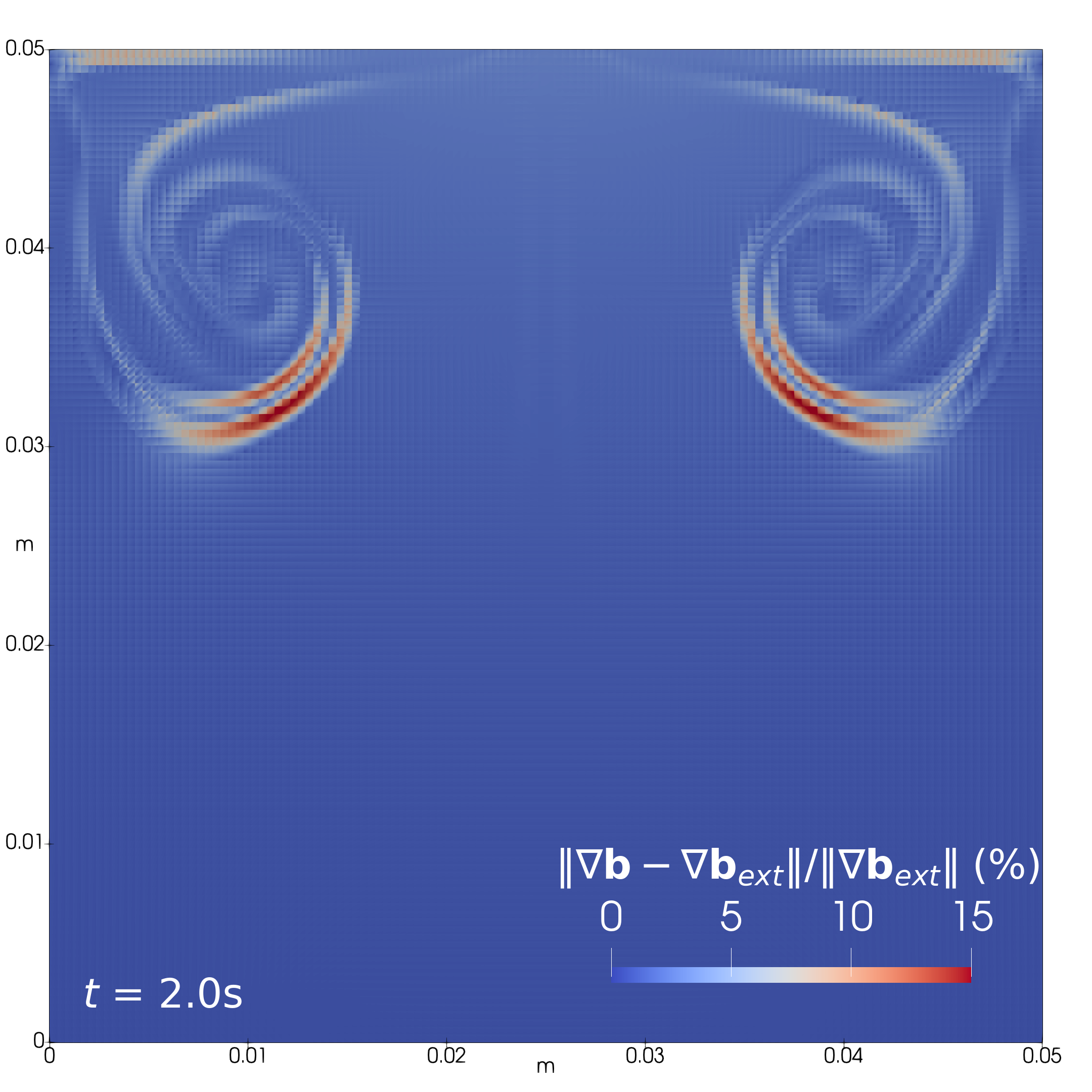}
    \includegraphics[width=0.24\textwidth]{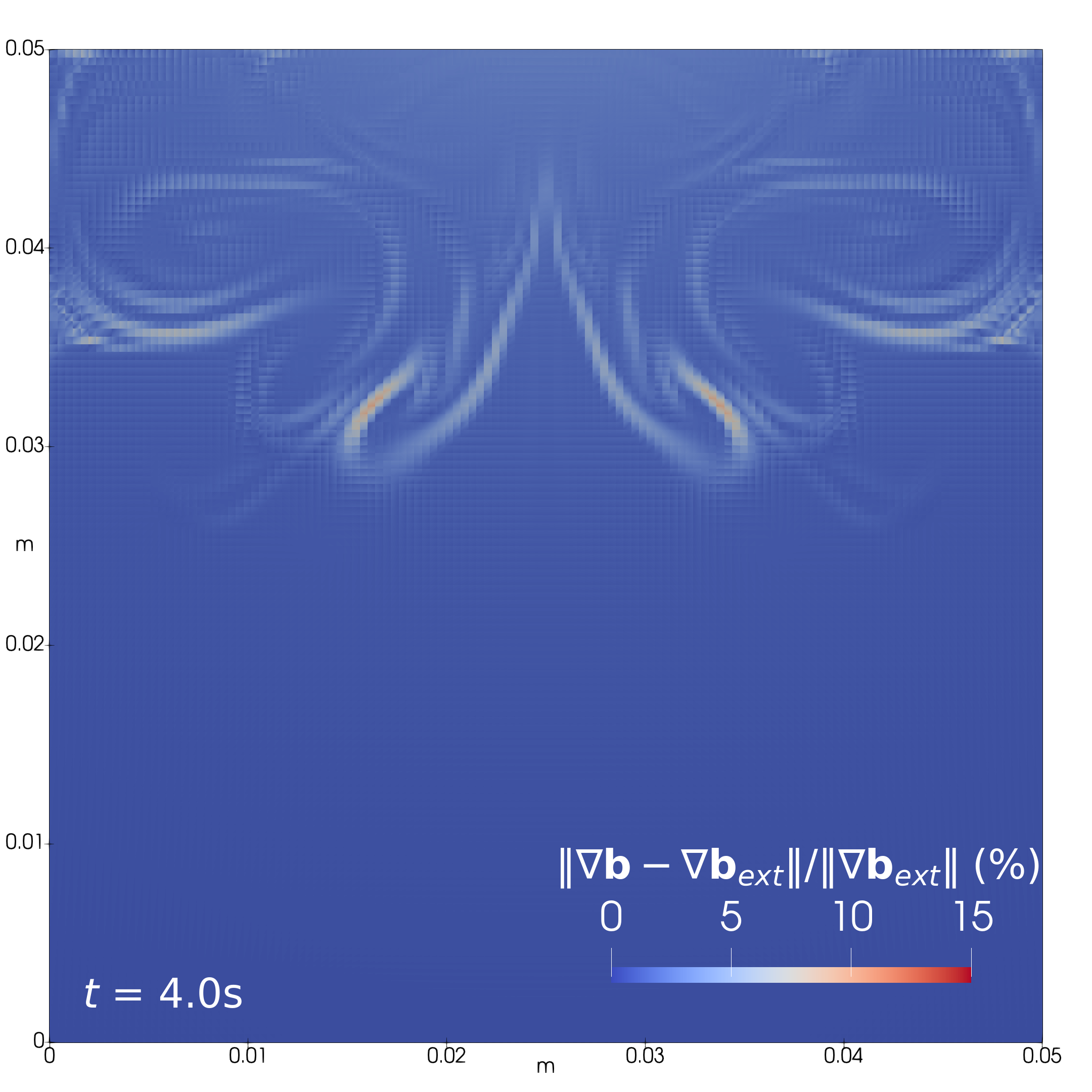}
    \includegraphics[width=0.24\textwidth]{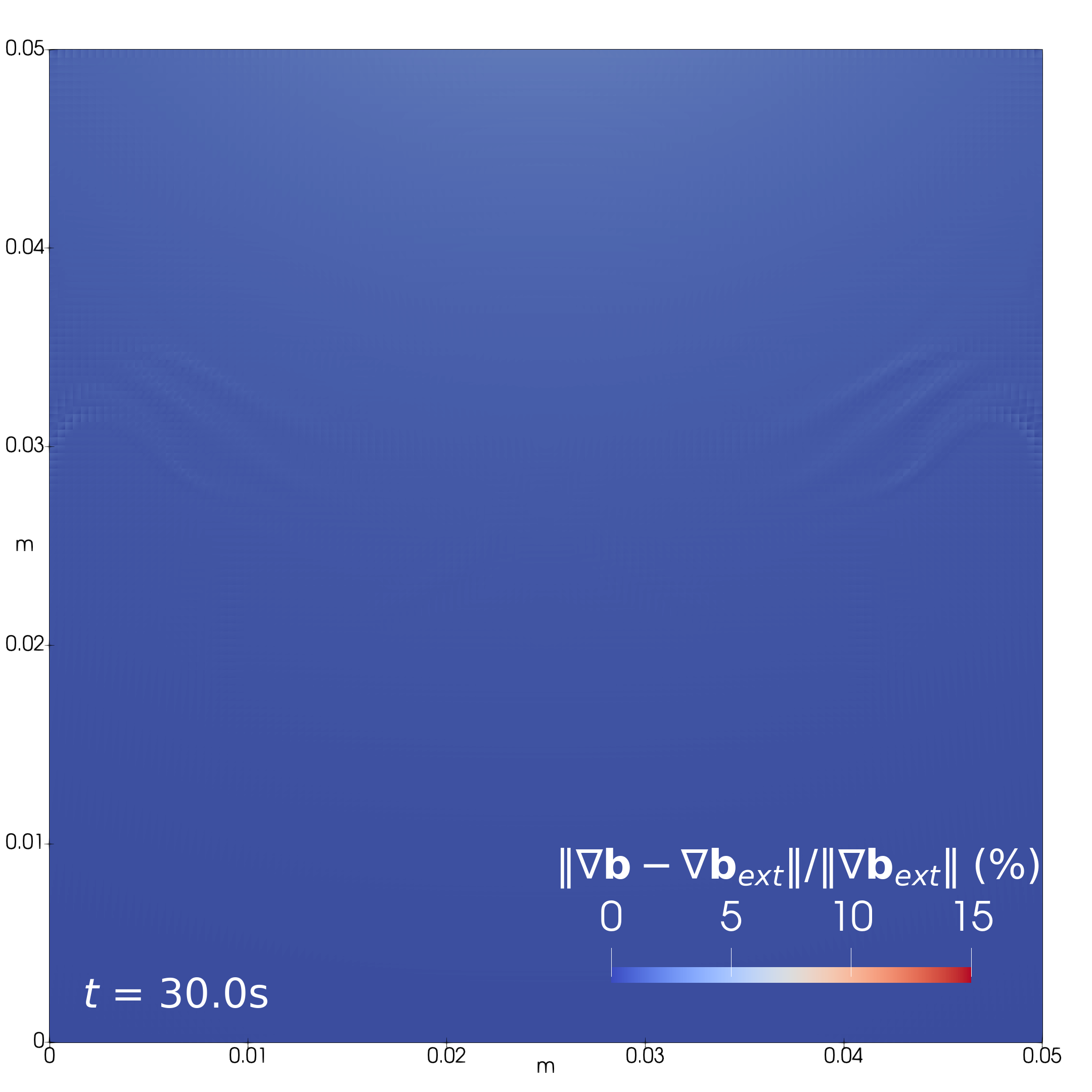}

    \includegraphics[width=0.24\textwidth]{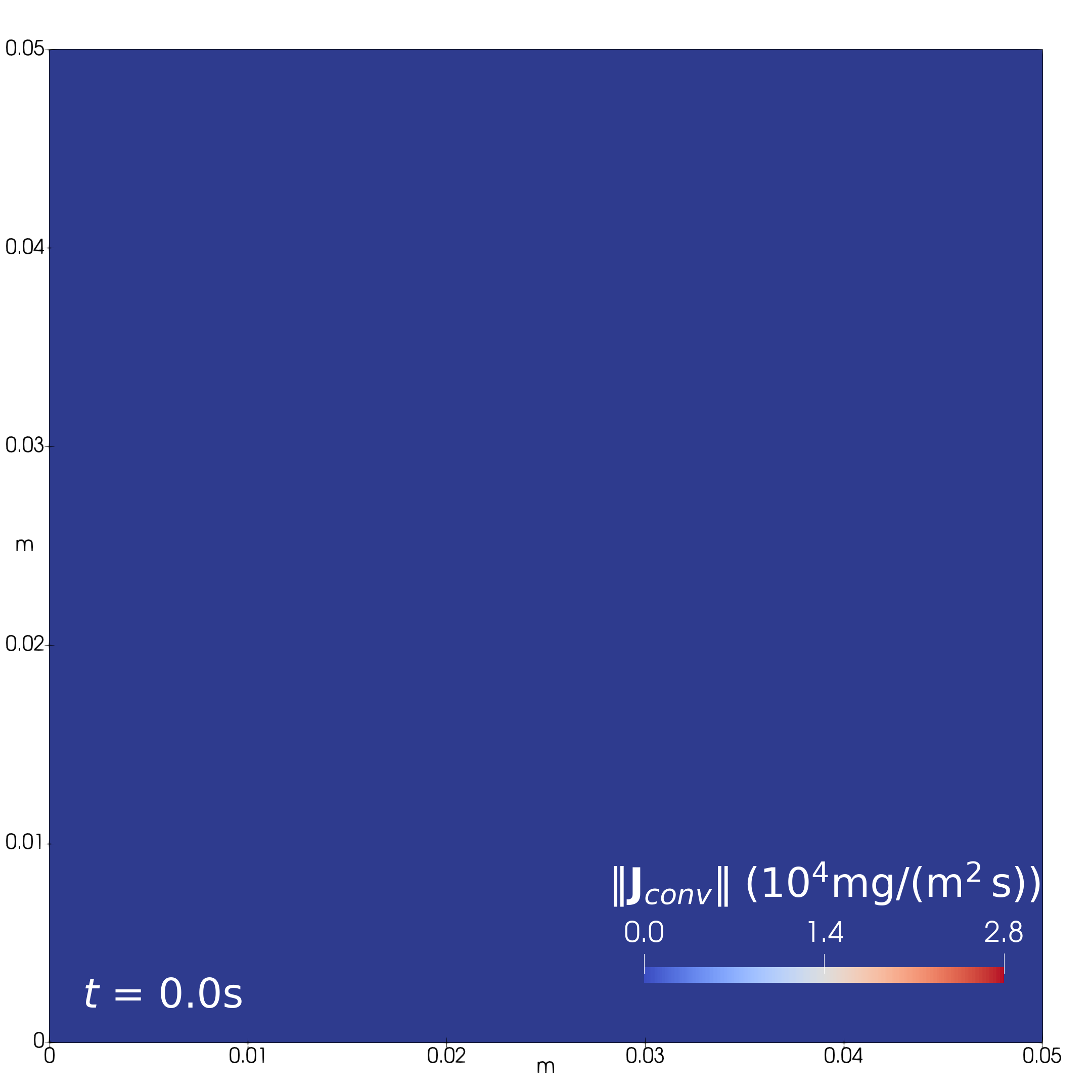}
    \includegraphics[width=0.24\textwidth]{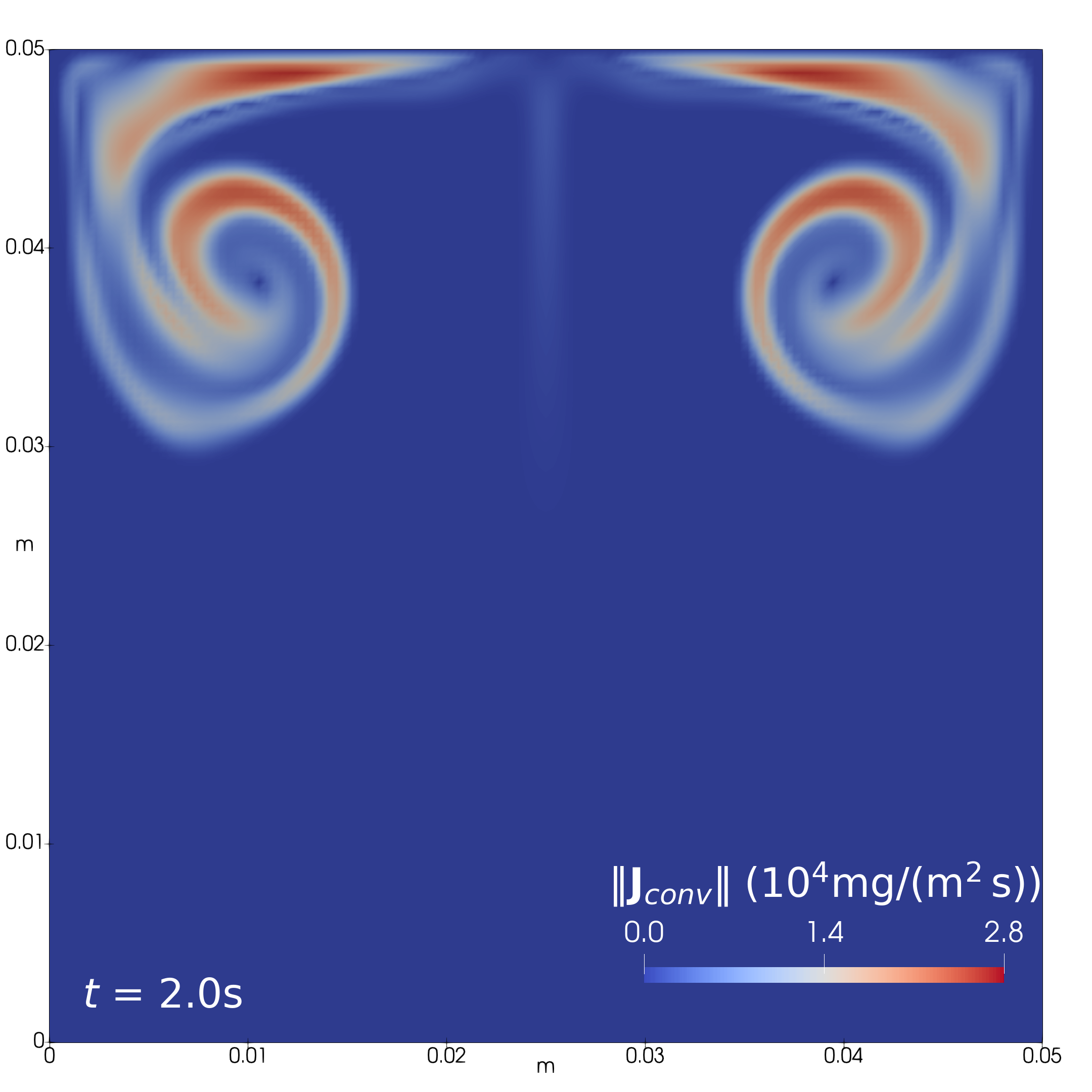}
    \includegraphics[width=0.24\textwidth]{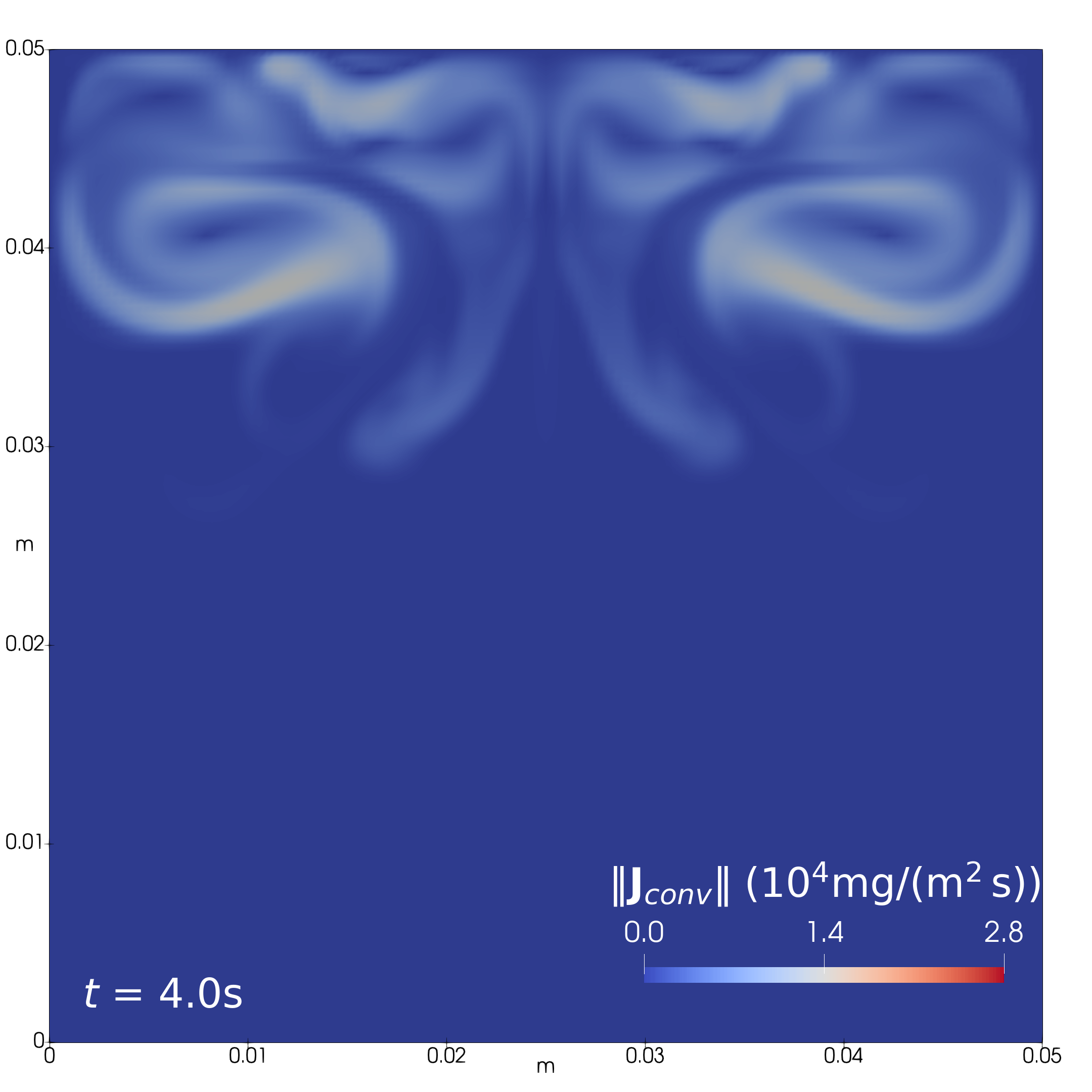}
    \includegraphics[width=0.24\textwidth]{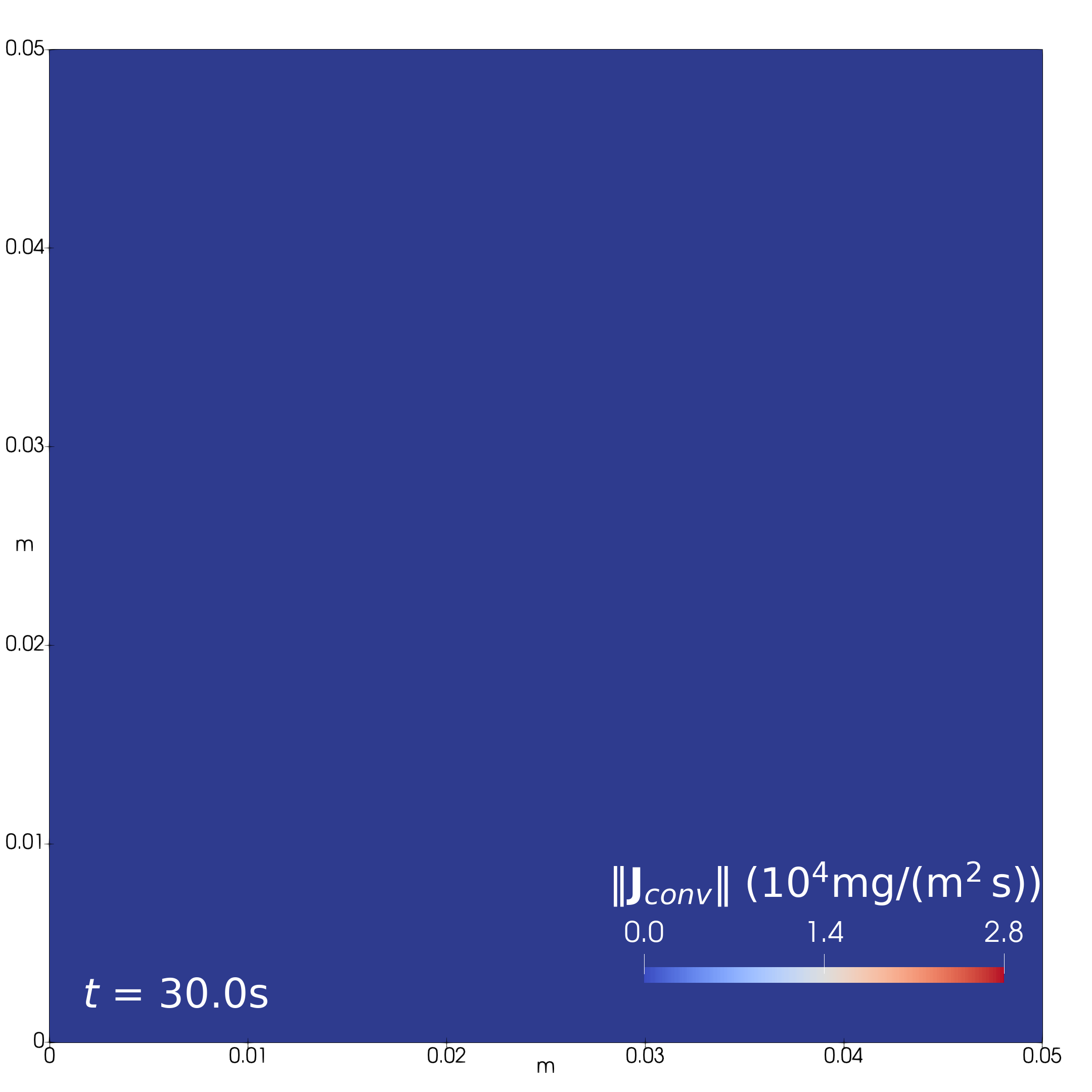}

    \includegraphics[width=0.24\textwidth]{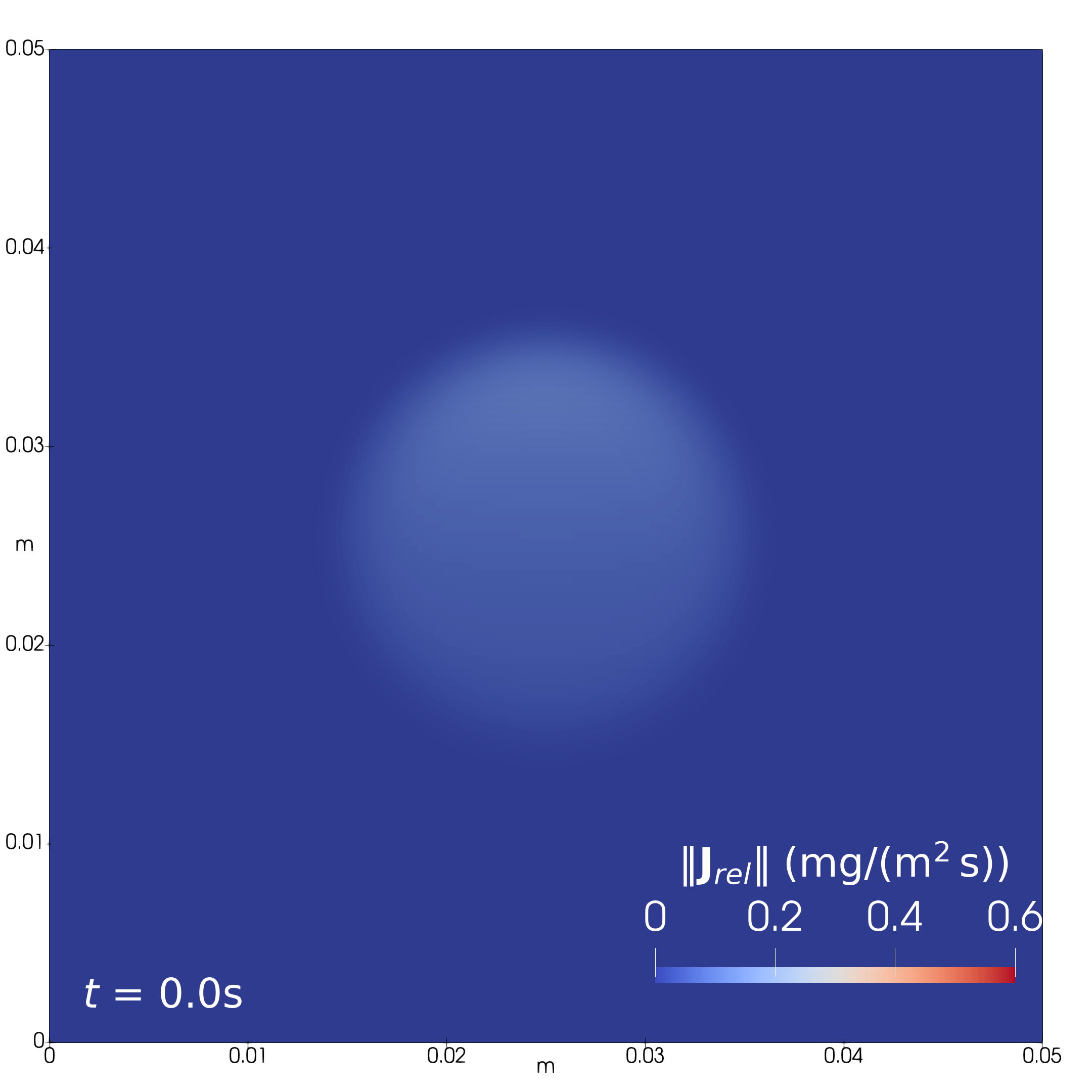}
    \includegraphics[width=0.24\textwidth]{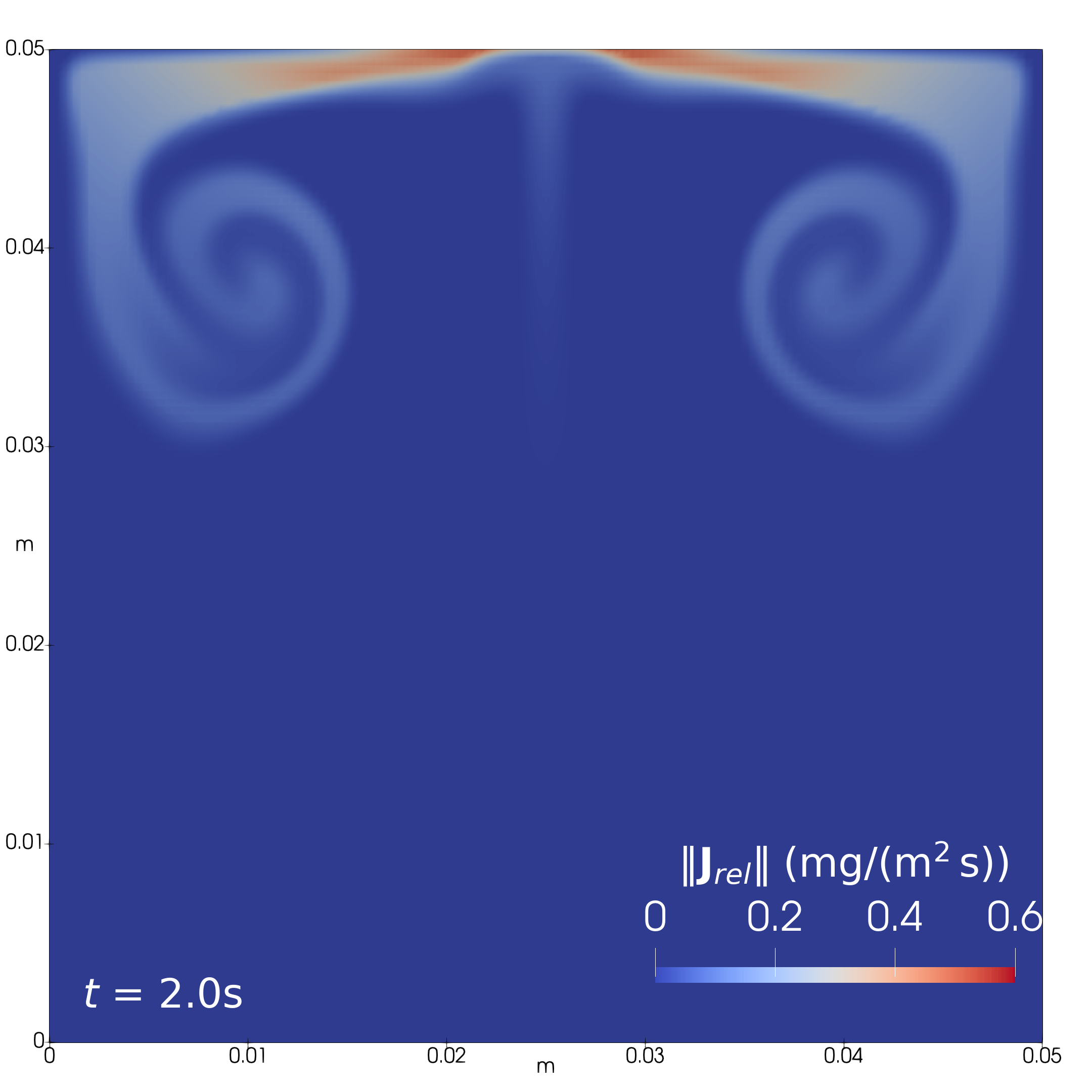}
    \includegraphics[width=0.24\textwidth]{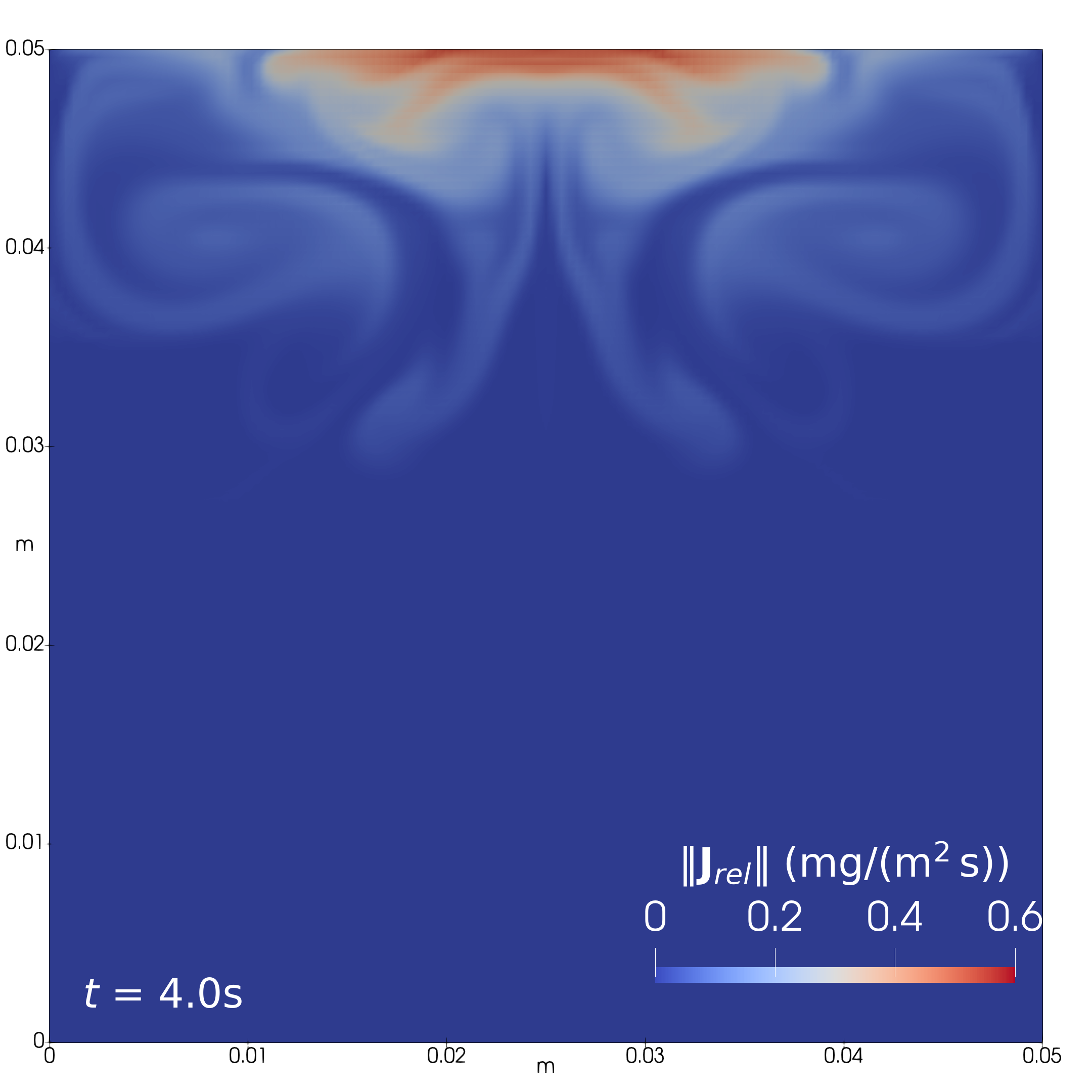}
    \includegraphics[width=0.24\textwidth]{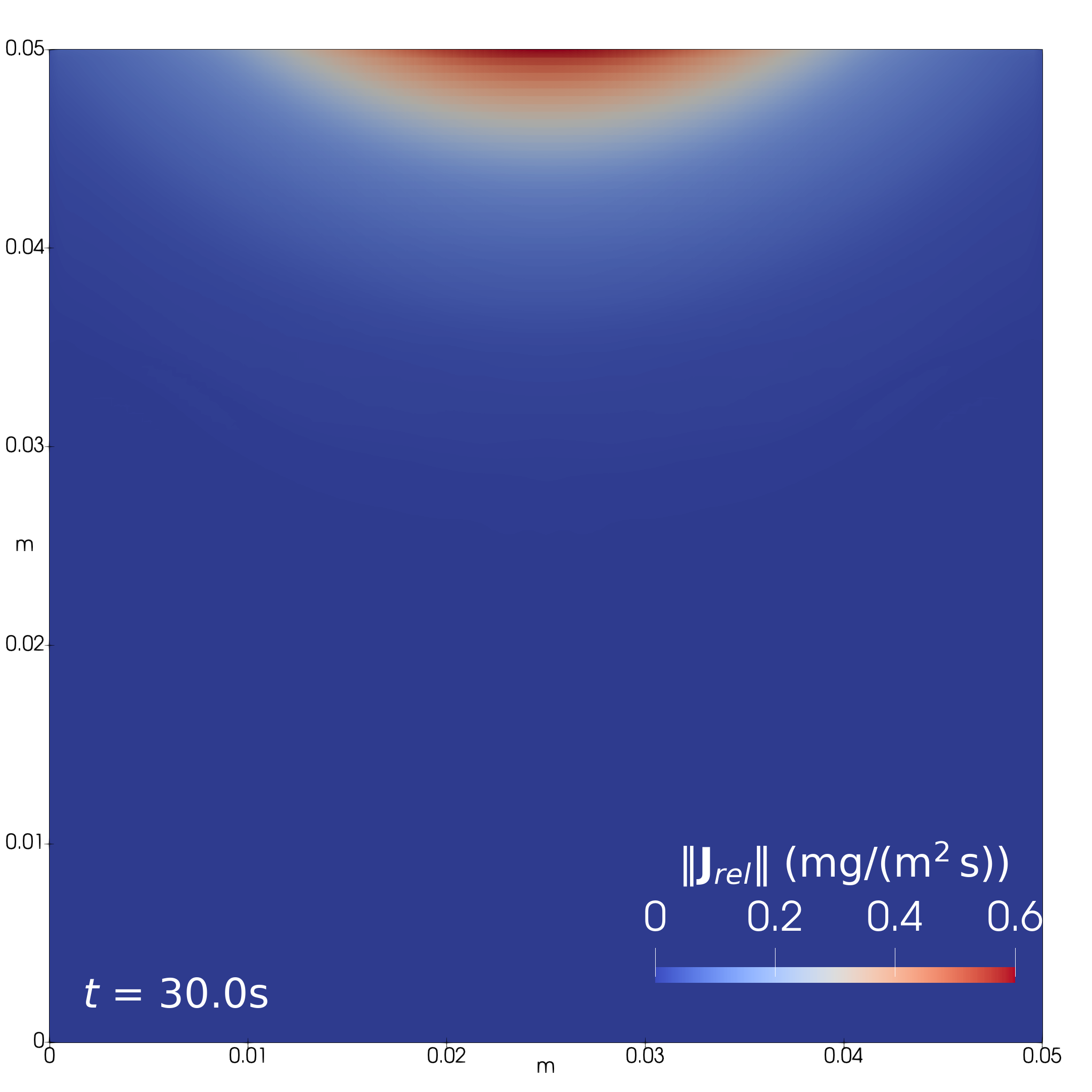}

    \includegraphics[width=0.24\textwidth]{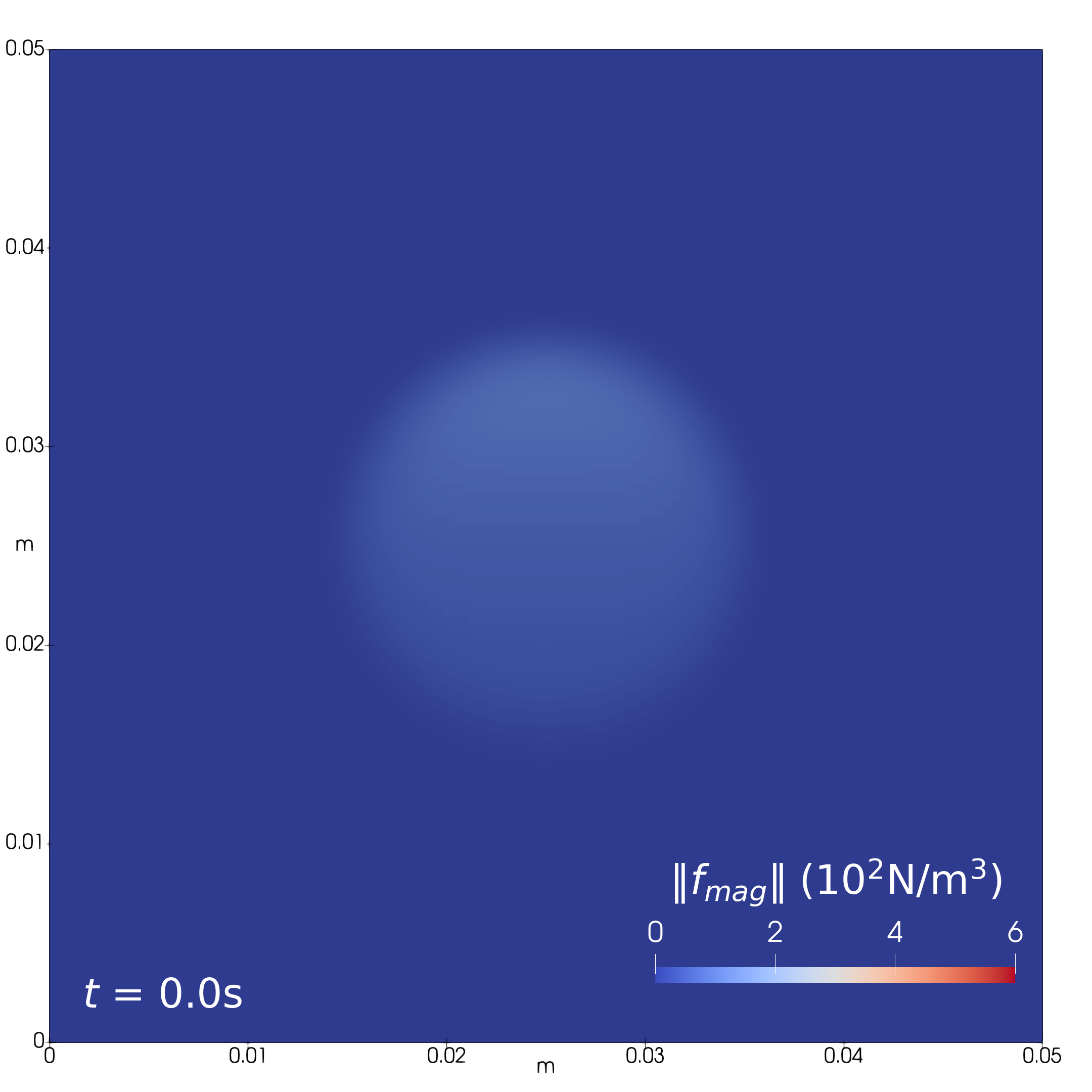}
    \includegraphics[width=0.24\textwidth]{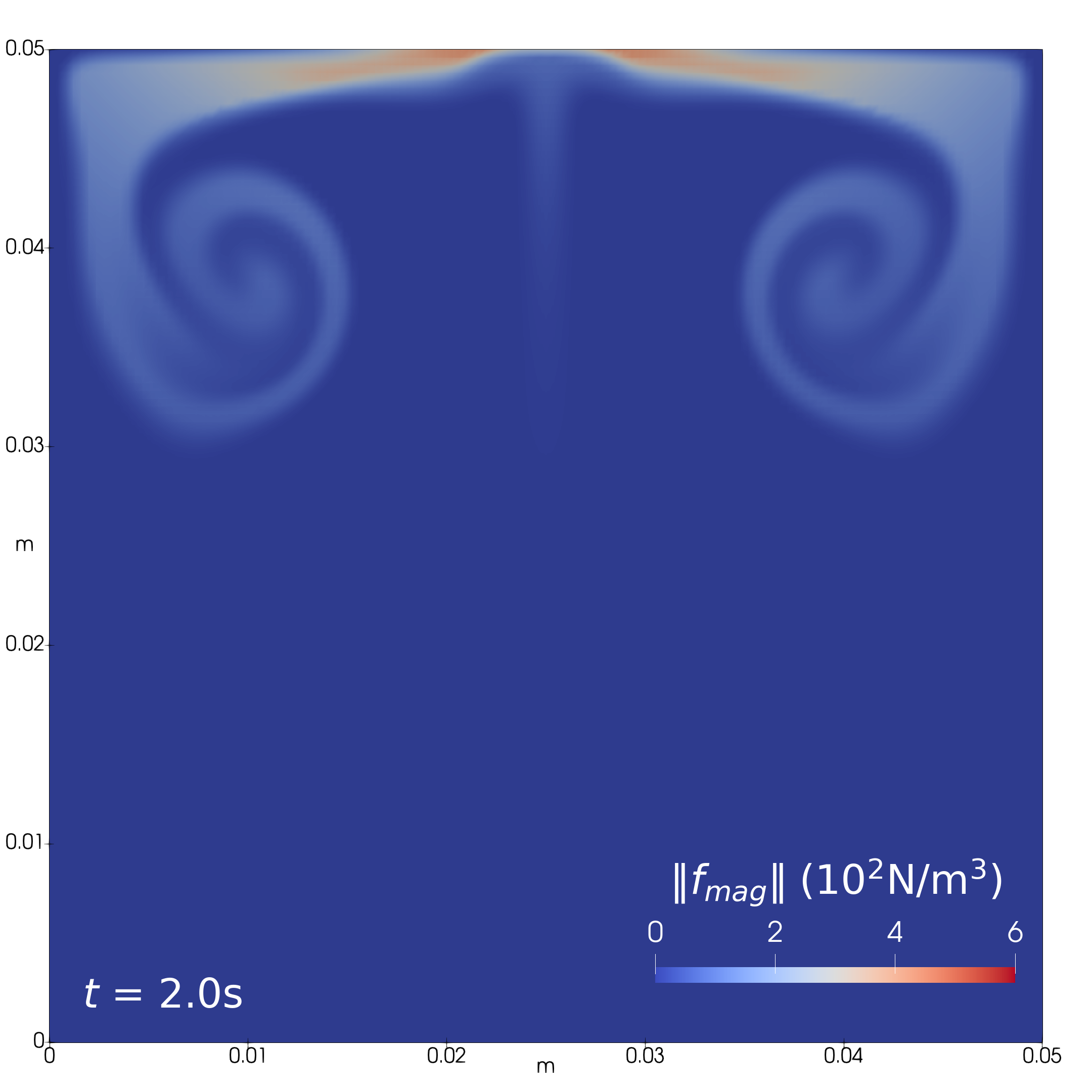}
    \includegraphics[width=0.24\textwidth]{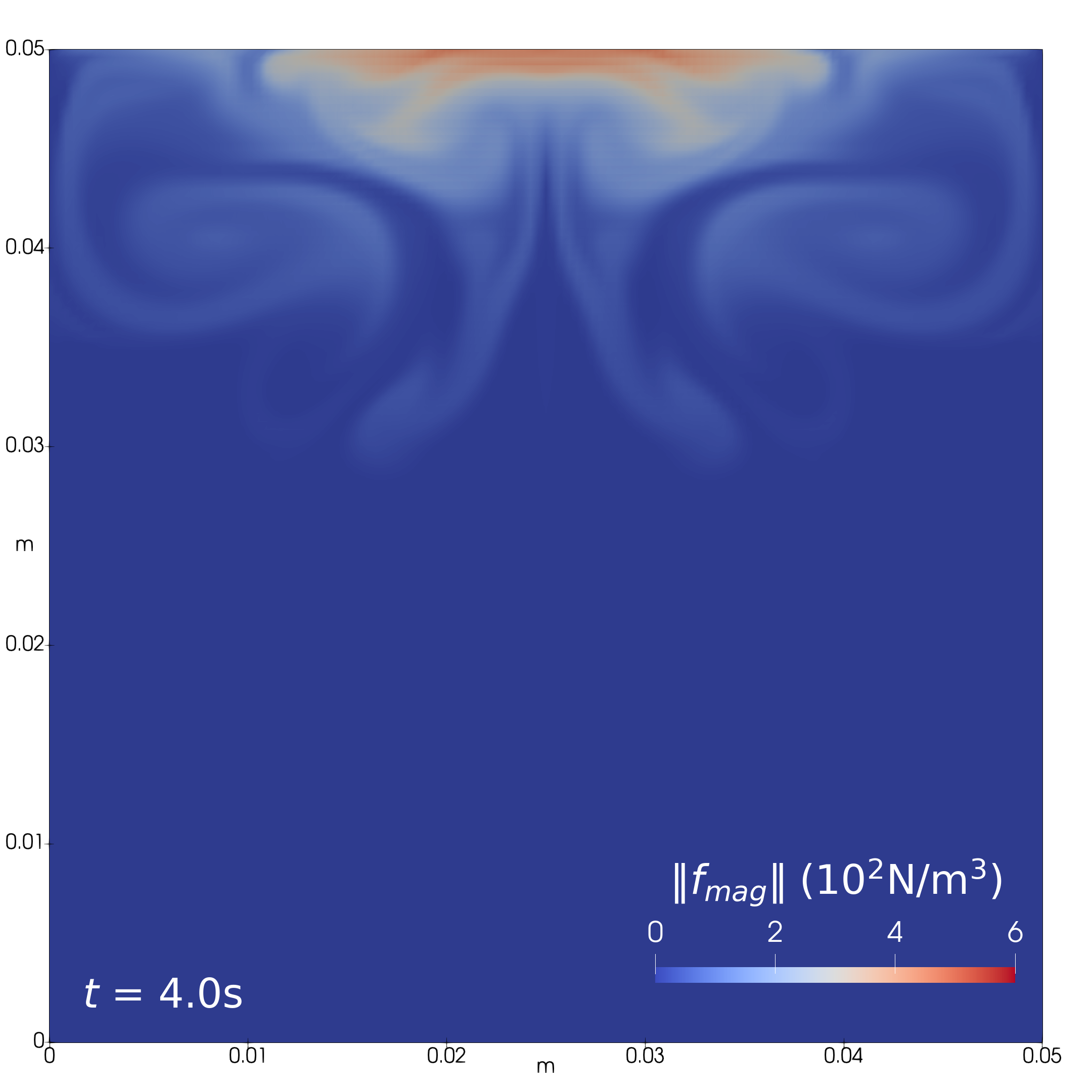}
    \includegraphics[width=0.24\textwidth]{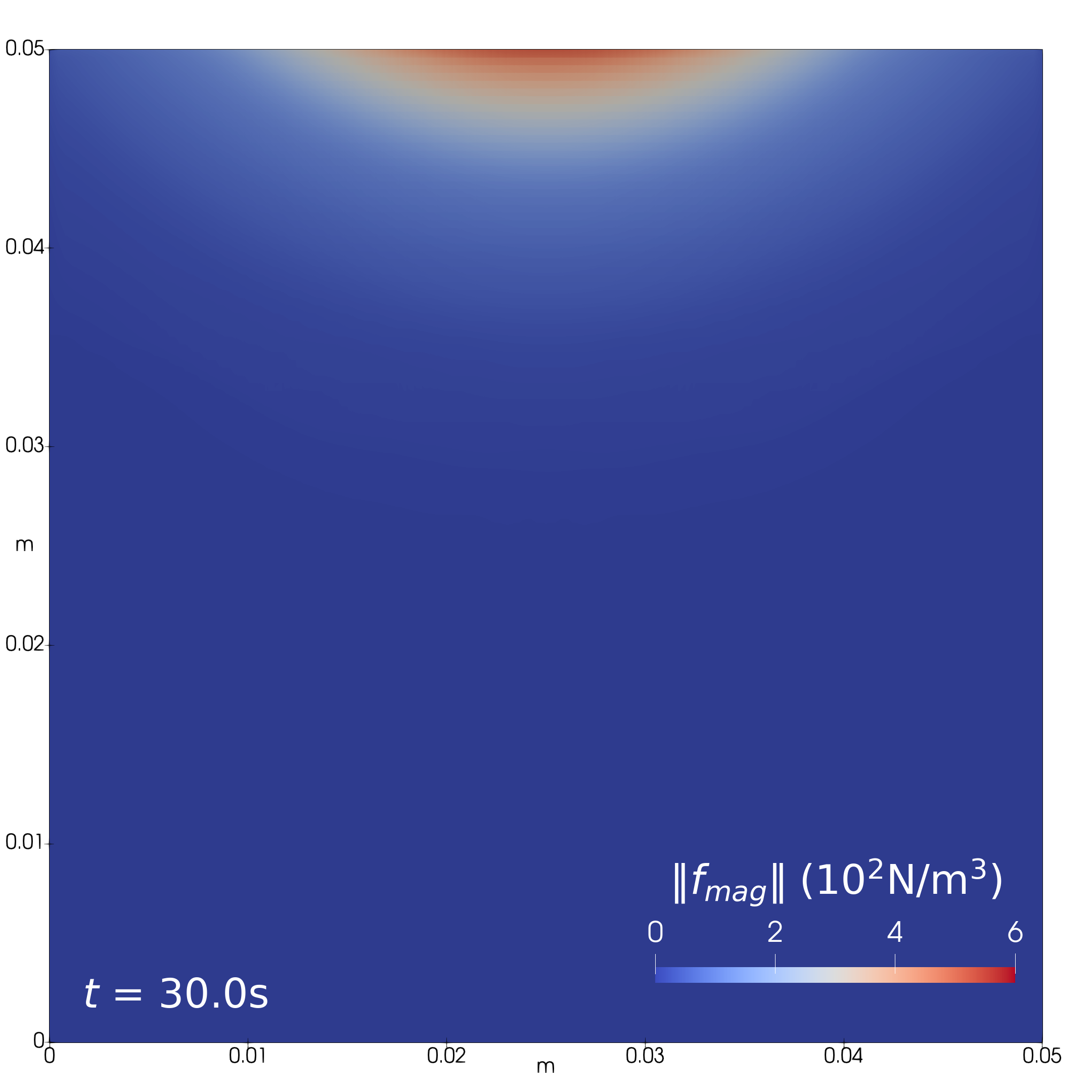}

    \caption{Snapshots at times $t=0\mathrm{s}, 2\mathrm{s}, 4\mathrm{s}, 30\mathrm{s}$ of several physical observables from the simulation of the magnetic attraction of the SPION dispersion. First row: normed relative difference between the computed magnetic field and the external magnetic field; second row: normed relative difference between the gradient of the computed magnetic field and the gradient of the external magnetic field; third row: norm of the convective part of the SPION flux, see Equation \eqref{eq:spion_flux_contributions}; fourth row: norm of the relative part of the SPION flux, consisting of diffusion and magnetophoretic drift, see Equation \eqref{eq:spion_flux_contributions}; fifth row: norm of the magnetic body force exerted by the SPION dispersion onto the mixture, see Equation \eqref{eq:magnetic_body_force}.}
    
    \label{fig:simulation-readout-figures}
\end{figure}

\begin{figure}[tbp]
    \centering
    \begin{subfigure}[t]{0.32\textwidth}
        \centering
        \begin{tikzpicture}[baseline=(current bounding box.south)]
    \begin{axis}[
    title={$\mathcal{E}_{\mathrm{mag},1}(t)$},
    title style={
    yshift=-2mm
    },
    width=\linewidth,
    xlabel={$\mathrm{s}$},
    ylabel={$\mathrm{J\,m^{-1}}$},
    axis y line*=left,
    axis x line*=bottom,
    xmin=0,
    xmax=10,
    enlarge x limits=false,
    grid=major,
    scaled y ticks=false,
    yticklabel={
    \pgfmathparse{\tick*1e1}
    \pgfmathprintnumber[
        std,
        precision=3
    ]{\pgfmathresult}
    },
    extra description/.code={
    \node[
        anchor=south west,
        inner sep=1pt,
        xshift=-5mm,   % weiter nach links
        yshift=-0.5mm    % weiter nach unten
    ] at (rel axis cs:0,1.02)
    {$\times 10^{-1}$};
    },
    ]
    \addplot[
        red,
        very thick,
        no marks,
        %each nth point=10,
    ]
    table[
        x index=0,
        y expr={
            \thisrowno{2}
        },
        col sep=space,
        skip first n=1,
    ]
    {simulationData/FCT_explicit/theta=0.5/ref=6/dt=0.00625/physicalObservables_thinned.txt};
    %\addlegendentry{$E_{\mathrm{mag,1}}(t)$}
    \end{axis}
\end{tikzpicture}
    \end{subfigure}
    \begin{subfigure}[t]{0.32\textwidth}
        \centering
        \begin{tikzpicture}[baseline=(current bounding box.south)]
    \begin{axis}[
    title={$\mathcal{E}_{\mathrm{mag},2}(t)$},
    title style={
    yshift=-2mm
    },
    width=\linewidth,
    xlabel={$\mathrm{s}$},
    ylabel={$\mathrm{J\,m^{-1}}$},
    axis y line*=left,
    axis x line*=bottom,
    xmin=0,
    xmax=10,
    enlarge x limits=false,
    grid=major,
    scaled y ticks=false,
    yticklabel={
    \pgfmathparse{\tick*1e4}
    \pgfmathprintnumber[
        std,
        precision=2
    ]{\pgfmathresult}
    },
    extra description/.code={
    \node[
        anchor=south west,
        inner sep=1pt,
        xshift=-5mm,   % weiter nach links
        yshift=-0.5mm    % weiter nach unten
    ] at (rel axis cs:0,1.02)
    {$\times 10^{-4}$};
    },
    ]
    \addplot[
        red,
        very thick,
        no marks,
        %each nth point=10,
    ]
    table[
        x index=0,
        y expr={
            \thisrowno{3}
        },
        col sep=space,
        skip first n=1,
    ]
    {simulationData/FCT_explicit/theta=0.5/ref=6/dt=0.00625/physicalObservables_thinned.txt};
    %\addlegendentry{$E_{\mathrm{mag,1}}(t)$}
    \end{axis}
\end{tikzpicture}
    \end{subfigure}
    \begin{subfigure}[t]{0.32\textwidth}
        \centering
        \begin{tikzpicture}[baseline=(current bounding box.south)]
    \begin{axis}[
    title={$\mathcal{E}_{\mathrm{mag},3}(t)$},
    title style={
    yshift=-2mm
    },
    width=\linewidth,
    xlabel={$\mathrm{s}$},
    ylabel={$\mathrm{J\,m^{-1}}$},
    axis y line*=left,
    axis x line*=bottom,
    xmin=0,
    xmax=10,
    enlarge x limits=false,
    grid=major,
    scaled y ticks=false,
    yticklabel={
    \pgfmathparse{\tick*1e1}
    \pgfmathprintnumber[
        std,
        precision=3
    ]{\pgfmathresult}
    },
    extra description/.code={
    \node[
        anchor=south west,
        inner sep=1pt,
        xshift=-5mm,   % weiter nach links
        yshift=-0.5mm    % weiter nach unten
    ] at (rel axis cs:0,1.02)
    {$\times 10^{-1}$};
    },
    ]
    \addplot[
        red,
        very thick,
        no marks,
        %each nth point=10,
    ]
    table[
        x index=0,
        y expr={
            \thisrowno{4}
        },
        col sep=space,
        skip first n=1,
    ]
    {simulationData/FCT_explicit/theta=0.5/ref=6/dt=0.00625/physicalObservables_thinned.txt};
    \end{axis}
\end{tikzpicture}
    \end{subfigure}
    \begin{subfigure}[t]{0.32\textwidth}
        \centering
        \begin{tikzpicture}[baseline=(current bounding box.south)]
    \begin{axis}[
        title={$\mathcal{E}_{\mathrm{mix}}(t)$},
        title style={yshift=-2mm},
        width=\linewidth,
        xlabel={$\mathrm{s}$},
        ylabel={$\mathrm{J\,m^{-1}}$},
        axis y line*=left,
        axis x line*=bottom,
        xmin=0,
        xmax=10,
        enlarge x limits=false,
        grid=major,
        scaled y ticks=false,
        yticklabel style={
            /pgf/number format/fixed,
            /pgf/number format/precision=2,
            /pgf/number format/fixed zerofill
        },
        extra description/.code={
            \node[
                anchor=south west,
                inner sep=1pt,
                xshift=-5mm,
                yshift=-0.5mm
            ] at (rel axis cs:0,1.02)
            {$\times 10^{-6}$};
        },
    ]
    \addplot[
        red,
        very thick,
        no marks,
        %each nth point=10,
    ]
    table[
        x index=0,
        y expr={1e6*\thisrowno{1}},
        col sep=space,
        skip first n=1,
    ]
    {simulationData/FCT_explicit/theta=0.5/ref=6/dt=0.00625/physicalObservables_thinned.txt};
    \end{axis}
\end{tikzpicture}
    \end{subfigure}
    \begin{subfigure}[t]{0.32\textwidth}
        \centering
        \begin{tikzpicture}[baseline=(current bounding box.south)]
    \begin{axis}[
        title={$\mathcal{E}_{\mathrm{kin}}(t)$},
        title style={yshift=-2mm},
        width=\linewidth,
        xlabel={$\mathrm{s}$},
        ylabel={$\mathrm{J\,m^{-1}}$},
        axis y line*=left,
        axis x line*=bottom,
        xmin=0,
        xmax=10,
        enlarge x limits=false,
        grid=major,
        scaled y ticks=false,
        yticklabel style={
            /pgf/number format/fixed,
            /pgf/number format/precision=2,
            /pgf/number format/fixed zerofill
        },
        extra description/.code={
            \node[
                anchor=south west,
                inner sep=1pt,
                xshift=-5mm,
                yshift=-0.5mm
            ] at (rel axis cs:0,1.02)
            {$\times 10^{-4}$};
        },
    ]
    \addplot[
        red,
        very thick,
        no marks,
        %each nth point=10,
    ]
    table[
        x index=0,
        y expr={1e4*\thisrowno{5}},
        col sep=space,
        skip first n=1,
    ]
    {simulationData/FCT_explicit/theta=0.5/ref=6/dt=0.00625/physicalObservables_thinned.txt};
    \end{axis}
\end{tikzpicture}
    \end{subfigure}
    \begin{subfigure}[t]{0.32\textwidth}
        \centering
        \begin{tikzpicture}[baseline=(current bounding box.south)]
    \begin{axis}[
        title={$\mathcal{D}_{\mathrm{kin}}(t)$},
        title style={yshift=-2mm},
        width=\linewidth,
        xlabel={$\mathrm{s}$},
        ylabel={$\mathrm{W\,m^{-1}}$},
        axis y line*=left,
        axis x line*=bottom,
        xmin=0,
        xmax=10,
        enlarge x limits=false,
        grid=major,
        scaled y ticks=false,
        yticklabel style={
            /pgf/number format/fixed,
            /pgf/number format/precision=2,
            /pgf/number format/fixed zerofill
        },
        extra description/.code={
            \node[
                anchor=south west,
                inner sep=1pt,
                xshift=-5mm,
                yshift=-0.5mm
            ] at (rel axis cs:0,1.02)
            {$\times 10^{-4}$};
        },
    ]
    \addplot[
        red,
        very thick,
        no marks,
        %each nth point=10,
    ]
    table[
        x index=0,
        y expr={1e4*\thisrowno{6}},
        col sep=space,
        skip first n=1,
    ]
    {simulationData/FCT_explicit/theta=0.5/ref=6/dt=0.00625/physicalObservables_thinned.txt};
    \end{axis}
\end{tikzpicture}
    \end{subfigure}
    \begin{subfigure}[t]{0.32\textwidth}
        \centering
        \begin{tikzpicture}[baseline=(current bounding box.south)]
    \begin{axis}[
        title={$\mathcal{D}_{\mathrm{drag}}(t)$},
        title style={yshift=-2mm},
        width=\linewidth,
        xlabel={$\mathrm{s}$},
        ylabel={$\mathrm{W\,m^{-1}}$},
        axis y line*=left,
        axis x line*=bottom,
        xmin=0,
        xmax=10,
        enlarge x limits=false,
        grid=major,
        scaled y ticks=false,
        yticklabel style={
            /pgf/number format/fixed,
            /pgf/number format/precision=2,
            /pgf/number format/fixed zerofill
        },
        extra description/.code={
            \node[
                anchor=south west,
                inner sep=1pt,
                xshift=-5mm,
                yshift=-0.5mm
            ] at (rel axis cs:0,1.02)
            {$\times 10^{-8}$};
        },
    ]
    \addplot[
        red,
        very thick,
        no marks,
        %each nth point=10,
    ]
    table[
        x index=0,
        y expr={1e8*\thisrowno{7}},
        col sep=space,
        skip first n=1,
    ]
    {simulationData/FCT_explicit/theta=0.5/ref=6/dt=0.00625/physicalObservables_thinned.txt};
    \end{axis}
\end{tikzpicture}
    \end{subfigure}
    \caption{Time series of the energy and dissipation functionals used to derive the model in Section \ref{sec:derivation_constitutive_relations} during the simulation of the magnetic attraction of a SPION dispersion.}
    \label{fig:energie_dissipation_functionals_over_time}
\end{figure}

\begin{figure}[tbp]
    \centering
    \begin{subfigure}[t]{0.46\textwidth}
        \centering
        \begin{tikzpicture}[baseline=(current bounding box.south)]
    % Total energy
    \begin{axis}[
        width=\linewidth,
        xlabel={Time (s)},
        ylabel={$\mathrm{J\,m^{-1}}$},
        axis y line*=left,
        axis x line*=bottom,
        xmin=0,
        xmax=30,
        enlarge x limits=false,
        grid=major,
        legend pos=north east,
        yticklabel style={
            /pgf/number format/fixed,
            /pgf/number format/precision=8
        },
    ]
    \addplot[
        blue,
        very thick,
        no marks,
        %each nth point=10,
    ]
    table[
        x index=0,
        y expr={
            \thisrowno{1}
            + \thisrowno{2}
            + \thisrowno{3}
            + \thisrowno{4}
            + \thisrowno{5}
        },
        col sep=space,
        skip first n=1,
    ]
    {simulationData/FCT_explicit/theta=0.5/ref=6/dt=0.00625/physicalObservables_thinned.txt};

    \addlegendentry{$\mathcal{E}_{\mathrm{tot}}(t)$}
    \addlegendimage{red,very thick}
    \addlegendentry{$\mathcal{D}_{\mathrm{tot}}(t)$}

    \end{axis}

    % ------------------------------------------------------------
    % Total dissipation
    % ------------------------------------------------------------
    \begin{axis}[
    width=\linewidth,
    axis y line*=right,
    axis x line=none,
    xmin=0,
    xmax=30,
    enlarge x limits=false,
    ylabel={$\mathrm{W\,m^{-1}}$},
    ylabel style={
        at={(axis description cs:1.40,0.5)},
        anchor=south,
    },
    xtick=\empty,
    % Display 0.0002 as 2, etc.
    scaled y ticks=false,
    yticklabel={
        \pgfmathparse{\tick*1e4}
        \pgfmathprintnumber{\pgfmathresult}
    },
    % Single common scaling factor
    extra description/.code={
        \node[anchor=south west]
        at (rel axis cs:1,1.02)
        {$\times 10^{-4}$};
    },
    ]

    \addplot[
        red,
        very thick,
        no marks,
        %each nth point=10,
    ]
    table[
        x index=0,
        y expr={
            \thisrowno{6}
            + \thisrowno{7}
        },
        col sep=space,
        skip first n=1,
    ]
    {simulationData/FCT_explicit/theta=0.5/ref=6/dt=0.00625/physicalObservables_thinned.txt};

    \end{axis}

\end{tikzpicture}
        \label{fig:energy-dissipation}
    \end{subfigure}
    \hfill
    \begin{subfigure}[t]{0.46\textwidth}
        \centering
        \begin{tikzpicture}[baseline=(current bounding box.south)]

% axis y line*=right,
%     axis x line=none,
%     xmin=0,
%     xmax=30,
%     enlarge x limits=false,
%     ylabel={$\mathrm{W\,m^{-1}}$},
%     ylabel style={
%         at={(axis description cs:1.40,0.5)},
%         anchor=south,
%     },
    
    \begin{axis}[
        width=\linewidth,
        xmin=0,
        xmax=30,
        enlarge x limits=false,
        xlabel={Time (s)},
        ylabel={$\mathrm{mg\,ml^{-1}}$},
        axis y line*=right,
        axis x line*=bottom,
        grid=major,
        legend columns=2,
        ylabel style={
        at={(axis description cs:1.40,0.5)},
        anchor=south,
    },
    ]
    % u_min
    \addplot[
        blue,
        very thick,
        no marks,
        %each nth point=10,
    ]
    table[
        x index=0,
        y index=8,
        col sep=space,
        skip first n=1,
    ]
    {simulationData/FCT_explicit/theta=0.5/ref=6/dt=0.00625/physicalObservables_thinned.txt};
    \addlegendentry{$\rho_{\mathrm{Fe},\min}(t)$}

    % u_max
    \addplot[
        red,
        very thick,
        no marks,
        %each nth point=10,
    ]
    table[
        x index=0,
        y index=9,
        col sep=space,
        skip first n=1,
    ]
    {simulationData/FCT_explicit/theta=0.5/ref=6/dt=0.00625/physicalObservables_thinned.txt};
    \addlegendentry{$\rho_{\mathrm{Fe},\max}(t)$}

    \end{axis}
\end{tikzpicture}
        \label{fig:minmax}
    \end{subfigure}
    \vspace{-5mm}
    \begin{subfigure}[t]{0.46\textwidth}
        \centering
        \begin{tikzpicture}
    \begin{axis}[
        width=\linewidth,
        xmin=0,
        xmax=30,
        enlarge x limits=false,
        xlabel={Time (s)},
        ylabel={$\mathrm{mg\,m^{-1}}$},
        axis y line*=left,
        axis x line*=bottom,
        grid=major,
        legend pos=south east,
    ]

    \addplot[
        blue,
        very thick,
        no marks,
        %each nth point=10,
    ]
    table[
        x index=0,
        y index=11,
        col sep=space,
        skip first n=1,
    ]{simulationData/FCT_explicit/theta=0.5/ref=6/dt=0.00625/physicalObservables_thinned.txt};

    \addlegendentry{$m_\text{target}(t)$}
    \addlegendimage{red,very thick}
    \addlegendentry{$\mu_2(t)$}    
    \end{axis}

    \begin{axis}[
        width=\linewidth,
        axis y line*=right,
        axis x line=none,
        xmin=0,
        xmax=30,
        enlarge x limits=false,
        ylabel={$\mathrm{m}$},
        ylabel style={
            at={(axis description cs:1.40,0.5)},
            anchor=south,
        },
        xtick=\empty,
        % Display 0.0002 as 2, etc.
        scaled y ticks=false,
        yticklabel={
        \pgfmathparse{\tick*1e2}
        \pgfmathprintnumber{\pgfmathresult}
    },
    % Single common scaling factor
    extra description/.code={
        \node[anchor=south west]
        at (rel axis cs:1,1.02)
        {$\times 10^{-2}$};
    },
    ]
    \addplot[
        red,
        very thick,
        no marks,
        %each nth point=10,
    ]
    table[
        x index=0,
        y index=12,
        col sep=space,
        skip first n=1,
    ]
    {simulationData/FCT_explicit/theta=0.5/ref=6/dt=0.00625/physicalObservables_thinned.txt};

    \end{axis}
\end{tikzpicture}
        \label{fig:target_mass_and_com}
    \end{subfigure}
    \hfill
    \begin{subfigure}[t]{0.46\textwidth}
        \centering
        \begin{tikzpicture}[baseline=(current bounding box.south)]

    \begin{axis}[
        width=\linewidth,
        xmin=0,
        xmax=30,
        enlarge x limits=false,
        xlabel={Time (s)},
        ylabel={$\mathrm{m^2}$},
        axis y line*=right,
        axis x line*=bottom,
        grid=major,
        legend columns=1,
        ylabel style={
            at={(axis description cs:1.40,0.5)},
            anchor=south,
        },
        scaled y ticks=false,
        yticklabel={
            \pgfmathparse{\tick*1e4}
            \pgfmathprintnumber{\pgfmathresult}
        },
        extra description/.code={
            \node[anchor=south west]
            at (rel axis cs:1,1.02)
            {$\times 10^{-4}$};
        },
    ]

    % sigma_1^2
    \addplot[
        blue,
        very thick,
        no marks,
        %each nth point=10,
    ]
    table[
        x index=0,
        y index=13,
        col sep=space,
        skip first n=1,
    ]{simulationData/FCT_explicit/theta=0.5/ref=6/dt=0.00625/physicalObservables_thinned.txt};

    \addlegendentry{$\sigma_1^2(t)$}

    % sigma_2^2
    \addplot[
        red,
        very thick,
        no marks,
        %each nth point=10,
    ]
    table[
        x index=0,
        y index=14,
        col sep=space,
        skip first n=1,
    ]{simulationData/FCT_explicit/theta=0.5/ref=6/dt=0.00625/physicalObservables_thinned.txt};

    \addlegendentry{$\sigma_2^2(t)$}

    \end{axis}

\end{tikzpicture}
        \label{fig:variance}
    \end{subfigure}
    \caption{Physical quantities during the simulation, computed from the discrete solution. Upper-left: total energy and total dissipation, see Equation \eqref{eq:E_tot_D_tot}; upper-right: minimal and maximal iron density in the domain, see Equation \eqref{eq:minmax}; lower-left: iron mass in target region and vertical center of mass of the SPION dispersion, see Equations \eqref{eq:target_mass} and \eqref{eq:center_of_mass}; lower-right: spatial variance of the SPION dispersion in horizontal and vertical direction, see Equation \eqref{eq:spatial_variance}.}
    \label{fig:simulation-summary}
\end{figure}

The simulation in this chapter is performed for the data described in the previous section with $T_\text{end}=30\mathrm{s}$, time step size $\tau = 0.00625\mathrm{s}$ and $\theta=0.5$. The finite element grid is obtained by uniformly refining the initial coarse grid in Figure \ref{fig:domain-target-region} six times, resulting in 66,049+33,282+132,098+16641+16641=264,711 degrees of freedom per time step for the scalar magnetic potential, the magnetic field strength, the velocity, the pressure and the SPION dispersion volume density, respectively. Further refinement of the spatial mesh or reduction of the time step size does not result in appreciable changes in the numerical solution. All physical material parameters, constants, reference quantities and characteristic values used in the simulation can be found in Tables \ref{tab:physical_parameters} and \ref{tab:reference_characteristic_values}, respectively. \par
Figure \ref{fig:simulation-snapshots} shows snapshots of the simulation arranged in chronological order from top left to bottom right. At the beginning, the initially circular SPION-rich region is moving towards the external magnetic field source located above the top edge of the domain, while becoming elongated in the direction of the magnetic field gradient. The SPION dispersion reaches the upper boundary within $1.5$ seconds and spreads laterally, forming an umbrella-like shape whose lateral extensions subsequently develop into two mirror-symmetric vortical structures. Up to this time point, the white velocity arrows, which are scaled according to their magnitude, are most clearly visible. A comparison of the convective flux with the relative flux of the SPIONs due to diffusion and magnetophoretic drift, see Figure \ref{fig:simulation-readout-figures}, row three and four, hints that the observed SPION dynamics during the initial phase of the simulation is dominated by convection with the mixture velocity field $\bfv$. At later times, the vortical structures gradually decay and the SPION dispersion becomes increasingly diffuse, while a persistent accumulation of SPIONs develops along the upper boundary of the domain, directly beneath the magnet. After about 15 seconds, almost all of the vortices have dissipated, and the density of the SPION accumulation at the upper edge decreases relatively smoothly with distance from the magnet. After $30$ seconds, the maximum iron density, at about $0.85 \,\mathrm{mg/ml}$, is still below the initial maximum density of $1 \,\mathrm{mg/ml}$. The streamlines of the magnetic field, visible in the images as black lines, did not change in any way visible to the naked eye over the course of the experiment. 
\par
Figure \ref{fig:simulation-readout-figures} allows for a more detailed look into the alterations of the magnetic field due to the SPIONs. In the first row of figures, the relative difference of the computed magnetic field, given by $\mathbf{b} = \mu_0(\bfh + \mathbf{m})$, to the external magnetic field $\mathbf{b}_\textrm{ext}$ is visualized. The spatial deviations are clearly linked to the support of the SPION dispersion volume density and a maximum difference of $1.2\%$ can be observed during the initial phase of the simulation. For later time points, the difference is smaller. A similar analysis is presented for the gradient of the magnetic field, $\nabla \mathbf{b}$ in the second row of figures. The deviation is significantly higher, up to $15\%$ at $t=2\textrm{s}$, even though it occurs only locally in regions with larger gradients of the SPION dispersion. Here, too, the magnitude of the deviation is significantly smaller for later times. The third and forth row of figures show the magnitude of the fluxes of the SPION dispersion, broken down into convective flux and relative flux due to diffusion and magnetophoretic drift, i.e., 
\begin{equation} \label{eq:spion_flux_contributions}
    \mathbf{J}_\mathrm{Fe, conv} = \rho_\text{Fe} \bfv, \quad \mathbf{J}_\mathrm{Fe,rel} = -D \nabla \rho_\mathrm{Fe} + \kappa (\tilde{\rho}_\mathrm{Fe,max} - \rho_\mathrm{Fe})\rho_\mathrm{Fe} M_\mathrm{s} \mathcal{L}(\xi|\bfh|) \frac{\bfh}{|\bfh|}.
\end{equation}
Since the mixture is at rest at the beginning of the simulation, the convective flux is initially zero, whereas a nonzero relative flux is visible at $t=0\mathrm{s}$. However, during the initial phase of the simulation ($0-5\mathrm{s}$), the SPION flux is clearly dominated by the convective contribution, lying approximately up to four orders of magnitude above the relative flux. This confirms the conclusion that the particle dynamics is initially dominated by convection with the velocity field $\bfv$ excited by the SPIONs, and that magnetophoretic drift only leads to the slow accumulation of SPIONs at the magnet over longer time intervals. In this context, the fifth row shows the magnitude of the magnetic volume force
\begin{equation} \label{eq:magnetic_body_force}
    \mathbf{f}_\text{mag} = \mu_0 M_\textrm{s} u_\mathrm{p} \frac{\mathcal{L}(\xi|\bfh|)}{|\bfh|} (\bfh \cdot \nabla) \bfh
\end{equation}
exerted by the nanoparticle phase on the mixture, which drives the convective processes.
\par
Some of the observations described above can also be understood by examining the temporal evolution of the energy and dissipation functionals shown in Figure \ref{fig:energie_dissipation_functionals_over_time}. Note that we display only the evolution during the initial 10 seconds of the numerical experiment for a better resolution of the dynamics. The curves remain almost constant during the interval $10\textrm{s}-30\textrm{s}$. Notably, the kinetic energy $\mathcal{E}_\text{kin}$ initially increases sharply, coinciding with the movement of the SPIONs towards the magnet, before gradually decaying towards zero. After approximately six seconds, the fluid is nearly at rest. The kinetic dissipation, $\mathcal{D}_\text{kin}$, on the other hand, exhibits two peaks. The first peak can likely be attributed to the direct attraction of the SPIONs during the 0–1.5-second interval and the associated localized energy transfer to the fluid, while the second peak arises from larger velocity gradients during the formation of the vortical structures. The magnetic energy of the SPION dispersion, $\mathcal{E}_\textrm{mag,2}$, initially decreases, which is consistent with the SPIONs' energetically favored position near the magnet. The minimum is reached when the umbrella-like shape is visible and a maximal amount of SPIONs is in the target region, see also Equation \eqref{eq:target_mass} and Figure \ref{fig:simulation-summary} (lower-left). Once the vortical structures start to form, transporting a fraction of the SPIONs away from the magnet, the energy is again increasing and on the long run levels off at a low level. Even though drag dissipation, $\mathcal{D}_\textrm{drag}$ is four orders of magnitude lower than kinetic dissipation, it is noteworthy that, unlike the kinetic dissipation, it does not return to its initial level of zero within 30 seconds, but instead remains constant at an elevated level. This suggests that the diffusion and magnetophoretic drift of the nanoparticles continue to proceed slowly, even though the mixture appears macroscopically at rest. We expect the particles to continue moving toward the magnet. However, since the force exerted by the magnetic field on the particles is relatively small, this secondary accumulation of SPIONs is difficult to observe within the short time span of the numerical experiment. A hint is given in Figure \ref{fig:simulation-summary} (upper-right) where it can be seen that after an initial decrease of the maximal iron density, it slowly starts to increase again after about 5 seconds once the mixing associated with the convective transport of the SPIONs is over.
\par
Following these observations, we present additional physical readouts from the simulation in Figure \ref{fig:simulation-summary}. The displayed quantities are, for all $t \in I$:
\begin{enumerate}
    \item The total energy and total dissipation of the system
    \begin{equation}\label{eq:E_tot_D_tot}
        \begin{aligned}
        \mathcal{E}_\text{tot}(t) &= \mathcal{E}_\text{mag,1}(t) + \mathcal{E}_\text{mag,2}(t) + \mathcal{E}_\text{mag,3}(t) + \mathcal{E}_\text{kin}(t) + \mathcal{E}_\text{mix}(t),\\
        \mathcal{D}_\text{tot}(t) &= \mathcal{D}_\text{kin}(t) + \mathcal{D}_\text{drag}(t).
        \end{aligned}
    \end{equation}
    \item The minimal and maximal iron density in the domain:
    \begin{equation}\label{eq:minmax}
        \rho_\text{Fe,min}(t) = \min_{x \in \Om} \rho_\text{Fe}(t,x), \quad\rho_\text{Fe,max}(t) = \max_{x \in \Om} \rho_\text{Fe}(t,x).
    \end{equation}
    \item The iron mass in target region $\Omega_\text{target} = [0\textrm{m}, 0.05\textrm{m}]\times[0.0375\textrm{m}, 0.05\textrm{m}]$:
    \begin{equation}\label{eq:target_mass}
        m_\text{target}(t) = \int_{\Om_\text{target}} \rho_\text{Fe}(t,x) \; \textrm{d}x    
    \end{equation}
    \item The coordinates of the center of mass of the SPION dispersion:
    \begin{equation}\label{eq:center_of_mass}
        \mu_i(t) = \frac{
            \int_{\Om} x_i\,\rho_{\text{Fe}}(t,x)\,\mathrm{d}x
            }{
            \int_{\Om} \rho_{\text{Fe}}(t,x)\,\mathrm{d}x
        }, \qquad i \in \{1,2\}.
    \end{equation}
    \item The spatial variance of the SPION dispersion:
    \begin{equation}\label{eq:spatial_variance}
        \sigma_i^2(t)
        =
        \frac{
            \displaystyle
            \int_{\Om}
            \bigl(x_i-\mu_i(t)\bigr)^2
            \rho_{\text{Fe}}(t,x)\,\mathrm{d}x
        }{
            \displaystyle
            \int_{\Om}
            \rho_{\text{Fe}}(t,x)\,\mathrm{d}x
        },
        \qquad i\in\{1,2\}.
    \end{equation}
\end{enumerate}
We note that, as in Figure \ref{fig:energie_dissipation_functionals_over_time}, extensive quantities are reported per unit length in the out-of-plane direction, namely in $\mathrm{J \; m^{-1}}$ for energy, $\mathrm{W \; m^{-1}}$ for dissipation, and $\mathrm{mg \;m^{-1}}$ for iron mass. 
\par
As shown in the upper-left  panel of Figure \ref{fig:simulation-summary}, the total energy of the system decreases throughout the simulation and approaches a nearly constant value after approximately 15 seconds. This behavior is consistent with the derived energy-dissipation law \eqref{eq:energy_balance} and, consequently, with the second law of thermodynamics. In particular, since the external magnetic field is time-independent, no external power is supplied through temporal variations of the magnetic field, and the total energy is non-increasing. The temporal evolution of the total dissipation is determined primarily by the kinetic dissipation, since drag dissipation is several orders of magnitude smaller (see Figure \ref{fig:energie_dissipation_functionals_over_time}). The minimal and maximal iron density of the SPION dispersion is depicted in the upper-right panel of Figure \ref{fig:simulation-summary}. Here, the density computed using the FCT method discussed in Section \ref{sec:discretization:transportStabilization} remains within physical bounds, that is, positive, up to machine precision. We stress that the observed maximal iron density of $1\, \textrm{mg/ml}$ is not the physical upper bound of the dispersion and higher values are in principle possible. Here, however, the initial maximum density first decreases due to mixing and then slowly increases again. We attribute this second effect to the slow magnetophoretic drift and expect that, in the long run, the maximum density near the magnet will reach values above $1\, \mathrm{mg/ml}$. However, the timescale of this process exceeds the 30-second simulation time used here. The lower-left panel shows the amount of iron in the target region $\Om_\text{target}$, defined as the upper quarter of the simulation domain $\Omega$. At the peak, which occurs at approximately $t = 1.44\textrm{s}$, 98.24\% of the SPION mass is concentrated in $\Om_\text{target}$. The value then levels off slightly below that and remains approximately constant at 83.1\%. The shape of the vertical center of mass curve $\mu_2$ is very similar to the target mass dynamics.
The spatial variance curves in the lower-right panel confirm the observation that SPION dispersion spreads more widely in the horizontal direction, whereas it tends to remain concentrated in the vertical direction.
\par
Finally we investigate the energy-dissipation balance from Equation \eqref{eq:energy_balance}  during the simulation. In the present simulation scenario, the externally imposed magnetic field $\mathbf b_{\mathrm{ext}}$ is time-independent, and hence it does not contribute an external power input to the system. Consequently, on the analytic level, the energy-dissipation balance reduces to
\begin{equation} \label{eq:energy_dissipation_with_zero_rhs}
    \mathcal{E}'_\text{kin}   + \mathcal{E}'_\text{mix} + \mathcal{E}'_\text{mag,1} + \mathcal{E}'_\text{mag,2} + \mathcal{E}'_\text{mag,3}+ 2\mathcal{D}_\text{kin} +  2\mathcal{D}_\text{drag} = 0.
\end{equation}
However, it cannot be expected that the balance directly transfers to the simulation since the cancellations used for its derivation on the continuous level do not carry over to the discrete level. The extent to which this energy balance is fulfilled by the numerical scheme can be observed in Figure \ref{fig:energy_balance}. There, we report the quantities 
\begin{equation} \label{eq:discrete_energy_dissipation}
    \begin{aligned}
        \delta_t \mathcal{E}^n &= \frac{1}{\tau}\left(\mathcal{E}^n_\textrm{tot} - \mathcal{E}^{n-1}_\textrm{tot} \right),\\
        \mathcal{D}^n_\theta &= \theta \mathcal{D}^n_\textrm{tot} + (1 - \theta) \mathcal{D}^n_\textrm{tot},\\
        \mathcal{R}^n_\theta & = \delta_t \mathcal{E}^n_\textrm{tot} + 2 \mathcal{D}^n_\theta,
    \end{aligned}
\end{equation}
for the time steps $n = 1,...,N$. Here, $\mathcal{E}^n_\textrm{tot}$ and $\mathcal{D}^n_\textrm{tot}$ denote the total energy and dissipation computed from the discrete solutions of the numerical scheme at time step $n$. Consequently, $\mathcal{R}^n_\theta$ is the discrete residual of the energy-dissipation balance \eqref{eq:energy_dissipation_with_zero_rhs}. The plot shows that the difference quotient of the total energy qualitatively follows twice the negative total dissipation, and that the residual is close to zero, particularly in the period after three seconds. Nevertheless, it is noticeable that there are two significant deviations at $t=1.2\textrm{s}$ and $t=2.4\textrm{s}$. Compared to the absolute value of $2\mathcal{D}^n_\theta$, the relative residual at these points is 10.89\% and 3.26\%, respectively. In the next section, we will numerically investigate the extent to which the deviation of the residual from zero can be controlled by the time step size $\tau$.

\begin{figure}
    \centering
    \begin{tikzpicture}[baseline=(current bounding box.south)]
    \begin{axis}[
        height=0.4\textwidth,
        width=0.8\linewidth,
        xlabel={$\mathrm{s}$},
        ylabel={$\mathrm{W\,m^{-1}}$},
        axis y line*=left,
        axis x line*=bottom,
        xmin=0,
        xmax=10,
        enlarge x limits=false,
        grid=major,
        scaled y ticks=false,
        yticklabel style={
            /pgf/number format/fixed,
            /pgf/number format/precision=1,
            /pgf/number format/fixed zerofill
        },
        legend pos=south east,
        % legend style={
        %     draw=none,
        %     fill=none
        % },
        extra description/.code={
            \node[
                anchor=south west,
                inner sep=1pt,
                xshift=-5mm,
                yshift=-0.5mm
            ] at (rel axis cs:0,1.02)
            {$\times 10^{-4}$};
        },
    ]

    % ------------------------------------------------------------
    % Residual
    % ------------------------------------------------------------
    \addplot[
        blue,
        very thick,
        no marks,
        %each nth point=10,
    ]
    table[
        x index=0,
        y expr={1e4*\thisrowno{1}},
        col sep=space,
        comment chars={\#},
    ]
    {simulationData/FCT_explicit/theta=0.5/ref=6/dt=0.00625/energyBalanceResidual.txt};
    \addlegendentry{$\mathcal{R}^n_\theta$}

    % ------------------------------------------------------------
    % d_t E
    % ------------------------------------------------------------
    \addplot[
        red,
        very thick,
        no marks,
        %each nth point=10,
    ]
    table[
        x index=0,
        y expr={1e4*\thisrowno{2}},
        col sep=space,
        comment chars={\#},
    ]
    {simulationData/FCT_explicit/theta=0.5/ref=6/dt=0.00625/energyBalanceResidual.txt};
    \addlegendentry{$\delta_t\mathcal{E}^n_\textrm{tot}$}

    % ------------------------------------------------------------
    % -2 D_theta
    % ------------------------------------------------------------
    \addplot[
        black,
        very thick,
        dashed,
        no marks,
        %each nth point=10,
    ]
    table[
        x index=0,
        y expr={1e4*\thisrowno{3}},
        col sep=space,
        comment chars={\#},
    ]
    {simulationData/FCT_explicit/convergenceStudy/energyBalanceResidual_dt=0.00625s_thinned.txt};
    \addlegendentry{$-2\mathcal{D}^n_\theta$}

    \end{axis}
\end{tikzpicture}
    \caption{Time series of the difference quotient of the total energy (red), the double negative discrete total dissipation (black, dashed) and the discrete residual of the energy-dissipation balance (blue) during the simulation, see also Equation \eqref{eq:discrete_energy_dissipation}.}
    \label{fig:energy_balance}
\end{figure}

\subsection{Temporal Convergence Study of the Energy and Dissipation Functionals and the Energy-Dissipation Balance}
\label{sec:convergence_study}
Here, we perform a temporal convergence study of the energy and dissipation functionals. Moreover, we investigate the convergence of the energy-dissipation residual
\begin{equation}
    \mathcal{R}_\theta(t) = \mathcal{E}'_\mathrm{tot}(t) + 2\mathcal{D}_\textrm{tot}(t)
\end{equation}
to zero with decreasing time step size. We distinguish the cases $\theta=0.5$ and $\theta=1$ for the numerical scheme presented in Section \ref{sec:discretization}.
To this end, we repeat the numerical experiment presented in the previous sections with $T_\text{end} = 10\mathrm{s}$ for the step sizes $\tau_i = 0.05 \cdot 2^{-i}\mathrm{s}$, $i = 0,1,2,3,4$ for $\theta \in \{0.5, 1\}$. To quantify the numerical convergence orders, we compute the errors
\begin{equation}
    e^\mathcal{Q}_{i} = \| \mathcal{Q}_{i-1} - \mathcal{Q}_{i}\|_{L^2(I)}, \qquad i = 1,2,3,4,
\end{equation}
for the energy and dissipation functionals $\mathcal{Q} \in \{\mathcal{E}_\text{mag,1}, \mathcal{E}_\text{mag,2}, \mathcal{E}_\text{mag,3}, \mathcal{E}_\text{mix}, \mathcal{E}_\text{kin}, \mathcal{D}_\text{drag}, \mathcal{D}_\text{kin} \}$. The discrete time series are compared by evaluating the corresponding quantities at the common time points of successive temporal refinement levels.
To evaluate the convergence of the energy-dissipation residual to zero, we consider
\begin{equation}
    e^{\mathcal{R}_\theta}_i = \|\mathcal{R}_{\theta,i} \|_{L^2(I)}, \qquad i = 0,1,2,3,4.
\end{equation}
Index $i$ in $\mathcal{Q}_i$ and $\mathcal{R}_{\theta,i}$ indicates the quantity being computed using time step size $\tau_i$. The resulting errors and estimated orders of convergence are presented in Table \ref{tab:energy-convergence}. Here, the nondimensional quantities are directly used for the analysis. The presented values indicate convergence of all quantities considered. Generally, the scheme with $\theta=0.5$ exhibits smaller differences between successive temporal refinement levels and higher experimental convergence orders than the scheme with $\theta=1$. For sufficiently small time step sizes, the estimated convergence orders are predominantly close to one for $\theta=0.5$, while remaining below one for $\theta=1$. In particular, the norm of the energy-dissipation residuum decreases consistently with decreasing time step size, with estimated orders of convergence ranging between $1.00 - 1.12$ for $\theta=0.5$ and $0.61-0.95$ for $\theta=1$. These results indicate that the physical energy-dissipation balance is increasingly well fulfilled by the discrete solutions from the numerical scheme as the time steps size is reduced, with a better convergence behavior observed for $\theta=0.5$. However, optimal convergence orders of classical one-step $\theta$ methods (order two for $\theta=0.5$, order one for $\theta\neq0.5$) are not reached. The reduced convergence orders may be attributed to lower-order splitting and coupling errors introduced by the semi-implicit treatment of the governing equations, as well as the explicit FCT stabilization of the SPION volume fraction.

\begin{table}[tbp]
    \centering
    \caption{Numerical convergence study for the energy and dissipation functional as well as the energy-dissipation residual.}
    \label{tab:energy-convergence}

    \begin{subtable}[t]{0.48\textwidth}
        \centering
        \caption{$\mathcal{E}_{\mathrm{mag,1}}$}
        \pgfplotstableread[
    col sep=tab,
    comment chars={\#}
]{simulationData/FCT_explicit/convergenceStudy/convergence_E_mag_1_ref=6_T=10s.txt}\EmagOne

{
\setlength{\tabcolsep}{3pt}

\pgfplotstabletypeset[
    columns={i,err_cn,eoc_cn,err_ie,eoc_ie},
    columns/i/.style={
        column name={$i$},
        string type
    },
    columns/err_cn/.style={
        column name={$e^{\mathcal{E}_\text{mag,1}}_{i}$},
        string type
    },
    columns/eoc_cn/.style={
        column name={EOC},
        string type
    },
    columns/err_ie/.style={
        column name={$e^{\mathcal{E}_\text{mag,1}}_{i}$},
        string type
    },
    columns/eoc_ie/.style={
        column name={EOC},
        string type
    },
    every head row/.style={
        before row={
            \toprule
            & \multicolumn{2}{c}{$\theta=0.5$}
            & \multicolumn{2}{c}{$\theta=1$} \\
            \cmidrule(lr){2-3}
            \cmidrule(lr){4-5}
        },
        after row=\midrule
    },
    every last row/.style={
        after row=\bottomrule
    }
]{\EmagOne}
}
    \end{subtable}
    \hfill
    \begin{subtable}[t]{0.48\textwidth}
        \centering
        \caption{$\mathcal{E}_{\mathrm{mag,2}}$}
        \pgfplotstableread[
    col sep=tab,
    comment chars={\#}
]{simulationData/FCT_explicit/convergenceStudy/convergence_E_mag_2_ref=6_T=10s.txt}\EmagTwo

{
\setlength{\tabcolsep}{3pt}

\pgfplotstabletypeset[
    columns={i,err_cn,eoc_cn,err_ie,eoc_ie},
    columns/i/.style={
        column name={$i$},
        string type
    },
    columns/err_cn/.style={
        column name={$e^{\mathcal{E}_\text{mag,2}}_{i}$},
        string type
    },
    columns/eoc_cn/.style={
        column name={EOC},
        string type
    },
    columns/err_ie/.style={
        column name={$e^{\mathcal{E}_\text{mag,2}}_{i}$},
        string type
    },
    columns/eoc_ie/.style={
        column name={EOC},
        string type
    },
    every head row/.style={
        before row={
            \toprule
            & \multicolumn{2}{c}{$\theta=0.5$}
            & \multicolumn{2}{c}{$\theta=1$} \\
            \cmidrule(lr){2-3}
            \cmidrule(lr){4-5}
        },
        after row=\midrule
    },
    every last row/.style={
        after row=\bottomrule
    }
]{\EmagTwo}
}
    \end{subtable}
        
    \vspace{1em}

    \begin{subtable}[t]{0.48\textwidth}
        \centering
        \caption{$\mathcal{E}_{\mathrm{mag,3}}$}
        \pgfplotstableread[
    col sep=tab,
    comment chars={\#}
]{simulationData/FCT_explicit/convergenceStudy/convergence_E_mag_3_ref=6_T=10s.txt}\EmagThree

{
\setlength{\tabcolsep}{3pt}

\pgfplotstabletypeset[
    columns={i,err_cn,eoc_cn,err_ie,eoc_ie},
    columns/i/.style={
        column name={$i$},
        string type
    },
    columns/err_cn/.style={
        column name={$e^{\mathcal{E}_\text{mag,3}}_{i}$},
        string type
    },
    columns/eoc_cn/.style={
        column name={EOC},
        string type
    },
    columns/err_ie/.style={
        column name={$e^{\mathcal{E}_\text{mag,3}}_{i}$},
        string type
    },
    columns/eoc_ie/.style={
        column name={EOC},
        string type
    },
    every head row/.style={
        before row={
            \toprule
            & \multicolumn{2}{c}{$\theta=0.5$}
            & \multicolumn{2}{c}{$\theta=1$} \\
            \cmidrule(lr){2-3}
            \cmidrule(lr){4-5}
        },
        after row=\midrule
    },
    every last row/.style={
        after row=\bottomrule
    }
]{\EmagThree}
}
    \end{subtable}
    \hfill
    \begin{subtable}[t]{0.48\textwidth}
        \centering
        \caption{$\mathcal{E}_{\mathrm{mix}}$}
        \pgfplotstableread[
    col sep=tab,
    comment chars={\#}
]{simulationData/FCT_explicit/convergenceStudy/convergence_E_mix_ref=6_T=10s.txt}\Emix

{
\setlength{\tabcolsep}{3pt}

\pgfplotstabletypeset[
    columns={i,err_cn,eoc_cn,err_ie,eoc_ie},
    columns/i/.style={
        column name={$i$},
        string type
    },
    columns/err_cn/.style={
        column name={$e^{\mathcal{E}_\text{mix}}_{i}$},
        string type
    },
    columns/eoc_cn/.style={
        column name={EOC},
        string type
    },
    columns/err_ie/.style={
        column name={$e^{\mathcal{E}_\text{mix}}_{i}$},
        string type
    },
    columns/eoc_ie/.style={
        column name={EOC},
        string type
    },
    every head row/.style={
        before row={
            \toprule
            & \multicolumn{2}{c}{$\theta=0.5$}
            & \multicolumn{2}{c}{$\theta=1$} \\
            \cmidrule(lr){2-3}
            \cmidrule(lr){4-5}
        },
        after row=\midrule
    },
    every last row/.style={
        after row=\bottomrule
    }
]{\Emix}
}
    \end{subtable}

    \vspace{1em}

    \begin{subtable}[t]{0.48\textwidth}
        \centering
        \caption{$\mathcal{E}_{\mathrm{kin}}$}
        \pgfplotstableread[
    col sep=tab,
    comment chars={\#}
]{simulationData/FCT_explicit/convergenceStudy/convergence_E_kin_ref=6_T=10s.txt}\Ekin

{
\setlength{\tabcolsep}{3pt}

\pgfplotstabletypeset[
    columns={i,err_cn,eoc_cn,err_ie,eoc_ie},
    columns/i/.style={
        column name={$i$},
        string type
    },
    columns/err_cn/.style={
        column name={$e^{\mathcal{E}_\text{kin}}_{i}$},
        string type
    },
    columns/eoc_cn/.style={
        column name={EOC},
        string type
    },
    columns/err_ie/.style={
        column name={$e^{\mathcal{E}_\text{kin}}_{i}$},
        string type
    },
    columns/eoc_ie/.style={
        column name={EOC},
        string type
    },
    every head row/.style={
        before row={
            \toprule
            & \multicolumn{2}{c}{$\theta=0.5$}
            & \multicolumn{2}{c}{$\theta=1$} \\
            \cmidrule(lr){2-3}
            \cmidrule(lr){4-5}
        },
        after row=\midrule
    },
    every last row/.style={
        after row=\bottomrule
    }
]{\Ekin}
}
    \end{subtable}
    \hfill
    \begin{subtable}[t]{0.48\textwidth}
        \centering
        \caption{$\mathcal{D}_{\mathrm{kin}}$}
        \pgfplotstableread[
    col sep=tab,
    comment chars={\#}
]{simulationData/FCT_explicit/convergenceStudy/convergence_D_kin_ref=6_T=10s.txt}\Dkin

{
\setlength{\tabcolsep}{3pt}

\pgfplotstabletypeset[
    columns={i,err_cn,eoc_cn,err_ie,eoc_ie},
    columns/i/.style={
        column name={$i$},
        string type
    },
    columns/err_cn/.style={
        column name={$e^{\mathcal{D}_\text{kin}}_{i}$},
        string type
    },
    columns/eoc_cn/.style={
        column name={EOC},
        string type
    },
    columns/err_ie/.style={
        column name={$e^{\mathcal{D}_\text{kin}}_{i}$},
        string type
    },
    columns/eoc_ie/.style={
        column name={EOC},
        string type
    },
    every head row/.style={
        before row={
            \toprule
            & \multicolumn{2}{c}{$\theta=0.5$}
            & \multicolumn{2}{c}{$\theta=1$} \\
            \cmidrule(lr){2-3}
            \cmidrule(lr){4-5}
        },
        after row=\midrule
    },
    every last row/.style={
        after row=\bottomrule
    }
]{\Dkin}
}
    \end{subtable}

    \vspace{1em}

    \begin{subtable}[t]{0.48\textwidth}
        \centering
        \caption{$\mathcal{D}_{\mathrm{drag}}$}
        \pgfplotstableread[
    col sep=tab,
    comment chars={\#}
]{simulationData/FCT_explicit/convergenceStudy/convergence_D_drag_ref=6_T=10s.txt}\Ddrag

{
\setlength{\tabcolsep}{3pt}

\pgfplotstabletypeset[
    columns={i,err_cn,eoc_cn,err_ie,eoc_ie},
    columns/i/.style={
        column name={$i$},
        string type
    },
    columns/err_cn/.style={
        column name={$e^{\mathcal{D}_\text{drag}}_{i}$},
        string type
    },
    columns/eoc_cn/.style={
        column name={EOC},
        string type
    },
    columns/err_ie/.style={
        column name={$e^{\mathcal{D}_\text{drag}}_{i}$},
        string type
    },
    columns/eoc_ie/.style={
        column name={EOC},
        string type
    },
    every head row/.style={
        before row={
            \toprule
            & \multicolumn{2}{c}{$\theta=0.5$}
            & \multicolumn{2}{c}{$\theta=1$} \\
            \cmidrule(lr){2-3}
            \cmidrule(lr){4-5}
        },
        after row=\midrule
    },
    every last row/.style={
        after row=\bottomrule
    }
]{\Ddrag}
}
    \end{subtable}
    \hfill
    \begin{subtable}[t]{0.48\textwidth}
        \centering
        \caption{$\mathcal{R}_\theta$}
        \pgfplotstableread[
    col sep=tab,
    comment chars={\#}
]{simulationData/FCT_explicit/convergenceStudy/energy_balance_absolute_norm_ref=6_T=10s.txt}\BalanceNorm

{
\setlength{\tabcolsep}{3pt}

\pgfplotstabletypeset[
    columns={i,err_cn,eoc_cn,err_ie,eoc_ie},
    columns/i/.style={
        column name={$i$},
        string type
    },
    columns/err_cn/.style={
        column name={$e^{\mathcal{R}_{\theta}}_i$},
        string type
    },
    columns/eoc_cn/.style={
        column name={EOC},
        string type
    },
    columns/err_ie/.style={
        column name={$e^{\mathcal{R}_{\theta,}}_i$},
        string type
    },
    columns/eoc_ie/.style={
        column name={EOC},
        string type
    },
    every head row/.style={
        before row={
            \toprule
            & \multicolumn{2}{c}{$\theta=0.5$}
            & \multicolumn{2}{c}{$\theta=1$} \\
            \cmidrule(lr){2-3}
            \cmidrule(lr){4-5}
        },
        after row=\midrule
    },
    every last row/.style={
        after row=\bottomrule
    }
]{\BalanceNorm}
}
    \end{subtable}
\end{table}
%\end{comment} 

\subsection{Investigation of the Relative Importance of Model Components}
\label{sec:relative_importance_study}
\begin{table}[htbp]
    \centering
    \caption{Modifications of model components used to determine their relative importance for the simulation outcome in Section \ref{sec:relative_importance_study}.}
    \renewcommand{\arraystretch}{1.5}
    \centering
    \begin{tabularx}{\linewidth}{llX}
        \hline
        \textbf{Index} & \textbf{Model component} & \textbf{Modification} \\
        \hline
        1
        & Variable density of the mixture
        & Replacement of the SPION-dependent density in
          Equation \eqref{eq:two-phase_flow_system_with_relations:c}
          by the constant density of the carrier fluid
          $\tilde{\rho}_\textrm{f}$ \\

        2
        & Variable viscosity of the mixture
        & Replacement of the SPION-dependent viscosity in
          Equation \eqref{eq:two-phase_flow_system_with_relations:c}
          by the constant viscosity of the carrier fluid
          $\tilde{\eta}_\textrm{f}$ \\

        3
        & Additional momentum transport
        & Neglect the additional inertia term in Equation \eqref{eq:two-phase_flow_system_with_relations:c} \\

        4
        & SPION-dependent magnetic body force on the mixture
        & Neglect the right-hand side in Equation \eqref{eq:two-phase_flow_system_with_relations:c} \\

        5
        & SPION-dependent magnetization
        & Neglect the influence of the SPION-dispersion on the magnetization in Equation \eqref{eq:two-phase_flow_system_with_relations:a}

    \end{tabularx}
    \label{tab:ablation_study}
\end{table}
\begin{figure}
    \centering
    \input{figures/code/heatmaps}
    \caption{Heatmaps for the relative effect of modifying model components as described in Table \ref{tab:ablation_study}, see Equation \eqref{eq:relative_average_effect}.}
    \label{fig:heatmaps}
\end{figure}
The magneto two-phase flow model includes several coupling terms describing the interactions between the magnetic field, the SPION dispersion, and the carrier fluid; see also Remark \ref{remark:couplings}. In this section, we systematically investigate the influence of these coupling mechanisms on the numerical solution by conducting a study in which individual model components are selectively activated or deactivated. The considered components and their corresponding modifications are summarized in Table \ref{tab:ablation_study}. We repeat the numerical experiment introduced in the preceding sections for all possible combinations of the five model components, resulting in a total of $2^5 = 32$ distinct model configurations. This full factorial design allows us to quantify the influence of each component while accounting for the different configurations of the remaining components. The complete model, in which all coupling mechanisms are retained, serves as the reference configuration. Each configuration can be characterized in terms of a binary vector of length five 
\begin{equation}
    \boldsymbol{\alpha} = (\alpha_1, \alpha_2, \alpha_3, \alpha_4, \alpha_5) \in \{0,1\}^5 =: \mathcal{A},
\end{equation}
where $\alpha_i=1$ indicates that mechanism $i$ is included in its original from for the simulation while $\alpha_i=0$ means it is modified according to Table $\ref{tab:ablation_study}$. Consequently, $\boldsymbol{\alpha} = (1,1,1,1,1)$ corresponds to the full model with all couplings being active. For each of the $32$ configurations, we evaluate the scalar observables
\begin{equation}
    m_{\mathrm{target}}(10\textrm{s}), \quad \ \rho_\mathrm{Fe,max}(10\mathrm{s}), \quad \max_{(t,x) \in I \times \Om} \|\bfv(t,x) \|.
\end{equation}
As a measure for the relative importance of a mechanism, we introduce the relative average effect of modifying model component $i$,
\begin{equation} \label{eq:relative_average_effect}
    \Delta_i q = \frac{100\%}{|q_{(1,1,1,1,1)}|} \left( \frac{1}{16} \sum_{\boldsymbol{\alpha} \in \mathcal{A}_{i,0}} [ q_{\boldsymbol{\alpha}} - q_{(\boldsymbol{\alpha} + \mathbf{e}_i)} ] \right), \qquad i \in \{ 1,2,3,4,5\},
\end{equation}
where $q \in \{ m_{\mathrm{target}}(10\textrm{s}), \;  \rho_\mathrm{Fe,max}(10\textrm{s}), \; \max_{(t,x) \in I \times \Om} \|\bfv(t,x) \| \}$, $\mathcal{A}_{i,0} = \{ \boldsymbol{\alpha} \in \mathcal{A}, \alpha_i=0 \}$ and $\mathbf{e}_i \in \R^5$ is the $i$-th unit vector. Here, $q_{\boldsymbol{\alpha}}$ means that $q$ is computed from the simulation with configuration $\boldsymbol{\alpha}$. A positive value of $\Delta_i q$ indicates that modifying model component $i$ increases the observable $q$ on average, whereas a negative value indicates a decrease. The magnitude $|\Delta_i q|$ quantifies the relative strength of this effect.
\par
The numerical study is run with $\rho_\mathrm{Fe,initial} = 1 \,\mathrm{mg/ml}$ and with $\rho_\mathrm{Fe,initial} = 10 \,\mathrm{mg/ml}$ and the relative averaged effects of the modifications are visualized using heatmaps in Figure \ref{fig:heatmaps}. It is evident that, for all evaluated physical readouts, the model's strongest response is associated with the omission of the magnetic volume force from the fluid equation in \eqref{eq:two-phase_flow_system_with_relations_non_dimensional:c}, for both tested initial SPION dispersions. This is consistent with the previously noted observation that the dynamics of the system during the time period simulated here are primarily driven by the excitation of the velocity field due to the SPIONs and their resulting convection. The reason the average maximum iron density $\rho_\mathrm{Fe}$ at $t=10\mathrm{s}$ increases is that the mixing of the phases due to the convection does not occur in the first place and, except for a very slow transport due to magnetophoretic drift and diffusion, the initial SPION distribution is essentially maintained even after $10\mathrm{s}$. Furthermore, it should be noted that the readouts are significantly less sensitive to the other modifications. When $\rho_\mathrm{Fe,initial} = 1 \,\mathrm{mg/ml}$, they are all below $1\%$; for the flow-associated terms, they are even well below that, confirming that the effects of differences in density, viscosity, and relative momentum flux are negligible at low nanoparticle concentrations. This statement remains valid for $\rho_\mathrm{Fe,initial} = 10 \,\mathrm{mg/ml}$, provided one is willing to accept errors of approximately $1\%$. Although the average deviations for mechanisms 1–3 increase in some cases by an order of magnitude or more, they remain well below $1\%$. The situation is different, however, when it comes to the coupling of the magnetization of the SPION dispersion. Here, the relative effect rises to over $5\%$ for the mass in the target region and the maximum iron density at $t = 10\mathrm{s}$. The maximum velocity also changes by an average of about $3\%$. In line with this trend, we expect the mechanisms examined here to become increasingly important as nanoparticle concentrations rise. This is particularly relevant in MDT applications where nanoparticles are specifically accumulated magnetically, for example, to maximize the amount of a drug in the vicinity of a tumor.
\par
For this section, we also wanted to include the linearization of the \textit{Langevin} function in the SPION dispersion magnetization
\begin{equation}
    \mathbf{m}(u_\mathrm{p},\bfh) = M_\textrm{s} u_\mathrm{p} \mathcal{L}(\xi|\bfh|) \frac{\bfh}{|\bfh|} \approx \frac{\xi}{3}M_\textrm{s} u_\mathrm{p} \bfh, \quad \xi|\bfh| \ll 1,
\end{equation}
as this term can also be found in the literature for mathematical models of Magnetic Drug Targeting; see also Figure \ref{fig:spion_sketches:c} for a graphical intuition. However, in our setting, this modification resulted in magnetic forces approximately three orders of magnitude larger than those obtained with the original nonlinear formulation. Consequently, the flow velocities induced by the magnetic body force were substantially higher, and the resulting dynamics occurred on a much shorter time scale, potentially requiring significantly smaller time steps for adequate temporal resolution. We therefore excluded this modification from our study, as the linearized formulation leads to qualitatively different dynamics that are no longer directly comparable to those of the original model.
\section{Conclusion and Outlook}
\label{sec:conclusion_outlook}
In this work, we have derived and numerically investigated a thermodynamically consistent two-phase flow model for the transport of superparamagnetic iron oxide nanoparticles (SPIONs) suspended in a carrier fluid under the influence of an external magnetic field. Starting from suitable energy and dissipation functionals, the governing equations were derived within the \textit{Onsager} framework, yielding a coupled system that satisfies a global energy-dissipation balance. For the  simulation of the derived system, we have developed a numerical scheme based in the finite element method and semi-implicit time-stepping involving a one-step $\theta$-method for the temporal discretization. The modular structure of the scheme allows the individual subproblems to be treated separately, facilitating its implementation and reducing the computational complexity of the coupled problem. To stabilize the convection-dominated transport equation, an explicit flux-corrected transport (FCT) method was employed, yielding nonnegative discrete SPION dispersion densities in the numerical experiments.
Numerical simulations were performed for the attraction of an initially localized SPION dispersion towards an external magnet in a closed domain, with an initially quiescent mixture. By evaluating a range of physical observables, we were able to identify and distinguish the mechanisms governing the evolution of the dispersion. During the initial phase ($0-5\mathrm{s}$), the dynamics are dominated by convective transport associated with fluid motion induced by the magnetic body force. At later times, the comparatively slow magnetophoretic drift becomes increasingly important and drives the continued accumulation of SPIONs near the magnet. In particular, the simulations suggest that the accumulation of SPIONs beyond their initial maximum concentration is primarily governed by magnetophoresis and occurs on a considerably longer time scale than the initial convective redistribution. Furthermore, we investigated to what extent the numerical scheme reproduces the thermodynamic energy-dissipation balance of the continuous model. Although the balance is not satisfied exactly at the discrete level, the numerical results show that the deviations decrease systematically under temporal refinement. The convergence study indicates approximately first-order convergence of the energy-balance residual to zero. Finally, an ablation study involving $32$ model configurations was conducted to asses the relative influence of five individual coupling mechanisms on selected physical observables. Here, the magnetic body force exerted by the SPIONs onto the carrier fluid was found to be by far the strongest influence on the simulation outcomes. Moreover, the pronounced overestimation of magnetic forces resulting from a linearized \textit{Langevin} magnetization law highlights the importance of accounting for nonlinear saturation under the considered magnetic field strengths.
Coupling the SPION magnetization to the calculated magnetic field leads to comparably small changes in our observables, even though the influence of the SPIONs can, for example, alter the magnitude of the magnetic field gradient locally by 15\% when the initial SPION dispersion iron density was as low as $1\, \mathrm{mg/ml}$. However, we have seen that all investigated coupling mechanisms, and in particular the coupling to the magnetic field, become more important for higher SPION concentrations and we expect their influence to become more important as the SPION concentrations rise, e.g. due to magnet-induced accumulation in the vicinity of a tumor.
\par
The proposed model and its numerical approximation provide several directions for future research. From a mathematical perspective, a rigorous analysis of the coupled PDE system, including the existence and uniqueness of solutions, constitutes an important next step. In this context, the global energy-dissipation balance provides a natural source of a priori estimates. At the same time, analyzing and further developing the numerical scheme, particularly with a view toward energy-stable numerical methods, is a possible next step. One challenge here is that, from a physical standpoint, the non-negativity of the SPION-dispersion volume density must be ensured.
From an application-oriented perspective, the next step is to compare numerical predictions with experimental measurements from magnetic drug targeting experiments. Such investigations require simulations in realistic three-dimensional geometries, where the computational cost of the coupled problem becomes particularly significant. The modular and computationally efficient structure of the present numerical scheme provides a suitable starting point for these more demanding simulations. A detailed numerical investigation of three-dimensional SPION injection and targeting experiments will be addressed in future work. In this context, further physical mechanisms may need to be incorporated to describe SPION transport under realistic in vitro and in vivo conditions. These include particle adhesion to vessel walls, transport across the vessel wall into surrounding tissue, particle aggregation and agglomeration, pulsatile blood flow, and interactions with red blood cells.
\paragraph{Acknowledgements.} {\small This research is supported by DyNano (Dynamics and Control of Superparamagnetic Iron Oxide Nanoparticles in Simple and Branched Vessels). DyNano is a collaborative project funded by the Deutsche Forschungsgemeinschaft (DFG, German Research Fundation) - Project-ID 518492286 - and provided in particular the funding for the research of the second author.
The authors thank Prof. Dr. Carsten Gräser (FAU Erlangen-Nürnberg) for valuable support in the implementation of the model using \textit{DUNE}.}
\bibliographystyle{plain}
\bibliography{literature}
\end{document}